\documentclass[11pt]{article}
\usepackage{amsfonts}
\usepackage{amscd}
\usepackage{amsmath}
\usepackage{epsfig}

\newlength{\dinwidth}
\newlength{\dinmargin}
\newtheorem{definition}{Definition}
\newtheorem{theorem}{Theorem}
\newtheorem{proposition}{Proposition}
\newtheorem{corollary}{Corollary}
\newtheorem{remark}{Remark}
\newtheorem{lemma}{Lemma}

\def\cy{s}
\def\ko{\zeta}

\def\be{\begin{equation}}
\def\ee{\end{equation}}
\def\ben{\begin{displaymath}}
\def\een{\end{displaymath}}
\def\baa{\begin{eqnarray}}
\def\eaa{\end{eqnarray}}

\def\ba{\begin{array}}
\def\ea{\end{array}}
\makeatletter
\@addtoreset{equation}{section}
\makeatother

\def\x{{\bf x}}
\def\y{{\bf y}}
\def\z{{\bf z}}

\def\CP1{{\bf CP}^1}

\def\Z{{\bf Z}}

\def\a{\alpha}
\def\g{\gamma}

\def\b{\beta}

\def\l{\lambda}

\def\s{\sigma}

\def\O{\Omega}

\def\Oc{{\cal O}}
\def\P{\wp}
\def\met{{\bf m}}
\def\Fcal{{\cal F}}

\def\f{\frac}
\def\la{\label}
\def\H{{\cal H}}
\def\L{{\cal L}}
\def\Acal{{\cal A}}
\def\p{\partial}
\def\B{{\bf B}}

\def\Be{{\bf w}}
\def\res{{\rm res}}
\def\s{{\bf s}}
\def\c{{\cal C}}

\def\xt{{\tilde{x}}}
\def\Lhat{{\widehat{{\cal L}}}}
\def\rb{{\bf r}}
\def\P{P}
\def\Pt{{\tilde{P}}}
\def\wt{\tilde{w}}

\def\met{{\me}}

\def\cdiff{{\cal C}}
\def\me{{\bf g}}

\def\Abel{{\cal A}}
\def\hd{v}
\def\cdiff{{\cal C}}
\def\Wcal{{\cal W}}
\def\Si{{\bf s}}
\begin{document}
\title{
Tau-functions on spaces of Abelian  differentials and
higher genus generalizations of Ray-Singer formula}

\author{A. Kokotov\footnote{e-mail: alexey@mathstat.concordia.ca}
 \;  and
D. Korotkin
\footnote{e-mail: korotkin@mathstat.concordia.ca}}

\maketitle

\begin{center}
Department of Mathematics and Statistics, Concordia University\\
7141 Sherbrook West, Montreal H4B 1R6, Quebec,  Canada
\end{center}
\vskip0.5cm {\bf Abstract.} Let $w$ be an Abelian differential on compact
Riemann surface of genus $g\geq 1$. We obtain an explicit holomorphic
factorization formula for $\zeta$-regularized determinant of the Laplacian in
flat conical metrics with trivial holonomy $|w|^2$, generalizing the
classical Ray-Singer result in $g=1$.

\vskip0.5cm

\tableofcontents

\section{Introduction}

The goal of this paper is to give a natural generalization of the Ray-Singer formula for analytic torsion
of flat elliptic curves (\cite{Ray}) to the case of higher genus.

Let $A$ and $B$ be two complex numbers such that $\Im\left(B/A\right)>0$.
Taking the quotient of the complex plane ${\mathbb C}$ by the lattice
generated by $A$ and $B$, we obtain an elliptic curve (a Riemann surface of
genus one) $\L$. Moreover, the holomorphic one-differential $dz$ on ${\mathbb
C}$ gives rise to an Abelian differential $w$ on $\L$, so we get a  pair
(Riemann surface of genus one, Abelian differential on this surface) and the
numbers $A, B$ provide the natural local coordinates on the space of such
pairs. In what follows we refer to the numbers $A, B$ as moduli.

The modulus square $|w|^2$ of the Abelian differential $w$ generates a smooth flat metric on $\L$. Define the determinant of the laplacian $\Delta^{|w|^2}$ corresponding to this metric via the standard $\zeta$-function regularization:
\begin{equation}\label{zetareg}{\rm det}\Delta^{|w|^2}=\exp\{-\zeta_{\Delta^{|w|^2}}'(0)\},\end{equation}
where   $\zeta_{\Delta^{|w|^2}}(s)$ is the operator zeta-function.
Now (a slight reformulation of) the Ray-Singer
 theorem \cite{Ray} claims that  there holds the equality:
\begin{equation}\label{Ray}
\frac{{\rm det}\Delta^{|w|^2}}{\Im(B/A){\rm Area}(\L, |w|^2)}=C |\eta(B/A)|^4,
\end{equation}
where ${\rm Area}(\L, |w|^2)=\Im(A\bar{B})$, $C$ is a moduli-independent constant (actually, $C=4$) and $\eta$ is the Dedekind eta-function
$$\eta(\sigma)=\exp\left(\frac{\pi i \sigma}{12}\right)\prod_{n\in {\mathbb N}}\big(1-\exp(2\pi i n\sigma)\big).$$

The main result of this paper is a  generalization of formula (\ref{Ray}) to the case of
 Riemann surfaces of genus $g>1$.
To explain our strategy we first reformulate the Ray-Singer theorem.

For any compact Riemann surface $\L$ we introduce the prime-form $E(P,Q)$ (for definition and properties of this object we refer the reader to Sect.\ref{vfsad})
 and the canonical meromorphic bidifferential
\begin{equation}\label{bergdef}
\Be(P,Q)= d_P d_Q\log E(P,Q)\end{equation}
(see \cite{Fay92}).
 The bidifferential $\Be(P,Q)$ has
the following
local behavior as $P\to Q$:
\be
\Be(P,Q)= \left(\f{1}{(x(P)-x(Q))^2}+\f{1}{6} S_B(x(P))+
o(1)\right)dx(P)dx(Q),
\la{defproco}
\ee
where $x(P)$ is a local parameter. The term
$S_B(x(P))$ is a projective connection which is called {\it the Bergman
projective connection}.
Let $w$ be an Abelian differential on $\L$ and, as before, let
 $x(P)$ be some local parameter on $\L$.
Denote by $S_w(x(P))$ the Schwarzian derivative $\left\{\int^P
w,\,x(P)\right\}$.
Then the difference  of two projective connections $S_B-S_w$ is a
 (meromorphic) quadratic differential on $\L$ \cite{Tyurin}. Therefore, the ratio $(S_B-S_w)/w$
is a (meromorphic) one-differential.
In the elliptic case, i. e. when the Riemann surface $\L$ and the Abelian differential $w$ are  obtained from the lattice $\{mA+nB\}$,  this one-differential is holomorphic and admits the following explicit expression in the local parameter $z$ (see \cite{Fay73}):
\begin{equation}\label{ellfay}
\frac{S_B-S_w}{w}=-24\pi i \frac{d \log\eta(\sigma)}{d\sigma}\frac{1}{A^2}dz,
\end{equation}
where $\sigma=B/A$.

Let $\{a, b\}$ be the canonical basis of cycles on the elliptic curve $\L$
(the sides of the fundamental parallelogram), such that the numbers $A$ and
$B$ are the corresponding $a$ and $b$-periods of the Abelian differential
$w$. Defining
\begin{equation}\label{rododin}\tau(A, B):=\eta^2(B/A),\end{equation}
we see from (\ref{ellfay}) that the function $\tau$ is subject to the system of equations
\begin{equation}\label{elltau1}
\frac{\partial\log\tau}{\partial A}=\frac{1}{12\pi i}\oint_b\frac{S_B-S_w}{w}\,,
\hskip1.0cm
\frac{\partial\log\tau}{\partial B}=-\frac{1}{12\pi i}\oint_a\frac{S_B-S_w}{w}.
\end{equation}
Now the Ray-Singer formula implies that the real-valued expression
\begin{equation}\label{Qui}
Q(A, B)=\frac{{\rm det}\Delta^{|w|^2}}{\Im(B/A){\rm Area}(\L, |w|^2)}
\end{equation}
 satisfies the same system:
\begin{equation}\label{rs1}
\frac{\partial \log Q}{\partial A}=\frac{1}{12\pi i}\oint_b\frac{S_B-S_w}{w}\,,
\hskip1.0cm
\frac{\partial \log Q}{\partial B}=-\frac{1}{12\pi i}\oint_a\frac{S_B-S_w}{w}.
\end{equation}
Clearly, if $\tau(A, B)$ and $Q(A, B)$ are (respectively) a holomorphic and a real-valued solutions of system (\ref{elltau1}), then  $Q(A, B)=C|\tau(A, B)|^2$ with some constant factor
$C$.
Thus,  the Ray-Singer result can be reformulated as follows:
\begin{theorem}\label{nol}
\begin{enumerate}
\item
The system (\ref{elltau1}) is compatible and has a holomorphic solution $\tau$ .
This solution can be found explicitly and is given by (\ref{rododin}).
\item
The variational formulas (\ref{rs1}) for the determinant of the laplacian $\Delta^{|w|^2}$ hold.
\item
The expression (\ref{Qui})
can be represented as the modulus square of a holomorphic function of
moduli $A, B$; this function coincides with the function $\tau$ up to a moduli-independent factor.
\end{enumerate}
\end{theorem}

In what follows we call the function $\tau$ (a {\it holomorphic}
solution to system (\ref{elltau1}) {\it the Bergman
tau-function}, due to its close link with the Bergman projective
connection.

Generalizing the statement 1 of Theorem \ref{nol} to higher genus, we define and explicitly compute the Bergman tau-function on different strata of the spaces ${\cal H}_g$
of Abelian differentials over Riemann
surfaces i.e. the spaces of pairs $(\L,w)$, where $\L$ is a compact
Riemann
surface of genus $g\geq 1$ and $w$ is a holomorphic Abelian
differential
(i.e.  a holomorphic 1-form) on $\L$.
In global terms, the ``tau-function'' is not a function, but a
section of a line bundle over the covering of a stratum of ${\cal
H}_g$ (the space of triples $(\L,w,\{a_\a,\,b_\a\})$, where $w$ has
fixed multiplicities of its zeros; $(a_\a,\,b_\a)$ is a canonical
basis of cycles).

An analog of the Bergman tau-function on spaces of holomorphic
differentials
was previously defined on  Hurwitz spaces (see \cite{MPAG, IMRN}), i.e. on the  spaces
of pairs $(\L, f)$, where $f$ is a meromorphic function on a compact
Riemann surface $\L$ with fixed multiplicities  of poles and zeros of
the differential $df$. In this case it coincides with the
isomonodromic Jimbo-Miwa tau-function for a class of Riemann-Hilbert
problems \cite{Annalen,Dub92}, this explains why we use the term ``tau-function'' in the
context of spaces ${\cal H}_g$.

Generalizing statement 2 of Theorem \ref{nol},
 we introduce the laplacian, $\Delta^{|w|^2}$ acting in the trivial
 line bundle over $\L$, corresponding to the  flat
singular metric $|w|^2$. Among other flat metrics with conical
 singularities metrics of this form are distinguished by the property
 that they have {\it trivial holonomy} along any closed loop on the
 Riemann surface.

 Since Abelian differentials on Riemann surfaces of genus $g>1$ do have zeros,
the metric $|w|^2$  has  conical
singularities and the laplacian is not essentially self-adjoint. Thus, one
has to choose a proper self-adjoint extension:  here we deal with the
Friedrichs   extension. It turns out that it is still possible to define the determinant of this laplacian via the regularization (\ref{zetareg}).
We derive formulas  for variations of $\det\Delta^{|w|^2}$ with respect to natural coordinates on the space of Abelian differentials. These formulas are direct analogs of system (\ref{rs1}).

Generalizing statement 3 of Theorem \ref{nol}, we get an explicit formula for the determinant of the laplacian $\Delta^{|w|^2}$:
\be
{\rm det}\Delta^{|w|^2}=C\,{\rm Area}(\L, |w|^2)\;\{{\rm det} \Im{\B}\}\,|\tau|^2\;,
\la{mainres0}\ee
where $\B$ is the matrix of $b$-periods of a  Riemann surface of genus $g$,
and the Bergman tau-function $\tau$ is expressed through theta-functions and prime-forms.
This formula can be considered as a natural generalization of the Ray-Singer formula to the higher genus case.

\begin{remark}\rm
The determinants of Laplacians in flat conical metrics first appeared in
works of string theorists (see, e. g., \cite{Knizhnik}). An attempt to
compute such determinants was made in \cite{Sonoda}. The idea was to make use
of Polyakov's formula \cite{Polyakov} for the ratio of determinants of the
Laplacians corresponding to two {\it smooth} conformally equivalent metrics.
If one of the metrics in Polyakov's formula has conical singularity, this
formula does not make sense, so  one has to choose some kind of
regularization of the arising divergent integral. This leads to an
alternative definition of the determinant of Laplacian in conical metrics:
one may simply take some smooth metric as a reference one and define the
determinant of laplacian in a conical metric  through properly regularized
Polyakov formula for the pair (the conical metric, the reference metric).
Such a way was chosen in
 \cite{Sonoda} (see also \cite{DH})
for metrics given by the modulus square of an Abelian differential (which is
exactly our case) and metrics given by the modulus square of a meromorphic
1-differential (in this case Laplacians have continuous spectrum and the
spectral theory definition of their determinants, if possible, must use
methods other than the Ray-Singer regularization). In \cite{Sonoda} the
smooth reference metric is chosen to be the Arakelov metric. Since the
determinant of Laplacian in Arakelov metric is known (it was found in
\cite{DS} and \cite{ABMNV}, see also \cite{Fay92}); such an approach leads to
a heuristic formula for $\det\Delta $ in a flat conical metric. This result
heavily depends on the choice of the regularization procedure. The naive
choice of the regularization leads to dependence of $\det\Delta $ in the
conical metric on the smooth reference metric which is obviously
unsatisfactory. More sophisticated (and used in \cite{Sonoda} and \cite{DH})
procedure of regularization eliminates the dependence on the reference metric
but provides an expression which behaves as a tensor with respect to local
coordinates at the zeros of the differential $w$ and, therefore, also can not
be considered as completely satisfactory. In any case it is unclear whether
this heuristic formula for $\det\Delta $ for conical metrics has something to
do with the determinant of Laplacian  defined via the spectrum of the
operator $\Delta$ in conical metrics.
\end{remark}

The paper is organized as follows. In Section 2 we derive variational
formulas of Rauch type on the spaces of Abelian differentials for basic
holomorphic differentials, matrix of $b$-periods, prime-form and other
relevant objects. In section 3 we  introduce and compute the Bergman
tau-function on the space of Abelian differentials over Riemann surfaces. In
Section 4 we give a survey of the spectral theory of the Laplacian on
surfaces with flat conical metrics (polyhedral surfaces) and derive
variational formulas for the determinants of Laplacians in such metrics. The
comparison of variational formulas for the tau-functions with variational
formulas for the determinant of Laplacian, together with explicit computation
of the tau-functions, leads to the explicit formulas for the
determinants.
We use our explicit formulas to derive the formulas of Polyakov type,
which show how determinant of laplacian depends on the choice of the
conical metric on a fixed Riemann surface.

\section{Variational formulas on spaces of Abelian  differentials
over Riemann surfaces}

\subsection{Coordinates on the spaces of Abelian differentials}
\la{Coospabe}

The space ${\cal H}_g$ of  holomorphic Abelian differentials over
Riemann surfaces of genus $g$ is
 the moduli space of pairs $(\L, w)$,
where $\L$ is a compact Riemann surface of genus $g>1$, and $w$ is a
holomorphic 1-differential on $\L$.
This space is stratified according to the multiplicities of zeros of
$w$.

 The corresponding strata may have several connected components.
The classification of these connected components is given in
 \cite{KonZor}. In particular, the stratum of
the space ${\cal H}_g$ having the
highest dimension (on this stratum  all the zeros of $w$ are simple) is
connected.

Denote by ${\cal H}_{g}(k_1,\dots, k_M)$ the stratum of ${\cal H}_g$, consisting of
differentials  $w$ which have $M$ zeros on $\L$ of
multiplicities $(k_1,\dots,k_M)$. Denote the zeros of $w$ by
$P_1,\dots,P_M$; then the divisor  of
differential $w$ is given by
$(w)=\sum_{m=1}^M k_m P_m$.
Let us  choose a  canonical basis $(a_{\a},b_{\a})$ in the homology
group $H_1(\L,{\bf Z})$.
Cutting the Riemann surface $\L$  along these cycles we
 get the fundamental polygon $\Lhat$ (the fundamental polygon is not
 simply-connected unless all basic cycles pass through one point). Inside of $\Lhat$
we choose $M-1$  paths
 $l_{m}$ which
connect the zero $P_1$ with other zeros $P_m$ of $w$, $m=2,\dots,M$.
The set of paths $a_{\a},b_{\a},l_m$ gives a basis in the relative
 homology group $H_1(\L;(w),\Z)$.
Then the local coordinates on ${\cal H}_{g}(k_1,\dots, k_M)$ can be
chosen as follows (\cite{KZ1}, p.5):
\be
A_\a:=\oint_{a_\a} w\;,\hskip0.6cm
B_\a:=\oint_{b_\a} w\;,\hskip0.6cm
z_{m} :=\int_{l_{m}} w\;,\hskip0.5cm
\a=1, \dots, g;\  m=2,\dots,M \;.
\la{coordint}
\ee
The area of the surface $\L$ in the metric $|w|^2$ can be expressed
 in terms of these coordinates as follows:
$${\rm Vol}(\L) =-\Im \sum_{\a=1}^g A_\a \bar{B_\a}\;.$$
If all zeros of $w$ are simple, we have $M=2g-2$; therefore, the
dimension of the highest stratum ${\cal H}_{g}(1,\dots, 1)$ equals
$4g-3$.

The Abelian integral $z(P)=\int_{P_1}^P w$ provides a local coordinate in a
neighborhood of any point $P\in \L$ except the zeros $P_1,\dots,P_M$. In a
neighborhood of $P_m$ the local coordinate can be chosen to be
$(z(P)-z_m)^{1/(k_m+1)}$. The latter local coordinate is often called the {
\it distinguished} local parameter.

The following construction  helps to visualize these coordinates in the case of the highest stratum $H_g(1, \dots, 1)$.

\begin{figure}[hbt]
\begin{center}
\epsfig{file=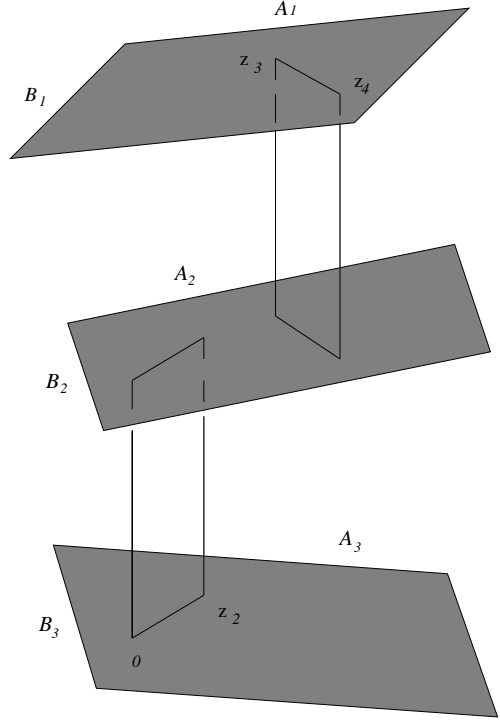}
\end{center}
\caption{Representation of a generic point of the stratum ${\cal
H}_3(1,1,1,1)$ by gluing three tori along cuts connecting zeros of
  $w$.}
\la{pseudocover}
\end{figure}

Consider  $g$ parallelograms $\Pi_1, \dots, \Pi_g$ in the complex plane with coordinate $z$ having the sides $(A_1, B_1)$, $\dots$,
$(A_g, B_g)$.  Provide these parallelograms with a system of cuts
$$[0, z_2],\ \ \  [z_3, z_4],\ \ \  \dots,\ \ \  [z_{2g-3}, z_{2g-2}]$$
(each cut should be repeated on two different parallelograms).
Identifying the opposite sides of the parallelograms and glueing the obtained $g$ tori along the cuts we get a compact Riemann surface $\L$ of genus $g$. (See figure 1 for the case $g=3$). Moreover, the differential $dz$ on the complex plane gives rise to a holomorphic differential $w$ on $\L$ which has
$2g-2$ zeros at the ends of the cuts. Thus, we get a point $(\L, w)$ from ${\cal H}_g(1, \dots, 1)$.
It can be shown that any generic point of ${\cal H}_g(1, \dots, 1)$ can be obtained via this construction; more sophisticated glueing is required to represent points of other strata, or non generic points of the stratum  ${\cal H}_g(1, \dots, 1)$.

The assertion about genericity follows from the theorem of  Masur and  Veech
(\cite{Masur}, \cite{Veech}, see also  \cite{KZ1})
stating  the ergodicity of the natural $SL (2, {\mathbb R})$-action on  connected components
of strata of the space of (normalized) Abelian differentials.
Namely, denote by ${\cal H}_g'(1, \dots, 1)$ the set of pairs $(\L, w)$ from
${\cal H}_g(1, \dots, 1)$ such that $\int_\L|w|^2=1$.
 Let a pair $(\L, w)$ from ${\cal H}_g'(1, \dots, 1)$  be obtained via the above construction.
Then  under the action of ${\bf A}\in SL (2, {\mathbb R})$ it goes to the pair $(\L_1, w_1)$ which is
obtained by gluing the parallelograms   ${\bf A}(\Pi_1), \dots, {\bf A}(\Pi_g)$
along the cuts $[0, {\bf A} z_2]$, $\dots$, $[{\bf A}z_{2g-3}, {\bf A}z_{2g-2}]$,
where the group $SL(2, {\mathbb R})$
acts on $z$-plane as follows
$${\bf A}:=\begin{pmatrix}
a & b\\
c & d
\end{pmatrix}:z\mapsto (a\Re z+b\Im z)+i(c\Re z+d\Im z)\,.$$
Thus, the set of pairs $(\L, w)$ from ${\cal H}_g'(1, \dots, 1)$ which can be glued from tori is invariant w. r. t. ergodic $SL(2, {\mathbb R})$-action, and, therefore, has the full measure.

To shorten the notations it is convenient to consider the coordinates
$\{A_\a,\,B_\a,\, z_m\}$ altogether. Namely,  in the sequel we shall
denote them by $\ko_k$, $k=1,\dots,2g+M-1$, where
\be\label{sokr}
\ko_\a:= A_\a\;,\hskip0.5cm \ko_{g+\a}:= B_\a\;,\hskip0.5cm \a=1,\dots, g\;,\hskip0.5cm
\ko_{2g+m}:= z_{m+1}\;\hskip0.5cm  m= 1,\dots, M-1
\ee

Let us also introduce corresponding cycles $\cy_k$,
$k=1,\dots,2g+M-1$, as follows:
\be
\cy_\a= - b_\a\;,\hskip0.5cm \cy_{g+\a}= a_\a\;,\hskip0.5cm \a=1,\dots, g\;;
\ee
the cycle $\cy_{2g+m}$, $m= 1,\dots, M-1$ is defined to be the small
circle with positive orientation around the point $P_{m+1}$.


Now we are going to  prove variational formulas (analogs of classical Rauch's formulas),
which describe dependence of basic holomorphic objects on Riemann surfaces
(the normalized
holomorphic differentials, the matrix of $b$-periods, the canonical
meromorphic bidifferential, the Bergman projective connection, the prime form, etc. ) on
coordinates (\ref{coordint}) on the spaces ${\cal H}_{g}(k_1,\dots,
k_M)$. We start from description of the objects we shall need in the sequel.

\subsection{Basic holomorphic objects on Riemann surfaces}

Denote by $\hd_\a(P)$ the basis of holomorphic 1-forms on $\L$ normalized
by
$\oint_{a_\a} \hd_\b=\delta_{\a\b}$. For a basepoint $P_0$ we define the Abel map
$\Acal_{\a}(P)=\int_{P_0}^P v_\a$ from the Riemann surface $\L$ to its Jacobian.

The matrix of b-periods of the surface $\L$ is given by
$\B_{\a\b}:=\oint_{b_\a} \hd_\b$.

Recall also the definition and properties of the  prime-form $E$, canonical meromorphic bidifferential $\Be$ and Bergman  projective connection $S_B$.

The prime form $E(P, Q)$ (see \cite{Fay73,Fay92})
 is an antisymmetric  $-1/2$-differential with respect to both $P$ and
$Q$. Let $\Theta [*]({\bf z})$ be the genus $g$ theta-function corresponding
to the matrix of $b$-periods ${\bf B}$ with some odd half-integer
characteristic $[*]$. Introduce the holomorphic differential
$q(P)=\sum_{\a=1}^g \Theta[*]_{z_\a}(0) v_{\a}(P)$. All zeros of this
differential are double and one can define the prime form on $\L$ by \be
E(P,Q)=\f{\Theta [*](\Acal(P)-\Acal(Q))}{\sqrt{q(P)}\sqrt{q(Q)}}\;;
\la{primeform} \ee this expression is independent of the choice of the odd
characteristic $[*]$.

The prime-form has the following properties (see \cite{Fay92}, p.4):
\begin{itemize}
\item Under tracing of $Q$ along the cycle $a_\a$ the prime-form
remains invariant; under
the tracing along  $b_\alpha$ it gains the factor
\begin{equation}\label{primetwist}
\exp(-\pi i \B_{\a\a}-2\pi i\int_P^Q v_\a)\;.
\end{equation}

\item On the diagonal $Q\to P$  the prime-form has first order zero (and no other zeros or poles)
with the following asymptotics:
$$E(x(P), x(Q))\sqrt{dx(P)}\sqrt{dx(Q)}=$$
\begin{equation}\label{primas}
(x(Q)-x(P))\left(1-\frac{1}{12}S_B(x(P))(x(Q)-x(P))^2+O((x(Q)-x(P))^3\right),
\end{equation}
where the subleading term $S_B$ is  called {\it Bergman projective connection} and $x(P)$ is an
arbitrary local parameter.
\end{itemize}

We recall that an arbitrary projective connection $S$ transforms under change of the local
coordinate $y\to x$ as follows:
\be
S(y)= S(x)\left(\f{d x}{d y}\right)^2 + \{x, y\}
\la{projcontrans}
\ee
where $\{x,y\}=\f{x'''}{x'}-\f{3}{2}\left(\f{x''}{x'}\right)^2$ is the Schwarzian derivative.
It is easy to verify that the term $S_B$ in (\ref{primas}) indeed transforms as
(\ref{projcontrans}) under change of the local coordinate. Difference of two projective
connections is a quadratic differential on $\L$.

The canonical meromorphic bidifferential $W(P,Q)$ is defined by $\Be(P,Q)=\p_P
\p_Q \log E(P,Q)$ (\ref{bergdef}). It is symmetric: $\Be(P,Q)=\Be(Q,P)$ and
has all vanishing $a$-periods with respect to both $P$ and $Q$; the only singularity
of $\Be(P,Q)$ is the second order pole on the diagonal $P=Q$ with biresidue
$1$. The subleading term in expansion of $\Be(P,Q)$ around diagonal is equal
to  $S_B/6$ (\ref{defproco}). The $b$-periods of $W(P,Q)$ with respect to any
of its arguments are given by the basic holomorphic differentials:
$\oint_{b_\a} W(P,\,\cdot\,)= 2\pi i \hd_\a(P)$.

 The prime-form can be expressed as follows in terms of $\Be(P,Q)$ (\cite{Fay92}, p. 3):
\begin{equation}\label{primrep}
E^2(P, Q)dx(P)dy(Q)=\lim_{P_0\to P, Q_0\to
Q}(x(P_0)-x(P))(y(Q)-y(Q_0))\exp\left(-\int_{P_0}^{Q_0}\int_P^Q \Be(\,
\cdot\,, \,\cdot\,)\right),
\end{equation}
where $x$ and $y$ are any local parameters near $P_0$ and $Q_0$,
respectively.
\begin{remark}\rm
Let us comment on the formula  (\ref{bergdef}) for $\Be(P,Q)$. Since
$E(P,Q)$ is a $-1/2$ differential with respect to $P$ and $Q$, this
formula should be understood as  $$\Be(P,Q)=\p_P \p_Q \{\log
E(P,Q)\sqrt{dx(P)}\sqrt{dy(Q)}\}\;,$$ where $x$ and $y$ are arbitrary local
parameters. Due to the presence of the operator $\p_P \p_Q$, this
expression is independent of the choice of these local parameters; therefore
it can be written in a shorter form (\ref{bergdef}), see \cite{Fay73,Mumford}.

In the same way we shall understand the formula for the normalized  (all
$a$-periods vanish) differential of the third king with poles at points $P$
and $Q$ and residues $1$ and $-1$, respectively (see \cite{Mumford}, vol. 2,
Chapter IIIb, Sect.1, p.212), which is extensively used below: \be
W_{P,Q}(R)=\p_R\log\f{E(R,P)}{E(R,Q)}\;. \la{Wintroduc} \ee This expression
should be rigorously understood as \be
W_{P,Q}(R)=\p_R\log\f{E(R,P)\sqrt{dx(P)}}{E(R,Q)\sqrt{dy(Q)}}\;, \la{Wdxdy}
\ee where $x$ and $y$ are arbitrary local coordinates; independence of
(\ref{Wdxdy}) of the choice of these local coordinates justifies writing it
in the short form (\ref{Wintroduc}).
\end{remark}

Denote by  $S_w(x(P))$ the projective connection given by the Schwarzian derivative $\left\{\int^P
w,\,x(P)\right\}$, where $x$ is a local parameter on $\L$.

The next object we shall need is the vector of Riemann constants:
\be
K^P_\alpha=\frac{1}{2}+\frac{1}{2}\B_{\a\a}-\sum_{\b=1,
\b\neq\a}^g\oint_{a_\b}\left(\hd_\beta\int_P^x\hd_\a\right)
\la{rimco}\ee
where the interior integral is taken along a path which does not
intersect
$\partial\widehat\L$.


Consider also the following multi-valued differential of  two
variables $\s(P,Q)$ ($P, Q\in\widehat \L$)
\begin{equation}\label{sigma}
\s(P, Q)=\exp\left\{
-\sum_{\a=1}^g\oint_{a_\a}\hd_\a(R)\log\frac{E(R, P)}{E(R, Q)}
\right\},
\end{equation}
where $E(R, P)$ is the prime-form (see \cite{Fay92}).
The right-hand side of (\ref{sigma}) is a non-vanishing holomorphic
$g/2$-differential on $\widehat\L$ with respect to $P$
and a non-vanishing holomorphic $(-g/2)$-differential with respect to
$Q$. Being lifted to the universal covering of
$\L$ it has  along the cycle $b_\a$  the automorphic factor
$\exp[(g-1)\pi i \B_{\a\a}+2\pi i K^P_\a]$ with respect to $P$ and the
automorphic factor $\exp[(1-g)\pi i \B_{\a\a}-2\pi i K_\a^Q]$ with respect to
$Q$.

 In what follows the pivotal role is played by the following holomorphic multivalued
$g(1-g)/2$-differential on $\Lhat$
\begin{equation}\label{c}
\cdiff(P)=\frac{1}{\Wcal[v_1, \dots, v_g](P)}\sum_{\a_1, \dots,
\a_g=1}^g
\frac{\partial^g\Theta(K^P)}{\p z_{\a_1}\dots \p z_{\a_g}}
v_{\a_1}\dots v_{\a_g}(P)\;,
\end{equation}
where
\be
\Wcal(P):= {\rm \det}_{1\leq \a, \b\leq g}||\hd_\b^{(\a-1)}(P)||
\la{Wronks}
\ee
is the Wronskian determinant of holomorphic differentials
 at the point $P$.

 It is easy to see that this differential
has multipliers
$1$ and $\exp\{-\pi i (g-1)^2\B_{\a\a}-2\pi i (g-1) K_\a^P\}$ along
basic
cycles $a_\a$ and $b_\a$, respectively.

The differential $\cdiff$ is an essential ingredient of the Mumford
measure
on the moduli space of Riemann surfaces of given genus
\cite{Fay92}. For $g>1$ the multiplicative differential $\Si$
(\ref{sigma}) is
expressed in terms of $\cdiff$ as follows \cite{Fay92}:
\be
\Si(P,Q)= \left(\f{\cdiff(P)}{\cdiff(Q)}\right)^{1/(1-g)}\;.
\la{Csi}
\ee
According to Corollary 1.4 from \cite{Fay92}, $\cdiff(P)$ does not have
any zeros.
Moreover, this object admits the following alternative representation:

\begin{equation}\label{repC}
\cdiff(P)=\frac{\Theta(\sum_{\a=1}^{g-1}\Acal_P(R_\a)+\Acal_Q(R_g)+K^P)\prod_{\a<\b}E(R_\a,
R_\b)\prod_{\a=1}^g\Si(R_\a, P)
}
{\prod_{\a=1}^g E(Q, R_\a)\,{\rm det}\,||v_\a(R_\b)||_{\a,
\b=1}^g\Si(Q, P)}\;,
\end{equation}
where $Q, R_1, \dots, R_g\in \L$ are arbitrary points of ${\cal L}$ and $\Acal_P$ is the Abel map with the base point $P$.

For arbitrary  points $P, Q, Q_0\in\L$ we introduce the following
multi-valued  $1$-differential
\begin{equation}\label{Om}
\Omega^{P}(Q)=\s^2(Q, Q_0)E(Q, P)^{2g-2}(w(Q_0))^g(w(P))^{g-1}
\end{equation}
(the $Q_0$-dependence of the right-hand side of (\ref{Om}) plays no
important role and is not indicated).

The differential $\Omega^{P}(Q)$ has automorphy factors $1$ and
$\exp(4\pi i K_\a^{P})$
along the basic cycles $a_\a$ and $b_\a$ respectively. The only zero
of the $1$-form $\Omega^{P}$ on
$\widehat\L$ is $P$; its multiplicity equals $2g-2$.
\begin{definition}
The projective connection $S_{Fay}^{P}$ on $\L$   given by the
Schwarzian derivative
\begin{equation}\label{Fc}
S_{Fay}^{P}(x(Q))=\left\{\int^Q\Omega^{P}, x(Q)\right\},
\end{equation}
where $x(Q)$ is a local coordinate on $\L$, is called the Fay
projective
connection (more precisely, we have here a family of projective
connections
parameterized by point $P\in\L$).
\end{definition}

Another projective connection we shall use below is associated to the
differential $w$ and given by the
Schwarzian derivative:
\be
S_w(x(Q)):=\left\{\int_{P_1}^Q w,\, x(Q)\right\}\;,
\la{defSw}
\ee
where $x(Q)$ is a local coordinate.

The difference of projective connections $S^P_{Fay}-S_w$ is a
quadratic differential.
The following lemma giving an expression for the 1-form $(S^P_{Fay}-S_w)/w$ is a simple corollary of above definitions:
\begin{lemma}\la{SFSW}
For any $Q\in\L$, $Q\neq P_m\;,\,\,m=1,\dots,M$
\begin{equation}
\f{1}{w}(S^P_{Fay}-S_w)(Q)=2\p_Q\left(\f{1}{w(Q)}\p_Q\log[\Si(Q,
    Q_0)E(Q,P)^{g-1}]\right)-
\f{2}{w(Q)}\left[\p_Q\log[\Si(Q, Q_0)E(Q, P)^{g-1}]\right]^2\;.
\la{SPSw}\end{equation}
\end{lemma}
{\it Proof.} We first notice that if one chooses the local parameter $x(Q)$
to coincide with $z(Q)$, then the projective connection $S_w$ vanishes:
$S_w(z(Q))=0$. Therefore, to find the left-hand side of (\ref{SPSw}) it is
sufficient to compute Fay's projective connection $S^P_{Fay}$ in the local
parameter $z(Q)$. From the definition (\ref{Fc}) of Fay's projective
connection and the definition (\ref{Om}) of multi-valued differential
$\Omega^{P}(Q)$ we get (\ref{SPSw}) taking into account that $d/d
z(Q)=w^{-1}(Q)\p_Q$.

$\square$

\begin{remark}\label{tensor}\rm In what follows we shall often treat tensor objects like $E(P, Q)$,
$\s(P, Q)$, etc as scalar functions of one of the arguments (or both). This
makes sense after fixing the local system of coordinates, which is usually
taken to be $z(Q)=\int^Q w$. In particular, the expression ``the value of the
tensor $T$ at the point $Q$ in local parameter $z(Q)$"  will mean the value
of the scalar $ T w^{-\alpha} $ at the point $Q$, where $\alpha$ is the
tensor weight of $T(Q)$. Very often one of the arguments (or sometimes both)
of the prime form coincide with a point  $P_m$ of the divisor $(w)$, in this
case we calculate the prime form in the corresponding distinguished local
parameter:
$$E(P, P_m):=E(P, Q)(dx_m(Q))^{1/2}|_{Q=P_m}\,.$$
\end{remark}

In the sequel we shall need the following theorem expressing the
differentials ${\bf s}(P,Q)$ and $\Omega^P(Q)$ in terms of
prime-forms. Since on Jacobian of the Riemann surface $\L$ the vectors
$\Abel_P((w))$ and $-2K^P$ coincide, there exist two vectors with
integer coefficients ${\bf r}$ and ${\bf q}$ such that
\be
\Abel_P\big((w)\big)+2K^{P}+{\bf B}{{\bf r}}+{\bf q}=0
\la{defvecpq}\ee
(here $(w):=\sum_{m=1}^M k_m P_m$ is the divisor of the differential $w$).

\begin{theorem}
The following expressions for  ${\bf s}(P,Q)$ and $\Omega^P(Q)$ hold:
\be
{\bf s}^2(P,Q)=\f{w(P)}{w(Q)}\prod_{m=1}^M\left\{\f{E(Q,P_m)}{E(P,P_m)}\right\}^{k_m}
e^{2\pi i\langle {\bf r},\,\Abel_P(Q)\rangle}
\la{sPQE}
\ee
\be
\Omega^P(Q)= E^{2g-2}(Q,P)\,w(Q) \{w(Q_0)w(P)\}^{g-1}\prod_{m=1}^M\left\{\f{E(Q_0,P_m)}{E(Q,P_m)}\right\}^{k_m}
e^{2\pi i\langle {\bf r},\,\Abel_P(Q)\rangle}
\la{OMPQPR}
\ee
\end{theorem}
{\it Proof.} We start from the following lemma:

\begin{lemma}\label{LL2}
The expression
\begin{equation}
\Fcal: = [w(P)]^{\f{g-1}{2}}e^{-{\pi i} \langle\rb,
K^P\rangle}
\left\{\prod_{m=1}^{M}[ E (P,P_m)]^{\f{(1-g)k_m}{2}}\right\} \c(P)
\la{const0}
\end{equation}
is independent of $P$.
\end{lemma}
{\it Proof.}
The tensor weight of $\Fcal$ with respect to $P$ is the sum of $(g-1)/2$ (from $w(P)$),
$\f{1-g}{2}\sum_{m=1}^M k_m$ (from the product of the prime-forms) and
$g(1-g)/2$ (from ${\cal C}(P)$), which equals $0$ since $\sum_{m=1}^M
k_m=2g-2$.  The zeros of $w(P)$ at $\{P_m\}$ are canceled against poles
arising from the
product of prime-forms.

Therefore, to prove that $\Fcal$ is constant with respect to $P$ it remains
to show that this expression does not have any monodromies along basic
cycles.  Because of  uncertainty of the sign choice if $(g-1)/2$ is
half-integer it is convenient to consider $\Fcal^2$. The only ingredient of
(\ref{const}) which changes under analytical continuation along the cycle
$a_\a$ is the vector of Riemann constants; the expression $\langle\rb,
K^P\rangle$ transforms to
 $\langle\rb, K^P\rangle+ (g-1)r_\a$, which, since $r_\a$ is an
integer, gives trivial monodromy of $\Fcal^2$ along $a_\a$.

Under analytical continuation along the cycle $b_\a$ the prime-form
$E(P,P_m)$ is multiplied with ${\rm exp} \{-\pi i \B_{\a\a}-2\pi i
(\Acal(P)-\Acal(P_m))\}$, and ${\cal C}(P)$ is multiplied with ${\rm exp} \{-\pi i
(g-1)^2\B_{\a\a}-2\pi i(g-1)K_\a^P\}$. Finally, the  expression $\langle\rb,
K^P\rangle$ transforms to
 $\langle\rb, K^P\rangle+ (g-1)(\B\rb)_\a$.

Collecting all these terms, we see that $\Fcal^2$ gets multiplied with
$$
{\rm exp}\big\{-2\pi i(g-1)[\Acal_\a\big( (w)\big)+2K_\a^P+ ({\bf B}\rb)_\a]\big\}
$$
which, due to (\ref{defvecr}), equals ${\rm exp}\{-2\pi i(g-1) q_\a\}=1$.

Therefore, $\Fcal^2$ is a holomorphic function on $\L$ with respect to $P$.
Hence, it is a constant, as well as $\Fcal$ itself.

$\Box$

Now the expression (\ref{sPQE}) follows from the link (\ref{Csi})
between ${\bf s}(P,Q)$ and ${\cal C}(P)$ and standard relation between
vectors of Riemann constants computed at different points:
$K^Q-K^P=(g-1)\Abel_P(Q)$.
The formula (\ref{OMPQPR}) follows from (\ref{sPQE})
and definition  (\ref{Om})  of $\Omega^P(Q)$.

$\Box$

\subsection{Variational formulas on spaces of holomorphic differentials}
\la{vfsad}

Variation of the coordinates $\{\zeta_k\}\equiv\{A_\a,\,B_\a,z_m\}$
generically changes the conformal structure of the Riemann surface $\L$. Here
we derive an analog of the Ahlfors-Rauch  formula for the variation of the
matrix of $b$-periods of $\L$ under variation of the coordinates
$\{\zeta_k\}$. Besides that, we find formulas for the variation of the
objects depending not only on the moduli of $\L$, but also on a point on $\L$
(as well as the choice of a local coordinate near this point), namely, the
basic holomorphic differentials $v_\a(P)$, the canonical bidifferential
$\Be(P,Q)$, the prime-form $E(P,Q)$, the differential $C(P)$ and other
objects described in the previous section.

We {\it define} the derivative of the basic holomorphic  differentials with
respect to $\zeta_k$ as follows:
\be \f{\p v_\a(P)}{\p\zeta_k}\Big|_{z(P)}:=
w(P) \f{\p}{\p\zeta_k}\Big|_{z(P)=const}\left\{\f{v_\a(P)}{w(P)}\right\}
\la{defderva}
\ee
where, as before, $z(P)=\int_{P_1}^P w$; $v_\a(P)/w(P)$ is
a meromorphic function on $\L$ with poles at $\{P_m\}$. Outside of the points
$P_m$ this function can  be viewed as a function of $z(P)$ and $\zeta_k$; the
derivative of this function with respect to $\zeta_k$ in the right-hand side
of (\ref{defderva}) is computed assuming that $z(P)$ is independent of
$\zeta_k$.

To introduce this definition in a more formal manner\footnote{We thank the
referee for
  mentioning this point.} consider the local universal family
$p:\;\;{\cal X}\to {\cal H}_g(k_1,\dots, k_M)$. Then the set
$(z:=\int_{P_1}^P w,\zeta_1,\dots,\zeta_{2g+M-1})$ gives a system of local
coordinates on ${\cal X}\setminus (w)$. A vicinity of a point $\{(\L, w),
P\}$ in the level set  $H_{z(P)}:= \{x\in {\cal X}\;,\; z(x)=z(P)\}$ is
biholomorphically mapped onto a vicinity of the point $(\L, w)$ of ${\cal
H}_g(k_1,\dots, k_M)$ via the projection $p\;:\; {\cal X}\to {\cal
H}_g(k_1,\dots, k_M)$. Then
$((p|_{H_{z(P)}})^{-1})^*\left\{\f{\hd_\a}{w}\right\}\big|_{H_{z(P)}}$ is a
  locally holomorphic function on ${\cal H}_g(k_1,\dots, k_M)$ and we
  denote
$$
 \f{\p}{\p\zeta_k}\Big|_{z(P)=const}\left\{\f{v_\a(P)}{w(P)}\right\}:=
 \f{\p}{\p\zeta_k}\left[((p|_{H_{z(P)}})^{-1})^*\left\{\f{\hd_\a}{w}\right\}\Big|_{H_{z(P)}}\right]\;.
$$
The differentiation with respect to $\{\zeta_k\}$ of other objects below (the bidifferential $W$,
the prime-form etc) will be understood in the same sense.

This differentiation looks very natural if $\L$ can be visualized
as a union of glued tori as in Figure \ref{pseudocover}. In this picture a
function $f(P)$ (depending also on moduli) on $\L$ is considered locally as a
function of $z$ and is differentiated with respect to $A_\a,B_\a$ and $z_m$
assuming that the projection $z(P)$ of the point $P$ on the $z$-plane remains
constant.

The derivatives
$\f{\p}{\p\zeta_k}\Big|_{z(P)}\left\{\f{v_\a(P)}{w(P)}\right\}$ are
meromorphic  in the fundamental polygon $\Lhat$, since  the map $P\mapsto
z(P)$ is globally defined in $\Lhat$; these derivatives are not necessarily
meromorphic functions globally defined on $\L$ since $z(P)$ is not
single-valued on $\L$. (Notice also that the map $P\mapsto z(P)$ is locally
univalent in $\L\setminus\{P_1, \dots, P_m\}$.)

The derivatives $\f{\p v_\a(P)}{\p\zeta_k}$ defined by (\ref{defderva}) are
therefore  meromorphic differentials of $(1,0)$ type defined within $\Lhat$;
they do not necessarily correspond to single-valued meromorphic differentials
on $\L$ itself.

Similarly, the derivatives of $\Be(P,Q)$ with respect to the moduli are
defined as follows: \be \f{\p \Be(P,Q)}{\p\zeta_k}\Big|_{z(P),z(Q)}:=
w(P)w(Q)\f{\p}{\p\zeta_k}\Big|_{z(P),z(Q)}\left\{\f{\Be(P,Q)}{w(P)w(Q)}\right\}
\la{defderbe} \ee Derivatives of other tensor objects depending not only on
moduli, but also on points of $\L$, are defined in the obvious analogy to
(\ref{defderva}) and (\ref{defderbe}).


\begin{remark}\rm
Our definition (\ref{defderva}) of the variation of $v_\a(P)$ with respect to
the coordinates on the space ${\cal H}(k_1,\dots,k_M)$ is different from the
variational scheme used by Fay (\cite{Fay92}, Chapter 3). In this scheme the
variation of $v_\a(P)$ in the direction defined by an arbitrary Beltrami
differential is computed assuming that the pre-image under the Fuchsian
uniformization map of the point $P$ on the upper half-plane (for $g\geq 2$)
is independent of the moduli. In this scheme the differential of the type
$(0,1)$ is present in the variational formula for $v_\a$, $\Be(P,Q)$ and
other objects (\cite{Fay92}, formula (3.21)). This $(0,1)$ contribution is
absent in our deformation framework by definition (\ref{defderva}),
(\ref{defderbe}). This difference makes it difficult to directly apply the
variational formulas for all interesting holomorphic objects which were
derived in \cite{Fay92} in our present context. However, many technical tools
of \cite{Fay92} can be used in our framework, too.

Actually, the deformation scheme we develop here is close to the Rauch
deformation of a branched covering via variation of a branch point
\cite{Rauch}. In particular, in the Rauch formulas for the basic holomorphic
differentials it is assumed that the projection of the argument of the
differential on the base of the covering is independent of the branch points.

\end{remark}

\begin{remark}{\rm
In what follows we very often deal with the derivatives with respect to
moduli of various integrals over a contour $\Gamma$  on the surface $\L$. In
this case calculations simplify under the assumption that the image of the
contour $\Gamma$ under the map $\L\ni P\mapsto z(P)=\int_{P_1}^Pw$ does not
vary under the variation of moduli. If the contour of integration coincides
with one of the cycles, say $a_1$, chosen to define the fundamental polygon
$\Lhat$, the map $P\mapsto z(P)$ and the local coordinates $\{A_\alpha,
B_\alpha, z_m\}$, then one can assume that the image of this contour does not
vary under variation of moduli $\{A_2, \dots, A_g, B_1, \dots, B_g, z_2,
\dots, z_m\}$ (and, of course, not $A_1$: in this case such an assumption is
no longer possible, in the sequel we shall consider expressions of the type
$\partial_{A_1}\oint_{a_1}$ in more detail).}
\end{remark}

\begin{theorem}\la{varfow}
The following variational formulas hold:
\be
\f{\p \hd_\a(P)}{\p \ko_k}\Big|_{z(P)} =
\f{1}{2\pi i}\oint_{\cy_k}\f{\hd_\a(Q)\Be(P,Q)}{w(Q)}\;,
\la{varw}
\ee
\be
\f{\p\B_{\a\b}}{\p \ko_k}=\oint_{\cy_k}\f{\hd_\a \hd_\b}{w}
\la{varB}
\ee
\be
\f{\p \Be(P,Q)}{\p \ko_k}\Big|_{z(P),\,z(Q)} =\frac{1}{2\pi
i}\oint_{\cy_k}\f{\Be(P,R)\Be(Q,R)}{w(R)}
\la{bevar}
\ee
\be\label{varprfo}
\f{\p }{\p \ko_k}\Big|_{z(P),\,z(Q)}\log \{E(P,Q)w^{1/2}(P)w^{1/2}(Q)\}=-\f{1}{4\pi i}
\oint_{\cy_k}\f{1}{w(R)}\left[\p_R\log\f{E(P,R)}{E(Q,R)}\right]^2
\ee
\be
\f{\p}{\p \ko_k}(S_B(P)-S_w(P))\Big|_{z(P)}=\f{3}{\pi
i}\oint_{\cy_k}\f{\Be^2(P,R)}{w(R)}\;,
\la{varcon}
\ee
where $k=1,\dots, 2g+M-1$;  we assume that the local coordinate
 $z(P)=\int_{P_1}^{P}w$
and  $z(Q)=\int_{P_1}^{Q}w$ are kept constant under
differentiation.
\end{theorem}
{\it Proof.} Let us prove first the variational formula  (\ref{varw}) for the
normalized holomorphic differential. As explained in Section \ref{Coospabe},
we use the Abelian integral $z(P)=\int_{P_1}^{P}w$ as a local coordinate in a
neighborhood of any point of $\L$ not coinciding with the zeros, $P_m$, of
the differential $w$. In a neighborhood of $P_m$ the local coordinate is
taken to be $x_m(P)=(z(P)-z_m)^{1/(k_m+1)}$, where $k_m$ is the multiplicity
of $P_m$. Consider now the derivative of $\hd_\a(P)$ with respect to $z_m$
($m\geq 2$) assuming that the coordinate $z(P)$ is independent of $z_m$. The
proof of the corresponding variational formula is completely parallel to the
proof of the standard Rauch formula on the Hurwitz spaces (see for example
Section 2.3 of \cite{MPAG}).

The differential $\p_{z_m}\hd_\a(P)|_{z(P)}$ is holomorphic outside of $P_m$
and has all vanishing $a$-periods (since the $a$-periods of $\hd_\a$ are
constant). Let us consider  the local behavior of
$\p_{z_m}\hd_\a(P)|_{z(P)}$ near $P_m$. We choose the local parameter near
$P_m$ to be $x_m=(z(P)-z_m)^{1/(k_m+1)}$. We have \be \hd_\a(x_m)= (C_0+C_1
x_m+\dots+C_{k_m}x_m^{k_m}+ O(|x_m|^{k_m+1})dx_m\;. \ee Differentiating this
expansion with respect to $z_m$ for fixed $z(P)$, we get:
$$
\f{\p}{\p {z_m}}\{\hd_\a(P)\}|_{z(P)}=
\left\{C_0\left(1-\f{1}{k_m+1}\right)\f{1}{x_m^{k_m+1}}+
C_1\left(1-\f{2}{k_m+1}\right)\f{1}{x_m^{k_m}}\right.
$$
\be \left.+\dots+
C_{k_m-1}\left(1-\f{k_m}{k_m+1}\right)\f{1}{x_m^2}+O(1)\right\}dx_m
\la{asderPm} \ee

Consider the set of standard meromorphic differentials of second kind with
vanishing $a$-periods: $W_{P_m}^{s+1}(P)$ with the only singularity
at the point $P_m$ of the form $x_m(P)^{-s-1}dx_m(P)$.
Since the differential (\ref{asderPm}) also has all vanishing
$a$-periods, it can be expressed in terms of these standard
differentials as follows:
$$
\f{\p}{\p {z_m}}\{\hd_\a(P)\}|_{z(P)}= C_0\left(1-\f{1}{k_m+1}\right)
W_{P_m}^{k_m+1}(P)+
C_1\left(1-\f{2}{k_m+1}\right) W_{P_m}^{k_m}(P)+\dots
$$
\be +C_{k_m-1}\left(1-\f{k_m}{k_m+1}\right)W_{P_m}^{2}(P)
\la{derhdOmega} \ee
Now, the differentials  $W_{P_m}^{s}(P)$  can be
expressed in terms of $\Be(P,Q)$ as follows:
\be W_{P_m}^{s}(P)=
\f{(-1)^{s}}{(s-1)!}\f{d^{s-2}}{d x_m^{s-2}(Q)}\Be(P,Q)\Big|_{Q=P_m}
\la{OmegaW} \ee

Using  (\ref{OmegaW}) we can rewrite (\ref{derhdOmega}) in the
following compact form:
$$
\f{\p \hd_\a(P)}{\p z_m}\Big|_{z(P)}
=\f{1}{(k_m+1)(k_m-1)!}\left(\f{d}{d x_m(Q)}\right)^{k_m-1}\left\{\f{\Be(P,Q)\hd_\a(Q)}{(dx_m(Q))^2}\right\}\Big|_{Q=P_m}
$$
or, equivalently,
 \be \f{\p \hd_\a(P)}{\p z_m}\Big|_{z(P)} = {\rm
res}\Big|_{Q=P_m}\f{\hd_\a(Q)\Be(P,Q)}{w(Q)}\;, \la{varwproof} \ee
 which leads to (\ref{varw}) for $k=2g+1,\dots,2g+M-1$.

Let us now prove formulas  (\ref{varw}) for $k=1,\dots 2g$. For example,
consider the derivative of $v_\a$ with respect to $B_\b$.

Denote by $U$ the universal covering of $\L$; let us choose the
fundamental cell (the ``fundamental polygon'' of $\L$) $\Lhat$ such that
 all the contours $l_m$ from the definition
(\ref{coordint}) of coordinates $z_m$ lie inside of $\Lhat$. The map $z(P)$
is a holomorphic function on $\Lhat$ with critical points at $\{P_m\}$.
Consider an arbitrary point in $\Lhat$ which does not coincide with any zero
of $w$; consider a neighborhood  $D\subset \Lhat$ of this point where $z(P)$
is univalent; denote by $\tilde{D}$ the image of $D$ under mapping $z(P)$:
$\tilde{D}=z[D]$.

Denote by  $T_{b_\b}$ the deck
transformation on $U$ which
 corresponds to the side $b_\b^+$  of the fundamental polygon.
 Consider the domain
 $T_{b_\b}[D]$ lying in the fundamental cell $T_{b_\b}[\Lhat]$ as well as
 its image in the $z$-plane $\tilde{D}_{b_\b}=\{z+B_\b|\; z\in \tilde{D}\} $.
We can always take sufficiently small domain $D$ such that
$\tilde{D}_{b_\b}\cap  \tilde{D}=\emptyset$. The holomorphic differential
$\hd_\a$ can be lifted from $\L$ to a holomorphic differential on $U$
invariant with respect to the deck transformations. Let us write
$v_\a(P)=f(z)dz$ for  $P\in D\cup D_{b_\b}$, $z:=z(P)\in
\tilde{D}\cup\tilde{D}_{b_\b} $. Since  $v_\a$ is invariant under the deck
  transformations, we have
\be f(z+B_\b)=f(z)\;,\hskip0.7cm z\in \tilde{D}\;.
\la{fzBfz} \ee
 Assuming $z$ to be
constant and differentiating this equality with respect to $B_\b$, and
taking into account that $ \f{\p f}{\p z}(z+B_\b)=\f{\p f}{\p z}(z)$ as a corollary of
(\ref{fzBfz}), we get:
\be \f{\p f}{\p B_{\b}}(z+B_\b) = \f{\p f}{\p B_{\b}}(z)- \f{\p f}{\p
z}(z)\;,\hskip0.7cm z\in \tilde{D}\;.
\la{perf} \ee
Let us denote
\be
\Phi(P):=\f{\p \hd_\a(P)}{\p B_\b}\Big|_{z(P)}\;,\hskip0.7cm P\in U\;;
\la{notderva} \ee
Since the coordinate $z(P)$ is single-valued on the
universal covering $U$,  the differential $\Phi$ is also single-valued and
holomorphic on $U$. Consider also the quadratic differential $\p[ \hd_a(P)]$,
which in a local coordinate $x$ is given by $(\p f/\p x) (dx)^2$ if
$\hd_\a(P)=f(x)dx$. Now we can rewrite (\ref{perf})  in a
coordinate-independent form:
\be T_{b_\b}[\Phi(P)]=\Phi(P)
-\f{\p[ \hd_\a(P)]}{w}\;,\hskip0.7cm P\in D\;.
\la{Tbhd} \ee
In complete analogy to
(\ref{Tbhd}) we can show that
\be T_{b_\g}[\Phi(P)]=\Phi(P)\;,\hskip0.7cm
\g\neq\b\;,\hskip0.7cm P\in D \la{Tbhd1} \ee
and
\be
T_{a_\g}[\Phi(P)]=\Phi(P)\;,\hskip0.7cm \g=1,\dots,g\;,\hskip0.7cm P\in D
\la{Tbhd2}
\ee
Since the formulas (\ref{Tbhd}), (\ref{Tbhd1}), (\ref{Tbhd2})
are valid in a neighborhood of any point of $\Lhat$ except  $\{P_m\}$, and
the differential $\Phi$ is holomorphic in $\Lhat$, we conclude that these
formulas are valid for any $P\in\Lhat$. Therefore, the differential $\Phi$
can be viewed as  a differential on $\L$  itself, which is holomorphic
everywhere except the cycle $a_\b$, where it has the additive jump  given by
$-{\p[ \hd_a(P)]}/{w}$. Moreover, it has all vanishing $a$-periods (this
condition of vanishing of all the $a$-periods obviously makes sense, since
all the $a$-periods of the ``jump differential" $-{\p[ \hd_a(P)]}/{w}$ also
vanish).

To write down an explicit formula for $\Phi$ we recall that on the complex
plane the contour integral $(1/2\pi i)\oint_C f(x)(x-y)^{-2}dx$ taken in
positive direction defines the functions $f^l$ and $f^r$ which are
holomorphic in the interior and the exterior of $C$, respectively, and  on
$C$  the boundary values of $f^r$ and $f^l$ (indeces $l(eft)$ and $r(right)$
refer to the side of the oriented contour $C$, where the boundary value is
computed)
 are
related by the Plemelj formula $f^r(y)-f^l(y)= -f_y(y)$.

This observation allows to write immediately the formula for the differential
$\Phi$ with discontinuity $-{\p[ \hd_a(P)]}/{w}$ on the cycle $a_\b$ and all
vanishing $a$-periods:
\be \Phi(P)=\f{1}{2\pi
i}\oint_{a_\b}\f{\hd_\a(Q)\Be(P,Q)}{w(Q)}\;;
\la{Phihdw}
\ee
the required discontinuity on the cycle $a_\b$ is implied by
singularity structure of $\Be(P,Q)$ and
Plemelj formula; vanishing of all the $a$-periods follows from  vanishing of
all the $a$-periods of bidifferential $\Be(P,Q)$.  Formula (\ref{Phihdw})
implies (\ref{varw})  for $k=M+g,\dots,M+2g-1$.

The formula for differentiation with respect to $A_\b$ has the different sign
due to the interchange of ``left" and ``right" in that case (due to the
asymmetry between the cycles $a_\b$ and $b_\b$ imposed by their intersection
index $a_\b\circ b_\b=-b_\b\circ a_\b=1$).

Integrating (\ref{varw}) over   $b$-cycles and changing the order of integration,  one gets (\ref{varB}). Formula (\ref{bevar}) can be proved in the same manner as (\ref{varw}).
Formula (\ref{varcon})  follows
from the variational formulas for the bidifferential $\Be(P,Q)$
(\ref{bevar})
in the limit $P\to Q$ if we write down these formulas with respect to
the local coordinate $z(P)$ (in this local coordinate the projective
connection $S_w$ vanishes) and take into account the definition
(\ref{defproco}) of the Bergman projective connection.

 The variational formula for the prime-form (\ref{varprfo})  follows from
the variational formula for $W(P,Q)$ (\ref{bevar}) and  the formula
$W(P,Q)=d_P d_Q\log \{E(P,Q)w^{1/2}(P) w^{1/2}(Q)\}$ defining $W(P,Q)$ in
terms of the prime-form. Namely, applying the second derivative   $d_Pd_Q$ to
(\ref{varprfo}) we arrive at (\ref{bevar}) (after squaring the integrand of
(\ref{varprfo}) and taking into account that the functions depending on $P$
or $Q$ only are annihilated by $d_Pd_Q$). Since (\ref{bevar}) is valid, we
see that (\ref{varprfo}) holds up to addition of a function of the form
$f(P)+g(Q)$, where $f(P)$ and $g(Q)$ are two functions holomorphic in
$\Lhat$. Since both left- and right-hand sides of (\ref{varprfo}) vanish at
$P=Q$, we have $g(Q)=-f(Q)$ and the additional term is of the form
$f(P)-f(Q)$. Furthermore, one can verify that the function $f(P)$ is
single-valued on $\L$. Namely, the left- and right-hand sides of
(\ref{varprfo}) have trivial monodromy along any $a$-cycle. Under analytical
continuation of variable $P$  along a cycle $b_{\a}$ the left-hand side of
(\ref{varprfo}) gains  due to (\ref{primetwist}) an additive term
$\p_{\zeta_k}\{-\pi i \B_{\a\a}- 2\pi i (U_{\a}(P)-U_{\a}(Q))\}$. By making
use of variational formulas (\ref{varw}), (\ref{varB}) it is easy to verify
that this term coincides with the additive term arising (due to
transformation law (\ref{primetwist})) in the right-hand side of
(\ref{varprfo}) under analytical continuation along $b_\a$ with respect to
variable $P$.

Therefore, the function $f(P)$ is a holomorphic single-valued function of
$P$; thus $f(P)=const$ and $f(P)-f(Q)=0$; therefore, the formula
(\ref{varprfo}) holds without any additional constants.

$\square$

In the sequel we shall also need to differentiate the prime-form $E(P,P_m)$
with respect to coordinate $z_m$ (this case is not covered by the variational
formula (\ref{varprfo}) since $z(P_m):=z_m$ can not be kept constant under
differentiation). Surprisingly enough, such formula still looks the same as
(\ref{varprfo}):
\begin{corollary}\la{varprfozm}
The following variational formula holds for any $m=2, \dots, M$:
$$
\f{\p \log \{E(P,P_m)w^{1/2}(P)\}}{\p z_m}\Big|_{z(P)}=-\f{1}{4\pi i}
\oint_{s_{2g+m-1}}\f{1}{w(R)}\left[\p_R\log\f{E(P,R)}{E(P_m,R)}\right]^2
$$
 \be\label{varprfo1}
\equiv -\f{1}{2}{\rm
res}\Big|_{R=P_m}\left\{\f{1}{w(R)}\left[\p_R\log\f{E(P,R)}{E(P_m,R)}\right]^2\right\};\
\ee
$$E(P,P_m):=E(P,Q)(d x_m(Q))^{1/2}\Big|_{Q=P_m}\;,$$
as before, $x_m(Q)=(z(Q)-z_m)^{1/(k_m+1)}\equiv\left(\int_{P_m}^Qw\right)^{1/(k_m+1)}$.
\end{corollary}

{\it Proof.}
 In what follows we shall use the simplified notation $r:=k_m+1$ and $C:=s_{2g+m-1}$.
 Let $Q$ be a point in a vicinity of $P_m$
whose $z$-coordinate is kept fixed, for $x_m$ coordinate of this point we
 shall use the simplified notation $x_m:=x_m(Q)$.
One has  $z(Q)-z_m=x_m^r$ and $\frac{\partial}{\partial
z_m}x_m=-\frac{1}{rx_m^{r-1}}$.

 Calculating $E(P, Q)$ in the local parameter $z(Q)$ and in the local parameter
$x_m$, one gets
$$E(P, Q)w(Q)^{1/2}=\left(E(P, Q)\sqrt{d x_m}\right)
\sqrt{\frac{dz}{dx_m}(Q)}=\left(E(P, Q)\sqrt{d x_m}\right)\sqrt{r x_m^{r-1}}$$
and
\be
\frac{\partial}{\partial z_m}\log (E(P, Q)w(Q)^{1/2}) =
\frac{\partial}{\partial z_m}\log\left(E(P, Q)\sqrt{d x_m}\right)
-\frac{r-1}{2rx_m^r}\,.
\la{Evraznyh}\ee

Applying to the left hand side of the last equality the variational formula
(\ref{varprfo}) for $\log E(P, Q)$ (an additional factor $w(P)^{1/2}$
in the left-hand side of (\ref{varprfo}) is inessential, since it is
assumed to be $z_m$-independent) one has
\begin{equation}\label{pf1}
-\frac{1}{4\pi i}\oint_{C}\f{1}{w(R)}\left[\p_R\log\frac{E(P, R)}{E(Q,
R)} \right]^2 =\frac{\p}{\p z_m}\log\left(E(P,
Q)\sqrt{d x_m}\right)
-\frac{r-1}{2r x_m^r}\,.
\end{equation}

Notice that the point $Q$ in the left hand side of (\ref{pf1}) lies outside
the contour $C$. Let $\tilde{C}$ be another contour encircling  $P_m$
such that
the point $Q$ and contour $C$ lie inside of $\tilde{C}$. Using the Cauchy theorem one gets
$$-\frac{1}{4\pi i}\oint_{C}\f{1}{w(R)}\left[\p_R\log\frac{E(P, R)}{E(Q,
R)} \right]^2$$
\be
=-\frac{1}{4\pi
i}\oint_{\tilde{C}}\f{1}{w(R)}\left[\p_R\log\frac{E(P, R)}{E(Q,
R)} \right]^2
-2\pi i\left(-\frac{1}{4\pi i}\right){\rm
Res}\Big|_{R=Q}\left\{ \f{1}{w(R)}\left[\p_R\log\frac{E(P, R)}{E(Q,R)}
  \right]^2\right\}
\la{F0000}\ee
Since the prime-form $E(Q,R)$ behaves as $[z(Q)-z(R)+
  O((z(Q)-z(R))^3]\sqrt{dz(Q)}\sqrt{dz(R)}$ as $R\to Q$, the residue
in (\ref{F0000}) is given by
$$\f{2}{w(Q)}\p_Q \log\{E(P,Q)w^{1/2}(Q)\}\;.$$
Writing down this expression in the local parameter $x_m$ we  rewrite
the right-hand side of
(\ref{F0000}) as follows:
\begin{equation}\label{E111}
-\frac{1}{4\pi
i}\oint_{\tilde{C}}\f{1}{w(R)}\left[\p_R\log\frac{E(P, R)}{E(Q,
R)} \right]^2 -\f{1}{r x_m^{r-1}}\f{d}{d x_m}\left(\log
\{E(P,Q)\sqrt{d x_m}\}
+\frac{r-1}{2}\log x_m\right)\,.
\end{equation}
Since the prime-form is holomorphic at $Q=P_m$, we have
$$\partial_R\left\{\log E(Q, R)\sqrt{dx_m}\right\}=\partial_R\log E(P_m, R)+O(x_m)\;,$$
and, therefore,
\begin{equation}
-\frac{1}{4\pi
i}\oint_{\tilde{C}}\f{1}{w(R)}\left[\p_R\log\frac{E(P, R)}{E(Q,
R)} \right]^2
=-\frac{1}{4\pi
i}\oint_{{C}}\f{1}{w(R)}\left[\p_R\log\frac{E(P, R)}{E(P_m,
R)} \right]^2+O(x_m)\,
\la{intQintPm}
\end{equation}
(the last integral in (\ref{intQintPm}) does not change if we
integrate over $\tilde{C}$ instead of $C$).
Now introducing the expansion
\begin{equation}\label{expan1}
\log\left(E(P, Q)\sqrt{d x_m}\sqrt{w(P)}\right)=e_0+e_1 x_m+\dots +e_r(x_m)^r+O(x_m^{r+1})\,,
\end{equation}
one rewrites (\ref{E111}) as
$$-\frac{1}{4\pi i}\oint_{C}\f{1}{w(R)}\left[\p_R\log\frac{E(P, R)}{E(Q,
R)} \right]^2=-\frac{1}{4\pi
i}\oint_{{C}}\f{1}{w(R)}\left[\p_R\log\frac{E(P_m, R)}{E(Q,
R)} \right]^2-\frac{r-1}{2r(x_m)^r}$$
\begin{equation}\label{E2}
-\frac{e_1+2e_2 x_m+\dots
+re_r x_m^{r-1}}{rx_m^{r-1}}+O(x_m)\,
\end{equation}

On the other hand by virtue of (\ref{expan1}) the right hand side of
(\ref{pf1}) can be rewritten as
\begin{equation}\label{sravn1}
\frac{\p}{\p z_m}\log\left(E(P,
Q)\sqrt{d x_m}\right)-\frac{r-1}{2r x_m^r}
=\frac{\partial
e_0}{\partial z_m}-e_1\frac{1}{rx_m^{r-1}}-e_2\frac{2x_m}{rx_m^{r-1}}-\dots
-e_r\frac{rx_m^{r-1}}{rx_m^{r-1}}+O(x_m)-\frac{r-1}{2rx_m^r}\,.
\end{equation}
Now from (\ref{pf1}), (\ref{E2}) and (\ref{expan1}) it follows that
$$-\frac{1}{4\pi i}\oint_{C}\f{1}{w(R)}\left[\p_R\log\frac{E(P, R)}{E(Q,
R)} \right]^2=\frac{\partial e_0}{\partial z_m}\,$$
which is equivalent to
the statement of the corollary.

$\square$

In the sequel we shall use the following Corollary of formulas
(\ref{varprfo}) and (\ref{varprfo1}):
\begin{corollary}\la{varprfozetak}
The following variational formulas hold:
\be \f{\p \log
\{E(P,P_n)\}}{\p \zeta_k}\Big|_{z(P)}=-\f{1}{4\pi i}
\oint_{\cy_k}\f{1}{w(R)}\left[\p_R\log\f{E(R,P)}{E(R,P_n)}\right]^2
\la{varprfoPPm} \ee
 \be
\f{\p \log \{E(P_l,P_n)\}}{\p \zeta_k}=-\f{1}{4\pi i}
\oint_{\cy_k}\f{1}{w(R)}\left[\p_R\log\f{E(P_l,R)}{E(P_n,R)}\right]^2\;,
\la{varprfoPkPn} \ee
 for any $k=1,\dots, 2g+M-1$, $l,n=1,\dots, M$, $l\neq n$;
 here $E(P,P_n)$ is defined in Corollary \ref{varprfozm};
 $$E(P_l,P_n):=E(P,Q)(d x_l(Q)\,d x_n(P))^{1/2}\Big|_{Q=P_l,\;P=P_n}\;,$$
$x_n(Q)=(z(Q)-z_n)^{1/(k_n+1)}$.
\end{corollary}
{\it Proof.} Notice that in (\ref{varprfo}) one can take $P=P_l$ and $Q=P_n$
with $l\neq n$ for $k=1, \dots, 2g$ and $P=P_l$, $P=P_n$ with $l\neq n$ and
$l, n\neq k-2g+1$ for $k=2g+1, \dots, 2g+M-1$. Namely, consider points $P$
and $Q$ in vicinities of $P_l$ and $P_n$ and apply to them (\ref{varprfo}).
One has $$ E(P, Q)\sqrt{w(P)}\sqrt{w(Q)}=E(P,
Q)\sqrt{dx_l(P)}\sqrt{dx_n(Q)}\sqrt{\frac{dz(P)}{dx_l(P)}}\sqrt{\frac{dz(Q)}{dx_n(Q)}}$$
and
$$\frac{\partial}{\partial \zeta_k}\log E(P, Q)=\frac{\partial}{\partial \zeta_k}\log \left\{E(P,
Q)\sqrt{dx_l(P)}\sqrt{dx_n(Q)}\right\}.$$ Sending $P\to P_l$ and $Q\to P_n$ one gets
the equality
$$\frac{\partial}{\partial \zeta_k}\log E(P_l, P_n)=-\frac{1}{4\pi
i}\oint_{s_k}\frac{1}{w(R)}\left[\p_R\log\frac{E(P_l, R)}{E(
P_n,R)}\right]^2\,.$$ The remaining equations stated in the Corollary can be
proved in the same manner.

$\square$

Dependence of the vector of Riemann constants and differential $\c(P)$
on coordinates $A_\a,
B_\a$ and $z_m$ is given by the following theorem:

\begin{theorem}
The following  variational formulas on the
space $\H(k_1,\dots,k_M)$ hold:
\be\label{varkp}
\f{\p K^P_\a}{\p \ko_k}\Big|_{z(P)}=\f{1}{2\pi i}\oint_{\cy_k}
\f{\hd_\a(R)}{w(R)}\p_R\log\frac{\s(R,Q_0)E(R,P)^{g-1}}{\sqrt{\hd_\a(R)}}
\ee
\be\la{varc}
\f{\p}{\p \ko_k}\log\{\c  w^{\f{g(g-1)}{2}} (P)\}\Big|_{z(P)}=
-\frac{1}{8\pi
i}\oint_{\cy_k}\f{1}{w}\left(S_B-S_{Fay}^{P}\right)
\ee
where $k=1,\dots, 2g+M-1$; the local parameter $z(P)$ is kept fixed under differentiation;  the value of the prime form $E(R, P)$
and the tensor $\s(R, Q_0)$ with respect to arguments $R$ and $Q_0$ respectively
 are calculated in the local parameter $z$.
\end{theorem}
(We notice that the product  of ${\cal C}$ by  a power of $w$  in the left-hand side of
(\ref{varc}) is a
scalar function (i.e. it has zero tensor weight) on $\widehat{\L}$, as well
as the right-hand side.)

{\it Proof.} These formulas are similar to  Fay's   formulas for
variations of $K^P$ and ${\cal C}(P)$ with
respect to variation of the conformal structure on $\L$ defined by an arbitrary Beltrami differential (\cite{Fay92}, pp. 57-59).
Unfortunately,  Fay's formulas do not directly imply  (\ref{varkp}), (\ref{varc}) due to essentially  different
 fixing of the argument $P$ which we use here. Nevertheless the general framework of \cite{Fay92}
is still applicable and we adopt it in the following proof.

From (\ref{rimco}), (\ref{varw}) and (\ref{varB}) one has
$$\frac{\p K^P_\a}{\p \ko_k}=\frac{1}{2}\oint_{\cy_k}\frac{\hd_\a^2}{w}-\sum_{\b\neq \a; \b=1, \dots, g}\delta_{k\b}\frac{\hd_\b(R^\b)}{w(R^\b)}\int_P^{R^\b}\hd_\a$$
$$-\frac{1}{2\pi i}\sum_{\b\neq \a; \b=1, \dots, g}\oint_{x\in a_\b}\left\{\oint_{Q\in \cy_k}\frac{\hd_\b(Q)\p_x\p_Q\log E(x, Q)}{w(Q)}\right\}\int_P^x\hd_\a$$
\begin{equation}\label{rc1}-\frac{1}{2\pi i}\sum_{\b\neq \a; \b=1, \dots, g}\oint_{x\in a_\b}\hd_\b(x)\int_{R=P}^{R=x}\oint_{Q\in\cy_k}\frac{\hd_\a(Q)\p_R\p_Q\log E(R, Q)}{w(Q)},\end{equation}
where $R^\b=a_\b\cap b_\b$.

Notice that
$$\oint_{x\in a_\b}\oint_{Q\in \cy_k}\left(\frac{\hd_\b(Q)\p_x\p_Q\log E(x, Q)}{w(Q)}\int_P^x\hd_\a\right)=$$
\begin{equation}\label{perstav}\oint_{Q\in \cy_k} \oint_{x\in a_\b}\left(\frac{\hd_\b(Q)\p_x\p_Q\log E(x, Q)}{w(Q)}\int_P^x\hd_\a\right)-2\pi i\delta_{k\b}\frac{\hd_\b(R^\b)}{w(R^\b)}\int_P^{R^\b}\hd_\a, \end{equation}
due to asymptotic expansion (\ref{defproco}) of the canonical meromorphic bidifferential.

\begin{remark}\label{treba} {\rm
Let us comment here on the appearance in the right hand sides of the two
formulas above the second terms which at the first sight look strange. To
differentiate an integral, say $\oint_{a_\beta}Gdz$, over the cycle $a_\beta$
with respect to the variable $A_\beta$ one cuts the surface along the basic
cycles and integrates along the contour $a_\beta$ which now is a part of the
boundary of the fundamental polygon $\Lhat$. Choose a finite cover of the
contour $a_\beta$ by the open intervals $I_k$ such that the map $P\mapsto
z(P)$ is univalent inside each interval and let $\{\chi_j\}$ be the
corresponding partition of unity. Then
$\oint_{a_\beta}Gw=\sum_j\int_{I_j}\chi_j(z)G(z)dz$ and the last integral in
the sum is an integral with variable upper limit: when the coordinate
$A_\beta$ gets an increment this upper limit gets the same increment. Thus,
after differentiation of the integral $\oint_{a_\beta}G$ an extra term
appears: the value of the integrand at the end point of the contour
$a_\alpha$ (that is the point $R^\beta$). It should be noted that the third
term in (\ref{rc1}) implicitly depends on the point $R^\b$: the iterated
integral $\oint_{a_\b}\oint_{b_\b}$ entering this term is singular at the
point of intersection of $a_\b$ and $b_\b$ and its value changes when we move
the contours inside their homology classes changing the point of their
intersection. On the other hand {\it the sum} of the second and the third
terms in the right hand side of (\ref{rc1}) does not depend on $R^\beta$ and
the concrete choice of the contours $a_\b$, $b_\b$ within their homology
classes.

To explain the appearance of the second term in the right hand side of
(\ref{perstav}) we observe that the integrand of the iterated integral
$\oint_{\a_\b}\oint_{b_\b}$ in the left hand side of (\ref{perstav}) has the
second order singularity at the point $R^\b$. Localizing the problem, i. e.
making the contours of integration locally coincide  with the subintervals of
real and imaginary axis containing the origin and writing the integrand as
$$i\frac{v_\b(R^\b)}{w(R^\b)}\left(\int_P^{R^\b}v_\a
\frac{1}{(x-iy)^2}+
O\left(\frac{1}{x-iy}\right)\right)dxdy$$
in a vicinity of $R^\b$,  one sees that after changing of the order of
integration the right hand side of (\ref{perstav}) gets the extra term
$$i\frac{v_\b(R^\b)}{w(R^\b)}\int_P^{R^\b}v_\alpha \left\{\int_{y=-a}^{y=a}dy\int_{x=-a}^{x=a}
\frac{dx}{(x-iy)^2}-\int_{x=-a}^{x=a}dx\int_{y=-a}^{y=a}\frac{dy}{(x-iy)^2}
\right\}=-2\pi i \frac{v_\b(R^\b)}{w(R^\b)}\int_P^{R^\b}v_\a\;,
$$
where we used the fact that the expression in the braces equals $-2\pi$. The
analytic  background of this fact is that the logarithmic expression arising
in the first iterated integral is computed assuming that the branch cut of
the logarithm goes from $0$ to $+i\infty$ along the imaginary axis; in the
second integral the branch cut of the logarithm is chosen along the real axis
from $0$ to $+\infty$. Equivalently, one calculates the first iterated
integral as $-2\int_{-a}^a\f{a\, dy}{y^2+a^2}=-\pi $, while the second
iterated integral gives $2\int_{-a}^a\f{a \,dx}{x^2+a^2}=\pi $. }
\end{remark}

Thus, after changing the order of integration and integration by parts the right-hand
side of (\ref{rc1}) reduces to
$$-\frac{1}{2\pi i}\oint_{\cy_k}\frac{1}{w(Q)}
\left\{-\pi i\hd_\a^2(Q)-\sum_{\b\neq\a}\oint_{a_\b}\left[\p_Q\log E(Q,
x)\hd_\b(Q)\hd_\a(x)- \hd_\b(x)\hd_\a(Q)\p_Q\log\frac{E(Q, x)}{E(Q,
P)}\right]\right\}.
$$
As it is explained in (\cite{Fay92}, p. 58) the quadratic differential in the braces coincides
with $$-\hd_\a(Q)\p_Q\log\frac{\s(Q,Q_0)E(Q,P)^{g-1}}{\sqrt{\hd_\a(Q)}}$$ which gives (\ref{varkp}).

To prove (\ref{varc}) we need the following  lemmas.

\begin{lemma}\label {lemma1} Let  the coordinates $z(P)$ and $z(Q)$ be kept fixed and all the tensor objects with arguments $P, Q$ and $Q_0$ are calculated in the local parameter $z$. Then
\begin{equation}\label{varsigma}
\frac{\p\log \s(P, Q)}{\p\ko_k}=\frac{1}{4\pi
i}\oint_{\cy_k}\frac{1}{w(R)}\p_R\log \frac{E(R, P)}{E(R,
Q)}\partial_R\log\left[\s^2(R, Q_0)E(R, P)^{g-1}E(R,
Q)^{g-1}\right]\;,\end{equation} where the values of  ${\bf s}$ and the
prime form are calculated in the local parameter $z$.
\end{lemma}

{\it Proof.} Assume for simplicity that none of the cycles $a_\a$, $\a=1, \dots, g$ has a nonzero intersection index with $\cy_k$.
(The case with intersections presents no serious difficulty, one should observe that the arising additional terms disappear after the change of order of integration -- cf. (\ref{rc1}) and (\ref{perstav}).)
Using (\ref{sigma}), (\ref{varprfo}) and (\ref{varw}), we get
$$ \frac{\p\log \s(P, Q)}{\p\ko_k}=-\frac{1}{2\pi i}\sum_{\b=1}^g\oint_{x\in a_\b}\oint_{R\in\cy_k}\frac{1}{w(R)}\p_x\p_R\log E(R, x)\hd_\b(R)
\log\frac{E(x, P)}{E(x, Q)}+$$
\begin{equation}\label{sig1}+\frac{1}{4\pi i}\sum_{\b=1}^g \oint_{x\in a_\b}\hd_\b(x)\oint_{R\in \cy_k}\frac{1}{w(R)}\left\{
\left(d_R\log\frac{E(x, R)}{E(P, R)}\right)^2-\left(d_R\log \frac{E(x,
R)}{E(Q, R)}\right)^2 \right\}=:\Sigma_1+\Sigma_2.\end{equation} To simplify
the first sum in (\ref{sig1}) we change the order of integration, integrate
by parts, rewrite the interior integral as an integral over the boundary of
the fundamental domain and (at the final step) apply the Cauchy theorem:
$$\Sigma_1=\frac{1}{2\pi i}\oint_{R\in \cy_k}
\frac{1}{w(R)}\sum_{\b=1}^g\oint_{x\in a_\b}\hd_\b(R)\p_R\log E(R,
x)\p_x\log\frac{E(x, P)}{E(x, Q)}=$$ $$-\frac{1}{(2\pi i)(4\pi i)
}\oint_{R\in\cy_k}\frac{1}{w(R)}\oint_{x\in \partial\Lhat}(\p_R\log E(R,
x))^2\p_x\log\frac{E(x, P)}{E(x, Q)}=$$ $$-\frac{1}{4\pi
i}\oint_{\cy_k}\frac{1}{w(R)}\left[(\p_R\log E(P, R))^2-(\p_R\log E(Q,
R))^2\right]-\frac{1}{4\pi
i}\oint_{\cy_{k}}\left(\frac{d}{dz(R)}\right)^2\log\frac{E(R, P)}{E(R,
Q)}dz(R)=$$
\begin{equation}\label{sig2}
-\frac{1}{4\pi i}\oint_{\cy_k}\frac{1}{w(R)}\left[(\p_R\log E(P,
R))^2-(\p_R\log E(Q, R))^2\right]\,.
\end{equation}

The second equality in the sequence of equalities above follows from
(\ref{primetwist}),  the single-valuedness  of the one-form
$$x\mapsto \p_x\log\frac{E(x, P)}{E(x, Q)}$$
on $\L$ and the relation
$$\oint_{\a_\b}\p_x\log \frac{E(x, P)}{E(x, Q)}=0,$$
which holds due to single-valuedness of the prime form along the $a$-cycles.
The last equality holds since
$$\oint_{\cy_{k}}\left(\frac{d}{dz(R)}\right)^2\log\frac{E(R, P)}{E(R,
Q)}dz(R)\equiv\oint_{\cy_{k}}\p_R\left\{\f{1}{w}\p_R\log\frac{E(R, P)}{E(R,
Q)}\right\} =0\,.$$

The second sum in (\ref{sig1}) transforms as follows
$$\Sigma_2=\frac{1}{4\pi i}\oint_{R\in\cy_k}\frac{1}{w(R)}\sum_{\b=1}^g\oint_{x\in \a_\b}\hd_\b(x)
 \left\{ 2\p_R\log E(x, R)\p_R\log\frac{E(Q, R)}{E(P, R)}+\right.$$
$$\left.+\p_R\log (E(P, R)E(Q, R))\p_R\log\frac{E(P, R)}{E(Q, R)}\right\}=$$
$$\frac{1}{\pi i}\oint_{R\in \cy_k}\frac{1}{w(R)}
\left[
-\frac{1}{2}\p_R\log\frac{E(P, R)}{E(Q, R)}\sum_{\b=1}^g\p_R\oint_{x\in a_\b}\hd_\b(x)\log\frac{E(x, R)}{E(x, Q_0)}+\right.$$
$$\left.+\frac{g}{4}\p_R\log (E(P, R)E(Q, R))\p_R\log\frac{E(P, R)}{E(Q, R)}\right]=
$$
\begin{equation}\label{sig3}
=\frac{1}{4\pi i}\oint_{\cy_k}\frac{1}{w(R)}\left[\p_R\log\frac{E(P, R)}{E(Q,
R)}\p_R\log(\s^2(R, Q_0)E^g(P, R)E^g(Q, R))\right]\,.
\end{equation}
The statement of the lemma follows from (\ref{sig1}), (\ref{sig2}) and
(\ref{sig3}).

$\square$

The next lemma describes the variation of the determinant ${\rm det}\,||\hd_\a(R_\beta)||$ from the denominator of expression (\ref{repC}).
\begin{lemma}\label{lemma2}
Assume that the $z$-coordinates of the points $R_1, \dots, R_g, P$ are moduli-independent. Then
\begin{equation}\label{denden}
\lim _{R_1, \dots, R_g\to P}\frac{\partial \log {\rm
det}\,||\hd_\a(R_\b)||}{\partial \ko_k}=-\frac{1}{2\pi i}\sum_{\a,
\b=1}^g\oint_{\cy_k}\frac{1}{w(R)}
\partial^2_{z_\a z_\b}\log \Theta(K^P-{\cal A}_P(R))\hd_\a(R)\hd_\b(R).
\end{equation}
\end{lemma}
{\it Proof.} Denoting the matrix $||v_\alpha(R_\beta)||$ by ${\mathbb V}$ and
using (\ref{varw}), one has
$$\frac{\partial \log {\rm det}\,{\mathbb V}}{\partial \ko_k}={\rm
Tr}\left\{{\mathbb V}^{-1}\left|\left|\f{1}{2\pi
i}\oint_{\cy_k}\f{\hd_\a(R)\Be(R_\beta,R)}{w(R)}\right|\right|\right\}=$$
$$\frac{1}{2\pi
i}\oint_{s_k}\frac{1}{w(R)}\sum_{\a, \b}({\mathbb V}^{-1})_{\a\b}\Be(R_\beta,
R)v_\a(R).
$$
Due to equation (35) from \cite{Fay73} this expression can be rewritten as
$$-\frac{1}{2\pi
i}\oint_{s_k}\frac{1}{w(R)}\sum_{\a,\b=1}^g\partial^2_{z_\a z_\b}
\log\Theta\left(\sum_{\a=1}^g{\cal A}_P(R_\a)-{\cal A}_P(R)
+K^P\right)v_\a(R)v_\beta(R)\;,$$
and one gets (\ref{denden}) sending $R_1, \dots, R_g$ to $P$, when all
${\cal A}_P(R_\a)\to 0$.

$\Box$


Similarly to \cite{Fay92}, we are to vary the logarithm of the right hand side of expression (\ref{repC}) and pass to the limit $R_1, \dots, R_g\to P$, and then $Q\to P_k$.
In what follows all the tensor objects with arguments $P, Q, Q_0, R_1, \dots, R_g$ are calculated in the local parameter $z$.
Using (\ref{varB}) we can represent the variation of the theta-functional term from the numerator of (\ref{repC}) as follows

$$\partial_{\ko_k}\log \Theta (\sum_{\gamma=1}^{g-1}\Acal_P(R_\gamma)+\Acal_Q(R_g)+K^P\, | \,{\B})=$$
\begin{equation}\label{summadvuh}
\sum_{\alpha=1}^g
\left[\partial_{\ko_k}\int_{Q+(g-1)P}^{\sum_{\gamma=1}^gR_\gamma}\hd_\a+\partial_{\ko_k}K^P_\a
\right]\frac{\partial \log\Theta}{\partial z_\a}+\sum_{\a, \b=1}^g
\frac{\partial \log \Theta}{\partial
\B_{\a\b}}\oint_{\cy_k}\frac{\hd_a(R)\hd_\b(R)}{w(R)}\, .
\end{equation}
We have
$$\partial_{\ko_k}\int_{Q+(g-1)P}^{\sum_{\gamma=1}^gR_\gamma}\hd_\a=\frac{1}{2\pi i}\oint_{\cy_k}\frac{1}{w(R)}\int_{Q+(g-1)P}^{\sum_{\gamma=1}^gR_\gamma}\partial_{R}\partial_x\log E(x, R)\hd_\a(R)=$$
\begin{equation}\label{vintegral} =\frac{1}{2\pi i}\oint_{\cy_k}\frac{1}{w(R)}\left\{\partial_{R}\log E(P, R) \hd_\a(R)-\partial_{R}\log E(Q, R) \hd_\a(R)\right\}+o(1) \end{equation}
as $R_1, \dots, R_g \to P$.

Now from (\ref{summadvuh}), (\ref{vintegral}), (\ref{rimco}), the heat equation for the theta-function
and the obvious relation
$$\partial_{R}\log\Theta(K^P-\Acal_P(R))=-\sum_{\a=1}^g(\log \Theta)_{z_\a}\hd_\a(R)$$
 it follows that
$$\lim_{R_1, \dots, R_g\to P} \partial_{\ko_k}\log \Theta (\sum_{\gamma=1}^{g-1}\Acal_P(R_\gamma)+\Acal_Q(R_g)+K^P\, | \,{\B})= $$
$$=-\frac{1}{2\pi i}\oint_{\cy_k}\frac{1}{w(R)}\left\{\partial_{R}\log \Theta(K^P-\Acal(R))\partial_{R}\log[ \s(R, Q_0)E^g(R, P)]
+\frac{(w(R)\partial_R)([w(R)]^{-1}\partial_R)\Theta(K^P-\Acal(R))}{4\Theta(K^P-\Acal_P(R))}\right.
$$
\begin{equation}\label{ochdlinn}
\left.+\sum_{\a=1}^g\partial_{z_\a}\log \Theta(K^P-\Acal_P(Q))\partial_{R}\log E(Q, R)\hd_\a(R)\right\}+o(1)
\end{equation}
as $Q\to R$.

The variation of remaining terms in the right hand side of (\ref{repC}) is much easier. One has
\begin{equation}\label{zero1}\lim_{R_1, \dots, R_g\to P}\partial_{\ko_k} \sum_{\a<\b}\log E(R_\a, R_\b)=0\, ,\end{equation}
\begin{equation}\label{zero2}\lim_{R_1, \dots, R_g\to P}\partial_{\ko_k}\sum_{\gamma=1}^g\log \s (R_\gamma, P)=0
\end{equation}
\begin{equation}\label{formaQ}
\lim_{R_1, \dots, R_g\to P}\partial_{\ko_k}\sum_{\gamma=1}^g\log E(Q, R_\gamma)=-\frac{g}{4\pi i}\oint_{\cy_k}\frac{1}{w(R)}\left(
\partial_{R}\log\frac{E(Q, R)}{E(P, R)}
\right)^2
\end{equation}
due to (\ref{varprfo}) and Lemma \ref{lemma1}.

Now using (\ref{repC}), summing up (\ref{varsigma}), (\ref{ochdlinn} - \ref{formaQ}) and (\ref{denden}),
cleverly rearranging the terms (as Fay does on p. 59 of \cite{Fay92}) and sending $Q\to R$, we get
$$\partial_{\ko_k}\log \cdiff(P)=$$$$\frac{1}{\pi i}\oint_{\cy_k}\frac{1}{w(R)}\left\{\frac{1}{4}\frac{(w(R)\partial_R)([w(R)]^{-1}\partial_R)\Theta(K^P-\Acal_P(R))}{\Theta(K^P-\Acal_P(R))}-\frac{1}{2}
\partial_{R}\log \Theta(K^P-\Acal_P(R))\partial_{R}\log [\s(R, Q_0)E^g(R, P)]\right.
$$
$$-\frac{1}{4}(w(R)\partial_R)([w(R)]^{-1}\partial_R)\log E(R, P)+\frac{1}{2}\partial_{R}\log\s(R, Q_0)\partial_{R}\log E(R, P)
+\frac{2g-1}{4}(\partial_{R}\log E(R, P))^2    $$
$$-\frac{1}{2}\left[\partial_{R}\log E(R, Q)\left(\sum_{\a=1}^g\partial_{z_\a}\log\Theta(K^P-\Acal_P(Q))\hd_\a(R)+
\partial_{R}\log[\s(R, Q_0)E^g(R, P)]\right)\right.
$$
\begin{equation}\label{pochtifin}
\left. -\frac{1}{2}\frac{(w(R)\partial_R)([w(R)]^{-1}\partial_R)E(R, Q)}{E(R,
Q)}-\sum_{\a, \b=1}^g\partial_{z_\a
z_\b}^2\log\Theta(K^P-\Acal_P(R))\hd_\a(R)\hd_\b(R)\right]_{Q=R}\Big\}\,.
\end{equation}
Due to (\ref{primas}), one has
$$\lim_{Q\to R}\partial_{R}\log E(R, Q)\left(\sum_{\a=1}^g\partial_{z_\a}\log\Theta(K^P-\Acal_P(Q))\hd_\a(R)+
\partial_{R}\log[\s(R, Q_0)E^g(R, P)]\right)$$
$$=\lim_{z(Q)\to z(R)}\frac{1}{z(Q)-z(R)}\left(\partial_{R}\log\frac{\s(R, Q_0)E^g(R, P)}{\Theta(K^P-\Acal(R))}+\right.$$
$$\left. +
\sum_{\a, \b=1}^g\partial_{z_\a z_\b}^2\log\Theta(K^P-\Acal_P(R))v_\a(R)\hd_\b(R)(z(Q)-z(R))+O((z(Q)-z(R))^2)\right)=
$$
$$=\sum_{\a, \b=1}^g\partial_{z_\a z_\b}^2\log\Theta(K^P-\Acal_P(R))\hd_\a(R)\hd_\b(R).$$
Here we made use of the fact that the function
\begin{equation}
\label{Rind}
R\mapsto \frac{\Si(R, Q_0)E^g(R, P)}{\Theta(K^P-\Acal_P(R))}
\end{equation}
for fixed $P$ is  holomorphic (since the zero of multiplicity $g$ at
$R=P$ is canceled by the zero of the same multiplicity of $E^g(R,P)$
while $\Si(R, Q_0)$ is non-singular in $\Lhat$) and single-valued  on $\L$ (using
(\ref{primetwist}) and information about the twists of $\Si$ given after
formula (\ref{sigma}),  one sees that all the monodromies of this
function along the basic cycles are trivial) and, therefore, a constant.
Using (\ref{primas}), we see that
$$\lim_{Q\to R}\frac{w(R)\partial_R([w(R)]^{-1}\partial_R)E(R, Q)}{E(R, Q)}=-\frac{1}{2}[S_B-S_w](R)$$
Thus, the last two lines of (\ref{pochtifin}) simplify to
$-\frac{(S_B-S_w)(R)}{8w(R)}$. Using the $R$-independence of expression
(\ref{Rind}) once again, we may rewrite the remaining part of
(\ref{pochtifin}) as
$$\frac{1}{4}w(R)\p_R\f{1}{w(R)}\partial_{R}\log[\s(R, Q_0)E(R, P)^{g-1}]-\frac{1}{4}(\p_{R}\log[\s(R, Q_0)E(R, P)^{g-1}])^2,$$
which coincides with $\frac{1}{8w}(S_{Fay}-S_w)$ due to relation (\ref{SPSw}).
Formula (\ref{varc}) is proved. $\square$

\begin{corollary}
The variational formula (\ref{varkp}) can be equivalently rewritten as follows:
\be\label{varkp1}
\f{\p K^P_\a}{\p \ko_k}\Big|_{z(P)}=\f{1}{2\pi i}\oint_{\cy_k}
\f{\hd_\a(R)}{w(R)}\p_R\log\frac{\s(R,Q_0)E(R,P)^{g-1}}{\sqrt{w(R)}}
\ee
or
\be\la{kzetak1}
\f{\p{K^P_\a}}{{\p \ko_k}}\Big|_{z(P)}=\f{1}{4\pi i}\oint_{\cy_k}
v_\a(R)\left\{\f{1}{w(R)}\p_R\log\prod_{m=1}^M\left(\f{E(R,P)}{E(R,P_m)}\right)^{k_m}
-2\pi i\f{\langle {\bf r},\,{\bf v}(R)\rangle}{w(R)}\right\}
\ee
where integer vector ${\bf r}$ is defined by (\ref{defvecpq})
\end{corollary}
{\it Proof.}
The difference between (\ref{varkp}) and (\ref{varkp1}) is, up to a constant factor, given by the integral
$$\oint_{\cy_k}\f{\hd_\a(R)}{w(R)}\p_R\log\f{\hd_\a(R)}{w(R)}=
\oint_{\cy_k}\p_R\f{\hd_\a(R)}{w(R)}.$$
Since $\hd_\a(R)/w(R)$ is a single-valued meromorphic function on $\L$, this integral vanishes.

The expression (\ref{kzetak1}) follows from representation (\ref{sPQE}) of the
differential ${\bf s}(P,Q)$ in terms of prime-forms.

$\Box$

\subsection{Relation to Teichm\"uller deformation}

Here we point out a close link of our deformation framework on the
moduli spaces of holomorphic differentials with  Teichm\"uller
deformation.
The existence and uniqueness
theorems of  Teichm\"uller state that any two points in Teichm\"uller space of Riemann surfaces
of given genus are related by so-called Teichm\"uller
deformation  (see for example  \cite{Abikoff}) defined by a holomorphic quadratic
differential $W$
and a real positive number $k$.
For our present purposes we assume that $W=w^2$, where $w$ is a holomorphic differential on $\L$
(for and arbitrary $W$ its ``square root'' $w$ is a holomorphic 1-form on two-sheeted ``canonical
covering'' of $\L$). The form $w$ defines local coordinate
$z(P)=\int_{P_0}^P w$ in a neighborhood of any point $P_0\in\L$. Introduce real coordinates
$(x,y)$: $z=x+iy$. Then Teichm\"uller deformation corresponds to stretching in horizontal direction
with some constant coefficient: $x\to\f{1+k}{1-k} x$, $y\to y$; such stretching is
defined globally on $\L$. The finite Beltrami differential corresponding to such finite variation of
conformal structure is given by $k\f{\bar{w}}{w}$ (\cite{Abikoff}, p.32). Infinitesimally,
when $k\to 0$, the stretching is given by $x\to (1+2k) x$ and  Beltrami
 differential defining  infinitesimal deformation $d/dk$ at $k=0$ is
\be \mu_w=\f{\bar{w}}{w} \la{muwwbar} \ee
Under infinitesimal deformation of
the complex structure by an arbitrary Beltrami differential $\mu$ the
variation of the matrix of $b$-periods is given by the Ahlfors-Rauch formula
(\cite{Nag}, p. 263): \be \delta_{\mu} {\bf B}_{\a\b} :=
\f{d}{dt}\Big|_{t=0}{\bf B}_{\a\b}=\int_{\L} \hd_\a\wedge (\mu_w\hd_\b)
\la{varBmu} \ee Therefore, according to (\ref{varBmu}), variation of the
matrix of $b$-periods under infinitesimal Teichm\"uller deformation  is given
by
\be
\f{\p {\bf B}_{\a\b}}{\p k}\Big|_{k=0} = \int_{\L}\f{\bar{w}
v_{\a}}{w}\wedge v_{\b} =-\int_{\L}\f{v_{\a} v_{\b}}{w}\wedge \bar{w}
\la{pBpK}
\ee
Applying Stokes theorem to the fundamental polygon $\Lhat$ with
deleted neighborhoods of zeros of differential $w$, we further rewrite
(\ref{pBpK}) as an integral over boundary: \be \left\{\oint_{\p\Lhat} -2\pi
i\sum_{m=1}^M {\rm res}|_{P_m}\right\}\f{v_\a v_\b}{w}(P)\int_{P_0}^P \bar{w}
\ee where $P_0$ is an arbitrary basepoint. Since both forms $\f{v_\a
v_\b}{w}$ and $\bar{w}$ are closed outside of zeros of $w$, in analogy to the
standard proof of Riemann bilinear relations (see, e. g., \cite{Nag}, p.
257), choosing $P_0$ to coincide with $P_1$, we rewrite this using the
coordinates (\ref{coordint}) as follows:


\be
\f{\p {\bf B}_{\a\b}}{\p k}\Big|_{k=0}=\sum_{\g=1}^g\left\{\bar{B}_\g \oint_{a_\g}\f{v_\a v_\b}{w}
-\bar{A}_\g \oint_{b_\g}\f{v_\a v_\b}{w}\right\}+
2\pi i \sum_{m=2}^M \bar{z}_m{\rm res}|_{P_m}\f{v_\a v_\b}{w}
\la{BKab}
\ee

On the other hand, we have $\int_{\L}\f{v_{\a} v_{\b}}{w}\wedge w=0$, which, repeating the same computation, implies,
\be
0=\sum_{\g=1}^g\left\{B_\g \oint_{a_\g}\f{v_\a v_\b}{w}
-{A}_\g \oint_{b_\g}\f{v_\a v_\b}{w}\right\}+
2\pi i \sum_{m=2}^M {z}_m{\rm res}|_{P_m}\f{v_\a v_\b}{w}
\la{tojd}
\ee
Adding up (\ref{BKab}) and (\ref{tojd}), we get:
\be
\f{\p {\bf B}_{\a\b}}{\p k}\Big|_{k=0}=2\sum_{\g=1}^g\left\{(\Re{B}_\g) \oint_{a_\g}\f{v_\a v_\b}{w}
-(\Re{A}_\g) \oint_{b_\g}\f{v_\a v_\b}{w}\right\}+
4\pi i \sum_{m=2}^M (\Re{z}_m){\rm res}|_{P_m}\f{v_\a v_\b}{w}
\la{BKab1}
\ee

Let us now verify that our variational formulas (\ref{varB})  for the matrix of b-periods
lead to the same result. Under Teichm\"uller deformation $\Im A_\a$, $\Im B_\a$ and
$\Im z_m$ remain unchanged, and corresponding real parts infinitesimally
multiply  with $1+2k$. Therefore,
\be
\f{\p {\bf B}_{\a\b}}{\p k}\Big|_{k=0}=2\sum_{\g}(\Re A_\g)\f{\p{\bf B}_{\a\b}}{\p (\Re A_\g)}+
2\sum_{\g}(\Re B_\g )\f{\p{\bf B}_{\a\b}}{\p (\Re B_\g)}
+2\sum_{m=2}^M (\Re z_m)\f{\p{\bf B}_{\a\b}}{\p (\Re z_m)}
\ee
in complete agreement  with (\ref{BKab1}) if we take into account that ${\bf B}_{\a\b}$ is independent of $\bar{A}_\a$, $\bar{B}_\a$ and $\bar{z}_m$ (i.e. for example
$\f{\p{\bf B}_{\a\b}}{\p (\Re A_\g)}=\f{\p{\bf B}_{\a\b}}{\p A_\g}$ etc)
 and  substitute here our variational formulas (\ref{varB}).

\section{Bergman  tau-function}

\begin{definition} The Bergman tau-function  $\tau(\L,w)$ on the
stratum
${\cal H}_{g}(k_1, \dots, k_M)$ of the
 space of Abelian differentials
is locally defined by the following
system of equations:
 \begin{equation}\label{tau1}
\frac{\partial \log\tau(\L,w)}{\partial \ko_k}=-\frac{1}{12\pi
i}\oint_{\cy_k}\frac{S_B-S_w}{w}\;,
\end{equation}
where  $k=1,\dots, 2g+M-1$; $S_B$ is the Bergman projective connection; $S_w(x):=
\left\{\int^P w,\,x\right\}$;
the difference between two projective connections $S_B$ and $S_w$ is a
meromorphic quadratic differential
with poles at the zeros of $w$.
\end{definition}

To justify this definition one needs to prove that the system of equations (\ref{tau1}) is compatible.
This follows in principle from the fact that in the sequel we find an
explicit expression for $\tau(\L,w)$. However, the computation of
$\tau(\L,w)$ is rather lengthy and technical, while the
straightforward verification of compatibility of equations
(\ref{tau1}) is simple, and we present it here.

Denote the right-hand sides of equations
(\ref{tau1})
by $H^{\zeta_k}$. In analogy with the construction of the Bergman tau-function on Hurwitz spaces (\cite{MPAG}) we call these quantities {\it Hamiltonians}.
Here it will be necessary to distinguish three groups of the coordinates on ${\cal H}(k_1, \dots, k_M)$, so we shall also use the self-explanatory notation
 $H^{{A}_\alpha}$, $H^{{B}_\alpha}$ and $H^{z_m}$ for these Hamiltonians.

We have to show that
$\frac{\partial H^{{A}_\alpha}}{\partial B_\beta}=\frac{\partial H^{{
B}_{\beta}}}{\partial
A_\alpha}$, $\frac{\partial H^{z_m}}{\partial A_\alpha}=\frac{\partial
H^{{ A}_\alpha}}
{\partial z_m}$, etc. Most of these equations immediately follow from Theorem \ref{varfow}
and the
symmetry of the bidifferential $\Be(P,Q)$.

For example, to prove that
\be
\frac{\partial H^{A_\alpha}}{\partial A_\b}=
\frac{\partial H^{A_{\b}}}{\partial A_\a}
\ee
for $\a\neq\b$
we write down the left-hand side as
\be
\frac{\partial H^{A_\a}}{\p
A_\b}=-\f{1}{4\pi^2}\oint_{a_\a}\oint_{a_\b}\f{\Be^2(P,Q)}{w(P) w(Q)}
\la{HaHb}
\ee
which  is obviously symmetric with respect to interchange of $A_\a$
and $A_\b$ since the cycles $a_\a$ and $a_\b$  always can be chosen
non-intersecting.
  Similarly, one can prove all other symmetry
relations
where  the integration contours don't intersect (interpreting the
residue at $P_m$ in terms of the integral over a small contour
encircling $P_m$).

The only equations which require interchange of the order of integration over
intersecting cycles  are \be \frac{\partial H^{A_\alpha}}{\partial B_\alpha}=
\frac{\partial H^{B_{\alpha}}}{\partial A_\alpha}\;. \la{AaBa} \ee To prove
(\ref{AaBa}) we denote the intersection point of $a_\a$ and $b_\a$ by $Q_\a$;
then we have: \be \f{\p H^{A_\a}}{\p B_\a}\equiv\f{1}{12\pi i}\f{\p}{\p
B_\a}\left\{\oint_{b_\a}\f{S_B-S_w}{w}\right\} =\f{1}{12\pi
i}\f{S_B-S_w}{w}(Q_a)-\f{1}{4\pi^2}\oint_{b_\a}\oint_{a_\a}\f{\Be^2(P,Q)}{w(P)
w(Q)} \la{HaBa} \ee where the value of $1$-form $\f{1}{w}(S_B-S_w)$ at the
point $Q_\a$ is computed in coordinate the $z(P)$. The additional term in
(\ref{HaBa})  arises from dependence of the cycle $b_\a$ in the $z$-plane on
$B_\a$ (the difference between the initial and endpoints of the cycle $b_\a$
in $z$-plane is exactly $B_\a$), which has to be taken into account in the
process of differentiation (cf. the arguments given in Remark \ref{treba}).

In the same way we find that \be \f{\p H^{B_\a}}{\p A_\a}\equiv-\f{1}{12\pi
i}\f{\p}{\p A_\a}\left\{\oint_{a_\a}\f{S_B-S_w}{w}\right\} =-\f{1}{12\pi
i}\f{S_B-S_w}{w}(Q_a)-\f{1}{4\pi^2}\oint_{a_\a}\oint_{b_\a}\f{\Be^2(P,Q)}{w(P)
w(Q)} \la{HaAa} \ee (note the change of the sign in front of the term
$\f{1}{w}(S_B-S_w)(Q_\a)$ in (\ref{HaAa}) comparing with (\ref{HaBa})).
Interchanging the order of integration in, say, (\ref{HaBa}) we come to
(\ref{HaAa}) after elementary analysis of the behavior of the integrals in a
neighborhood of the singular point $Q_\a$. (Notice that near the diagonal
$P=Q$ one has $$\Be^2(z(P),
z(Q))=\frac{1}{(z(P)-z(Q)^4}+\frac{S_B(z(P))}{3}\frac{1}{(z(P)-z(Q))^2}+O\left(\frac{1}{z(P)-z(Q)}\right)$$
and only the second term gives a nontrivial input into the difference
$$\left(\oint_{a_\a}\oint_{b_\a}-\oint_{b_\a}\oint_{a_\a}\right)\f{\Be^2(P,Q)}{w(P)
w(Q)};$$ cf. Remark \ref{treba}.)

 This
completes the proof of existence of the Bergman tau-function defined by
(\ref{tau1}).

\subsection{Global definition of the Bergman tau-function}

The right-hand side of formulas (\ref{tau1}) depends not only on the choice
of the  canonical basis of absolute homologies on the surface $\L$, but also
on mutual positions of the basic cycles and the points of the divisor $(w)$,
i.e. it depends on the choice of the basis  $(a_\a,b_a,l_m)$ in $H_1(\L,
\{P_1, \dots, P_M\}; {\mathbb Z})$.

However, it turns out that dependence on the choice of contours $\{l_m\}$ is
in fact absent, and one possible global definition of the tau-function could be as a
horizontal section of some (flat) line bundle ${\cal T}$ over the covering
${\widehat {\cal H}}_{g}(k_1, \dots, k_M)$ of the space ${\cal H}_{g}(k_1,
\dots, k_M)$. Here  ${\widehat {\cal H}}_{g}(k_1, \dots, k_M)$ is the space
of triples $(\L, w, \{a_\a, b_\a\})$, where  $\{a_\a, b_\a\}$ is a canonical
basis in the first homologies $H_1(\L, {\mathbb Z})$. In the trivial line
bundle ${\widehat {\cal H}}_{g}(k_1, \dots, k_M)\times {\mathbb C}$ introduce
the connection
\begin{equation}\label{connection}
d_B=d-\sum_{k=1}^{2g+M-1}H^{\zeta_k}d\zeta_k.
\end{equation}
 (Here $d$
is the external differentiation having both ``holomorphic" and
``antiholomorphic" components.) The  Lemma \ref{neza} below shows that this connection is well-defined on ${\widehat {\cal H}}_{g}(k_1, \dots, k_M)$ i. e.
expression (\ref{connection}) is independent of
the choice of contours $l_m$ connecting the zeros $P_1$ and $P_m$.

Let two systems of  cuts on $\L$: $\{a_\a, b_\a\}$ and $\{ a_\a',  b_\a'\}$
 define the same canonical basis in $H_1(\L, {\mathbb Z})$. Notice that the cycles $a_\a$ and
$a_\a'$ (as well as $b_\a$ and $b_\a'$)
are not necessarily equivalent as  elements of $H_1(\L, \{P_1, \dots, P_M\}; {\mathbb Z})$.
Let  $\Lhat$ and $\Lhat'$ be the corresponding fundamental polygons and let
$\{\zeta_k\}=\{A_\alpha, B_\alpha, z_m\}$, $\{\zeta_k'\}=\{ A_\alpha',  B_\alpha',  z_m'\}$
be the corresponding systems of local coordinates on ${\cal H}_{g}(k_1, \dots, k_M)$.
We recall that when defining the coordinate $z_m$ (or $ z_m'$) we integrate the differential
$w$    over a contour $l_m$ (or $l_m'$) connecting the zeros $P_1$ and $P_m$ and lying {\it inside} the
fundamental polygon $\Lhat$ (or $\Lhat'$).
Let also $H^{\zeta_k}$
and $H^{\zeta_k'}$ be the corresponding Hamiltonians.
\begin{lemma}\label{neza}
The following equality holds
\begin{equation}\label{forma}
\sum_{k=1}^{2g+M-1}H^{\zeta_k}d\zeta_k=\sum_{k=1}^{2g+M-1}H^{\zeta_k'}d\zeta_k'\, .
\end{equation}
\end{lemma}
{\it Proof.} We may deform one system of cuts
(keeping it defining the same canonical basis in $H_1(\L, {\mathbb Z})$) into another
 through a sequence of elementary moves: each elementary move corresponds to
passing of a  chosen zero $P_k$ of $w$ from one shore of some cut to another.
 It is sufficient to show that (\ref{forma}) holds if the system of cuts $\{ a_\a',  b_\a'\}$ can be obtained from
the system $\{a_\a, b_\a\}$ via one elementary move.

Let  the zero $P_k$ pass from the right shore of the (oriented) cut $a_\gamma$ to  its left shore.
Due to the Cauchy theorem we have
\begin{equation}\label{Ham111} H^{ B_\gamma'}=H^{B_\gamma}+H^{z_k}\end{equation}
and all other Hamiltonians do not change.
The coordinate $z_k$ transforms to
\begin{equation}\label{coord111}
z_k'= z_k- B_\gamma\,
\end{equation}
and all other coordinates do not change.
Equation (\ref{forma}) immediately follows from (\ref{Ham111}) and (\ref{coord111}).

Let  the zero $P_k$ pass from the right shore of the (oriented) cut $b_\gamma$
to  its left shore. Then
\begin{equation}\label{Ham22} H^{ A_\gamma'}=H^{A_\gamma}-H^{z_k}\end{equation}
and all other Hamiltonians do not change.
The coordinate $z_k$ transforms to
\begin{equation}\label{coord22}
z_k'=z_k+A_\gamma\,
\end{equation}
and all other coordinates do not change.
Equation (\ref{forma}) again holds.  $\square$

The compatibility of equations (\ref{tau1}) provides flatness of  connection (\ref{connection}).

The flat connection $d_B$ determines a character of the fundamental
group of ${\widehat {\cal H}}_{g}(k_1, \dots, k_M)$ i.e. the
representation
\begin{equation}\label{repr}
\rho:\pi_1\big({\widehat {\cal H}}_{g}(k_1, \dots, k_M)\big)\rightarrow {\mathbb C}^*\;.
\end{equation}

Denote by  ${\cal U}$  the universal covering of ${\widehat {\cal H}}_{g}(k_1, \dots, k_M)$; then
the group $\pi_1\big({\widehat {\cal H}}_{g}(k_1, \dots, k_M)\big)$ acts on the direct product
${\cal U}\times {\mathbb C}$ as follows:
$$g(u, z)=(gu, \rho(g)z)\;,$$
where $u\in {\cal U}$, $z\in {\mathbb C}$, $g\in \pi_1\big({\widehat {\cal
H}}_{g}(k_1, \dots, k_M))\big)$. The factor space $\;\Big({\cal
U}\times{\mathbb C}\Big)/\pi_1\big({\widehat {\cal H}}_{g}(k_1, \dots,
k_M)\big)$ has the structure of a holomorphic line bundle over ${\widehat
{\cal H}}_{g}(k_1, \dots, k_M)$; we denote this bundle by ${\cal T}$. Now the
local definition \ref{tau1} of the Bergman tau-function can be reformulated
as follows:
\begin{definition}\rm
{The flat holomorphic line bundle ${\cal T}$ equipped with the flat
connection $d_B$ is called
the Bergman line bundle over the  space ${\widehat {\cal H}}_{g}(k_1, \dots, k_M)$.
The (unique up to a multiplicative constant) horizontal holomorphic
section of the bundle  ${\cal T}$ is called
the Bergman $\tau$-function.}
\end{definition}

\subsection{Explicit formula for the Bergman tau-function}

Here we are going to give an explicit formula for the Bergman
tau-function.
As the first step we rewrite the definition of the tau-function
(\ref{tau1}) can be rewritten as follows:
 \begin{equation}\label{tau1mod}
\frac{\partial \log\tau(\L,w)}{\partial \ko_k}=-\frac{1}{12\pi
i}\oint_{\cy_k}\frac{S_B-S_{Fay}^{P}}{w}-
\frac{1}{12\pi i}\oint_{\cy_k}\frac{S_{Fay}^{P}-S_w}{w}\;,
\end{equation}
where $S_{Fay}^P$ is Fay's projective connection (\ref{Fc}).
The first term in (\ref{tau1mod}) can be integrated in terms of
differential ${\cal C}(P)$ (\ref{c}) using the variational formula
(\ref{varc}).

To formulate the theorem giving the antiderivative of the second term
in (\ref{tau1mod}) we introduce two vectors ${\bf r}$ and ${\bf q}$ with integer
coefficients such that for a given choice of the fundamental cell $\Lhat$
\be
\Abel_P\big((w)\big)+2K^{P}+{\bf B}{{\bf r}}+{\bf q}=0\;.
\la{defrq}\ee
\begin{theorem}\la{derivG}
For any point $P\in\L$ not coinciding with any $P_m$ introduce the
following function ${\cal G}(P)$ on $\Lhat$:
\be
{\cal G}(P)=e^{8\pi i\langle{\bf r},K^P\rangle+ 2\pi i \langle{\bf r},
  {\bf B}{{\bf r}}\rangle}
[w(P)]^{(2g-2)^2}\left\{\prod_{m=1}^M
E^{k_m}(P,P_m)\right\}^{4g-4}\prod_{m,n=1\,m< n}^M E^{-2 k_m k_n} (P_m,P_n)
\la{defG}
\ee
Then the following variational formulas hold:
\be
\frac{\partial \log {\cal G}(P)}{\partial \ko_k}\Big|_{z(P)}=\frac{1}{\pi i}\oint_{\cy_k}\frac{S_{Fay}^{P}-S_w}{w}
\la{varforQ}
\ee
\end{theorem}
{\it Proof.}
To simplify our computation in this proof for any $Q\in \L$ we
introduce the 1-forms $f_Q$ (these forms are meromorphic on $\Lhat$,
but their combinations arising below are all meromorphic one-forms on
$\L$ itself). If $Q$ does not coincide with any $P_m$, $f_Q(R)\equiv \p_R\log\{E(R,Q)w^{1/2}(R) w^{1/2}(Q)\}$.
For $Q=P_m$ we define $f_{P_m}(R)=\p_R\log\{E(R,P_m) w^{1/2}(R)\}$.

To compute the left-hand side of (\ref{varforQ}) we use variational
formulas ,(\ref{varprfoPPm}), (\ref{varprfoPkPn}) (\ref{varB}) and
(\ref{varkp})  for the prime-form, $K^P$ and ${\bf B}$. Using
(\ref{varprfoPPm}) and  (\ref{varprfoPkPn})
 we get:
$$
\frac{\partial \log {\cal G}(P)}{\partial \ko_k}\Big|_{z(P)}=
-\f{1}{4\pi i}\oint_{\cy_k}\f{1}{w}\left\{(4g-4)\sum_{m=1}^M
k_m(f_P-f_{P_m})^2- 2\sum_{m<n}k_m k_n (f_{P_n}-f_{P_m})^2\right\}
$$
$$
+8\pi i\left\langle {\bf r}, \,\f{\p{K^P}}{{\p \ko_k}}\right\rangle + 2\pi i \left\langle
{\bf r},\,\f{\p{{\bf B}}}{{\p \ko_k}}{\bf r}\right\rangle
$$
For $\p{{\bf B}}/{\p \ko_k}$ we shall use the variational formula
(\ref{varB}); for $\p{K^P_\a}/{\p \ko_k}$ we shall use the formula (\ref{kzetak1}).

From (\ref{kzetak1}) we have:
\be
\f{\p{\langle K^P,\,{\bf r}\rangle}}{{\p \ko_k}}=\f{1}{4\pi
  i}\oint_{\cy_k}
\f{\langle {\bf r},\,{\bf
    v}(R)\rangle}{w(R)}\left\{(2g-2)f_{P}-\sum_{m=1}^M k_m f_{P_m}
-2\pi i \langle {\bf r},\,{\bf v}(R)\rangle\right\}
\ee
Taking into account (\ref{varB}), we get
\be
\f{\p\langle
{\bf r},\,{\bf B}{\bf r}\rangle
}{{\p \ko_k}}= \oint_{\cy_k}\f{\langle {\bf r},\,{\bf
    v}(R)\rangle^2}{w}
\ee
Let us observe now that the first term in (\ref{kzetak1}) can be
rewritten as
$$
-\f{1}{2\pi i}\oint_{\cy_k}\f{1}{w(R)}\left\{(2g-2)f_P-\sum_{m=1}^M
k_m f_{P_m}\right\}^2
$$
Now  (\ref{kzetak1}) can be rewritten as follows:
$$
\frac{\partial \log {\cal G}(P)}{\partial \ko_k}\Big|_{z(P)}
=-\f{1}{2\pi i}\oint_{\cy_k}\f{1}{w(R)}\left\{(2g-2)f_P-\sum_{m=1}^M
k_m f_{P_m}\right\}^2
$$
$$
+2\oint_{\cy_k}\f{\langle {\bf r},\,{\bf
  v}(R)\rangle}{w(R)}\left\{(2g-2)f_P-\sum_{m=1}^M k_m f_{P_m}\right\}-2\pi i
\oint_{\cy_k}\f{\langle {\bf r},\,{\bf
    v}(R)\rangle^2}{w(R)}
$$
\be
=-\f{1}{2\pi i}\oint_{\cy_k}\f{1}{w(R)}\left\{(2g-2)f_P-\sum_{m=1}^M
k_m f_{P_m}-2\pi i \langle {\bf r},\,{\bf v}(R)\rangle\right\}^2
\la{dGdzetacal}
\ee

\vskip0.5cm
Consider now the right-hand side of (\ref{varforQ}).
Using formula (\ref{sPQE}) for the differential ${\bf s}$ we have:
\be
\f{1}{w(R)}\p_R\log\left\{{\bf s}^2(R,Q_0)E^{2g-2}(R,P)\right\}=
\f{1}{w(R)}\p_R\log\prod_{m=1}^M\left\{\f{E(R,P)}{E(R,P_m)}\right\}^{k_m}
-2\pi i \f{\langle {\bf r},\,{\bf v}\rangle}{w}
\ee
Substituting this expression into  representation   (\ref{SPSw}) of
the 1-form $(S_{Fay}^P-S_w)/w$,
we get
$$
\frac{1}{\pi i}\oint_{\cy_k}\frac{S_{Fay}^{P}-S_w}{w}=-\f{1}{2\pi i}
\oint_{\cy_k}\f{1}{w(R)}\left\{\p_R\log\prod_{m=1}^M
\left\{\f{E(R,P)}{E(R,P_m)}\right\}^{k_m}-2\pi i \langle {\bf r},\,{\bf v}\rangle\right\}^2
$$
\be +\f{1}{\pi i}\oint_{\cy_k}{\p_R}\left\{\f{1}{w(R)}\p_R\log\prod_{m=1}^M
\left[\f{E(R,P)}{E(R,P_m)}\right]^{k_m}-2\pi i \f{\langle {\bf
r},\,{\bf v(R)}\rangle}{w(R)}\right\} \la{rhscalcul} \ee The first integral
in the right-hand side of (\ref{rhscalcul}) coincides with the right-hand
side of (\ref{dGdzetacal}). The second
integral in the right-hand side of (\ref{rhscalcul}) vanishes, since it is an
integral of the derivative of the meromorphic function in the braces   over a
closed contour. The theorem is proved. 

$\square$

Now from variational formula for differential ${\cal C}$ (\ref{varc})
and Theorem \ref{derivG} we get the formula for Bergman tau-function:
\be
\tau(\L,w)= ({\cal G}(P))^{-1/12}\left({\cal C}(P) \{w(P)\}^{g(g-1)/2}\right)^{2/3}
\la{taupredv}
\ee
We notice that the expression if the right-hand side of
(\ref{taupredv}) is in fact  independent of the choice of point $P$.
Taking into account expression for ${\cal G}(P)$ (\ref{defG}), we come
to the following theorem:
\begin{theorem}\la{tauquadr}
The Bergman tau-function on the space ${\cal H}_g(k_1, \dots, k_M)$ is
 given by the following formula:
\be
 \tau(\L,w)= {\cal F}^{2/3} e^{-\f{\pi i}{6}\langle {\bf r},\, {\bf
 B}{\bf r}\rangle}\prod_{m,n,\,m<n}
\{E(P_m,P_n)\}^{k_m k_n/6}
\la{*}
\ee
where the function ${\cal F}$ defined by (\ref{const0}):
\be
{\cal F}= [w(P)]^{\f{g-1}{2}}e^{-{\pi i} \langle\rb,
K^P\rangle}
\left\{\prod_{m=1}^{M}[ E (P,P_m)]^{\f{(1-g)k_m}{2}}\right\} \c(P)
\la{const}
\ee
is independent of $P$;
the integer vector ${\bf r}$ is defined by the equality
\begin{equation}
\Abel\big((w)\big)+2K^{P}+{\bf B}{\bf r}+{\bf q}=0\;;
\la{defvecr}
\end{equation}
${\bf q}$ is another integer vector, $(w)$ is the divisor of the
differential $w$,  the initial point of the Abel map $\Abel$
coincides
with $P$ and all the paths are chosen inside the same fundamental
polygon $\widehat{\L}$.
\end{theorem}

The expression (\ref{*}), (\ref{const}) for the Bergman tau-function
can be slightly simplified
for the case of the highest stratum ${\cal H}_g(1,\dots,1)$.

\begin{lemma}\la{fundcell}
Let all the zeros of the Abelian differential $w$ be simple.  Then the
fundamental cell
$\Lhat$ can always be chosen such that $\Acal((w))+2K^P=0$.
\end{lemma}
{\it Proof.}  For an arbitrary choice of the fundamental cell we can
claim that
 the vector $\Acal((w))+2K^P$  coincides with $0$ on the Jacobian of
the surface $\L$, i.e.
there exist two integer vectors ${\bf r}$ and ${\bf q}$ such that
$\Abel((w))+2K^P+{\bf B}{\bf r}+{\bf q}=0$.
Fix some zero $P_k$ of $w$; according to our assumption this zero is
simple.
By a smooth deformation of a cycle $a_\a$ within a given homological
class we can stretch this
cycle in such a way that the point $P_k$ crosses this cycle; two possible
directions of the crossing
correspond to the jump of component ${\bf r}_\a$ of the vector ${\bf
r}$ to $+1$ or $-1$.
Similarly, if we deform a cycle $b_\a$ in such a way that it is crossed by the point
$P_k$, the component
${\bf q}_\a$ of the vector ${\bf q}$  also jumps by $\pm 1$. Repeating
such procedure, we
come to fundamental domain where ${\bf r}={\bf q}=0$.

$\Box$

From the proof it is clear that even a stronger statement is true: the
choice of the fundamental
domain such that $\Abel((w))+2K^P=0$ is always possible if the
differential $w$ has {\it at least one}
simple zero.

\begin{corollary} \la{Corotvtau}
Consider the highest stratum ${\cal H}(1,\dots,1)$ of the space ${\cal
H}_g$
containing Abelian differentials $w$ with simple zeros.
Let us choose  the fundamental cell $\widehat{\L}$  such that
$\Abel((w))+2K^P=0$.
Then the
Bergman tau-function on ${\cal H}(1,\dots,1)$ can be written as
follows:
\begin{equation}\label{otvtau}
\tau(\L,w) = {\Fcal}^{2/3}  \prod_{m,l=1\;\; m < l}^{2g-2} [E(P_m,P_l)]^{{1}/{6}}
\end{equation}
where expression
\be
\Fcal := [w(P)]^{\f{g-1}{2}}\c(P)
\prod_{m=1}^{2g-2}[ E (P,P_m)]^{\f{(1-g)}{2}}
\ee
does not depend on $P$; all prime-forms are evaluated at the points
$P_m$ in the distinguished local parameters $x_m(P)=\left(\int_{P_m}^P
w\right)^{1/2}$.
\end{corollary}

The following corollary describes the dependence of the Bergman
tau-function on the choice of holomorphic differential assuming that
he Riemann surface remains the same. For simplicity we assume that all
zeros of both holomorphic differentials are simple, and none of the
zeros of the first differential coincides with a zero of the second
differential.
This corollary will be used below in deriving formulas of Polyakov
type, which describe the dependence of the determinant of Laplacian on
the choice of flat conical metric on a fixed Riemann surface.

\begin{corollary}\la{theotautau1}
Let $w$ and $\wt$ be two holomorphic 1-forms with simple zeros on the
same Riemann surface $\L$; assume that all of these zeros are different. Introduce their divisors
$(w):=\sum_{m=1}^{2g-2} P_m$ and $(\wt):=\sum_{m=1}^{2g-2}
\Pt_m$. Then
\be
\f{\tau(\L,w)}{\tau(\L,\wt)}= \prod_{m=1}^{2g-2}
\left\{\f{\res|_{\Pt_m}\{w^2/\wt\}}
 {\res|_{\P_m}\{\wt^2/w\}}\right\}^{1/24}\;.
\la{zamena}
\ee
\end{corollary}
{\it Proof.}
The distinguished local parameter in a neighborhood of $P_m$ is
$x_m(P):=\left[\int_{P_m}^P w\right]^{1/2}\;$;
in a neighborhood of $\Pt_m$ the distinguished local parameter is
$\xt_m(P):=\left[\int_{\Pt_m}^P w\right]^{1/2}$.
Then the formula (\ref{zamena}) can be alternatively rewritten as
follows:
$$
{\tau(\L,w)}\prod_{m=1}^{2g-2}\wt^{1/12}(\P_m)=
{\tau(\L,\wt)}\prod_{m=1}^{2g-2} w^{1/12}(\Pt_m)\;,
$$
where we use the  standard convention  for evaluation of
the differentials $w$ and $\wt$ at their zeros:
\be
\wt (\P_m):= \f{\wt(P)}{dx_m (P)}\Big|_{P=P_m}\;,\hskip0.7cm
w(\Pt_m):= \f{w(P)}{d\xt_m (P)}\Big|_{P=\Pt_m}\;.
\ee
Let us  assume that the fundamental cell $\hat{\L}$ is chosen in such
a way that the Abel maps of divisors
$(w)$ and $(\wt)$ equal $2K^P$; this choice is always possible (see
Lemma \ref{fundcell}) in
our present case, when all
points of these divisors have multiplicity $1$.
Then we get, according to the formulas
(\ref{otvtau}) (all products below are taken from $1$ to $2g-2$):
\be
\f{{\tau^{12}(\L,w)}\prod_{m}\wt(\P_m)}{{\tau^{12}(\L,\wt)}\prod_{m}
w(\Pt_m)}
=\prod_{m}\f{\wt(\P_m)}{w(\Pt_m)}\prod_{m<n}\f{E^2(\P_m,\P_n)}{E^2(\Pt_m,\Pt_n)}
\left\{\f{w(P)\prod_{m} E(P,\Pt_m)}{\wt(P)\prod_{m}
E(P,\P_m)}\right\}^{4g-4}\;.
\la{ratio}
\ee
Since this expression is independent of $P$, we can split the power
$4g-4$ of the expression in the braces
into product over arbitrary $4g-4$ points, in particular, into product
over $\P_1,\dots,\P_{2g-2}$ and
$\Pt_1,\dots,\Pt_{2g-2}$. Then most of the terms in (\ref{ratio})
cancel each other. The only terms left are
due to the fact that the prime-forms vanish at coinciding arguments;
this compensates vanishing of $w$ and $\wt$
at their zeros. As a result we can rewrite (\ref{ratio}) as follows:
\be
\prod_m\left\{\lim_{P\to
\P_m}\f{w(P)}{E(P,\P_m)(dx_m(P)^{3/2}}\lim_{P\to
\Pt_m}\f{E(P,\Pt_m)(d\xt_m(P))^{3/2}}{\wt(P)}\right\}\;,
\ee
which equals $1$, since, say,  in a neighborhood of $\P_m$ we have
$w(P)=2 x_m(P) dx_m(P)$ and
$E(P,\P_m)=x_m(P)/\sqrt{dx_m(P)}$.

$\Box$

\begin{remark}{\rm
In the early version of this paper Theorem \ref{derivG} (which is the key
point of the proof of the explicit expression for the tau-function) was
proved in an indirect way, parallel to the proof of the formula for Bergman
tau-function on Hurwitz spaces \cite{IMRN}. Namely, it was shown that the
modulus square $|{\cal G}(P)|^2$ of the function ${\cal G}$ from (\ref{defG})
up to a moduli independent constant coincides with the properly regularized
Dirichlet integral
$${\mathbb D}=\frac{1}{\pi}\iint_\L|\partial \phi|^2,$$
where $\phi=\log\big|\frac{\Omega^P}{w}\big|^2$ and the one-form $\Omega^P$
is defined in (\ref{Om}). This explains how one can guess expression
(\ref{defG}): this guess is based on general idea (coming from string theory) that
Dirichlet and Liouville integrals arise in integrating projective
connections. After that via a rather complicated calculation it was shown
that the Dirichlet integral ${\mathbb D}$ satisfies the same system of
equations (\ref{varforQ}) as the function ${\cal G}$.}
\end{remark}

\section{Determinants of Laplacians in the metrics $|w|^2$}
\subsection{Laplacians on polyhedral surfaces. Basic facts}

Any holomorphic Abelian differential $w$
 defines a natural flat metric  on the Riemann surface $\L$ given by
$|w|^2$. This metric has  conical
 singularities at the zeroes of $w$.
The cone angle of the metric $|w|^2$ equals $2(k+1)\pi$ at the zero of $w$ of
multiplicity $k$. The surface $\L$ provided with metric $|w|^2$ is a special
case of a {\it compact polyhedral surface}, i. e. a two dimensional compact
Riemannian manifold provided with flat metric with conical singularities (any
such surface can be glued from Euclidean triangles, see \cite{Troyanov}).

Here we give a short  self-contained survey of some basic facts from the
spectral theory of Laplacian on compact polyhedral surfaces. We start with
recalling the (slightly modified) Carslaw construction (1909) of the heat
kernel on a cone. Then we describe all self-adjoint extensions of
conical Laplacian (these results are complementary to Kondratjev's study
\cite{Kondr} of elliptic equations on conical manifolds and are well-known,
being in the folklore since sixties; their generalization to the case of
Laplacians acting on $p$-forms can be found in \cite{Mooers}). Finally, we
establish the precise heat asymptotics for the Friedrichs extension of the
Laplacian on a compact polyhedral surface. More
general results on the heat asymptotics for Laplacians acting on $p$-forms on
piecewise flat pseudomanifolds can be found in \cite{Cheeger}.

\subsubsection{The heat kernel on infinite cone}

 We start from the standard heat kernel
  \begin{equation}
\label{h1}H_{2\pi}( \x,  \y; t)=\frac{1}{4\pi
 t}\exp\left\{\f{-(\x-\y)\cdot(\x-\y)}{4t}\right\}\;\hskip0.7cm {\bf
 x},\,{\bf y}\in {\mathbb R}^2
\end{equation}
 in  ${\mathbb R}^2$ which we consider as the cone with conical angle $2\pi$.
  Introducing the polar coordinates $(r, \theta)$ and $(\rho,
 \psi)$ in the $\x$ and $\y$-planes respectively, one can rewrite (\ref{h1}) as the contour integral
\begin{equation}\label{h2}
H_{2\pi}( \x,  \y;
t)=\frac{1}{16\pi^2it}\exp\left\{\f{-(r^2+\rho^2)}{4t}\right\}\int_{C_{\theta,
\psi}}\exp\left\{\f{r\rho\cos(\alpha-\theta)}{2t}\right\}\cot
\frac{\alpha-\psi}{2}\,d\alpha\;,
\end{equation}
where $C_{\theta, \psi}$
denotes the union of a small positively oriented circle centered at
$\alpha=\psi$ and the two vertical lines, $l_1=(\theta-\pi-i\infty,
\theta-\pi+i\infty)$ and $l_2=(\theta+\pi+i\infty, \theta+\pi-i\infty)$,
having mutually opposite orientations.

To prove (\ref{h2}) one has to notice that

1) $\Re \cos(\alpha-\theta)<0$ in vicinities of the lines $l_1$ and $l_2$
and, therefore, the integrals over these lines converge.

2) The integrals over the lines cancel due to $2\pi$-periodicity of the
integrand and the remaining integral over the circle coincides with
(\ref{h1}) due to the Cauchy Theorem.

Observe that one can deform the contour $C_{ \theta, \psi}$ into the union,
$A_{\theta}$, of two contours lying in the open domains
$\{\theta-\pi<\Re\alpha<\theta+\pi\,,\, \Im \alpha>0\}$ and
$\{\theta-\pi<\Re\alpha<\theta+\pi\,,\,  \Im \alpha<0\}$ respectively. The
first contour goes from $\theta+\pi+i\infty$ to $\theta-\pi+i\infty$, the
second one goes from $\theta-\pi-i\infty$ to $\theta+\pi-i\infty$. This leads
to the following alternative integral representation for the heat kernel $H_{2\pi}$:
\begin{equation}\label{h3} H_{2\pi}( \x,  \y;
t)=\frac{1}{16\pi^2it}\exp\left\{\f{-(r^2+\rho^2)}{4t}\right\}\int_{A_
\theta}\exp\left\{\f{r\rho\cos(\alpha-\theta)}{2t}\right\}\cot
\frac{\alpha-\psi}{2}\,d\alpha\;.
\end{equation}

The latter representation admits natural generalization to the case of the
cone $C_\beta$ with conical angle $\beta$, $0<\beta<+\infty$:
$$C_\beta=\{(r, \theta) \ :\  r\in [0, \infty), \ \theta\in {\mathbb
R}/\beta{\mathbb Z}\}/(0, \theta_1)\sim (0, \theta_2)$$ equipped with the
metric $(dr)^2+r^2(d\theta)^2$; notice here that in case $0<\beta\leq 2\pi$
the cone $C_\beta$ is isometric to the surface $
z_3=\sqrt{(\frac{4\pi^2}{\beta^2}-1)(z_1^2+z_2^2)}$.

Namely, introducing the polar coordinates on $C_\beta$, we see that the
following expression represents the heat kernel on $C_\beta$:
\begin{equation}\label{h4}
H_{\beta}(r, \theta, \rho, \psi; t)=\frac{1}{8\pi\beta
it}\exp\left\{-\f{r^2+\rho^2}{4t}\right\}\int_{A_
\theta}\exp\left\{\f{r\rho\cos(\alpha-\theta)}{2t}\right\}\cot
\frac{\pi(\alpha-\psi)}{\beta}\,d\alpha\,.
\end{equation}

Clearly, expression (\ref{h4}) is symmetric with respect to $(r, \theta)$ and
$(\rho, \psi)$ and is $\beta$-periodic with respect to the angle variables
$\theta, \psi$. Moreover, it satisfies the heat equation on $C_\beta$.
Therefore,  to verify that $H_{\beta}$ is in fact the heat kernel on
$C_\beta$ it remains to show that $H_\beta(\cdot, y, t)\longrightarrow
\delta(\cdot-y)$ as $t\to 0+$. To this end deform the contour $A_\psi$ into
the union of the lines $l_1$ and $l_2$ and (possibly several) small circles
centered at the poles of $\cot \frac{\pi(\cdot-\psi)}{\beta}$ in the strip
$\theta-\pi<\Re\alpha<\theta+\pi$. The integrals over all the components of
this union except the circle centered at $\alpha=\psi$ vanish in the limit as
$t\to 0+$, whereas the integral over the latter circle coincides with
$H_{2\pi}$.

\subsubsection{The heat asymptotics near the vertex}
\begin{proposition}\label{ugol}
Let $R>0$ and $C_{\beta}(R)=\{\x\in C_{\beta}: {\rm dist}(\x, {\cal O})<
R\}$, where ${\cal O}$ is the conical point.
Let also $d\x$ denote the area element on $C_\beta$. Then for some
$\epsilon>0$
\begin{equation}\label{as1}\int_{C_\beta(R)}H_{\beta}(\x, \x; t)\,d\x=\frac{1}{4\pi
t}{\rm Area}(C_{\beta}(R))+
\frac{1}{12}\left(\frac{2\pi}{\beta}-\frac{\beta}{2\pi}
\right)+O(e^{-\epsilon/t})\end{equation} as $t\to 0+$.
\end{proposition}

{\it Proof}\,(cf. \cite{Fursaev}, p. 1433). Make in (\ref{h4}) the change of
variable $\gamma=\alpha-\psi$ and
 deform the contour
$A_{\theta-\psi}$ into the contour $\Gamma^-_{\theta-\psi}\cup
\Gamma^{+}_{\theta-\psi}\cup \{ |\gamma|=\delta\}$, where the oriented curve
$\Gamma^{-}_{\theta-\psi}$ goes from $\theta-\psi-\pi-i\infty$ to
$\theta-\psi-\pi+i\infty$ and intersects the real axis at $\gamma=-\delta$,
the oriented curve $\Gamma^+_{\theta-\psi}$ goes from
$\theta-\psi+\pi+i\infty$ to $\theta-\psi+\pi-i\infty$ and intersects the
real axis at $\gamma=\delta$, the circle $\{|\gamma|=\delta\}$ is positively
oriented and $\delta$ is a small positive number. Calculating the integral
over the circle $\{|\gamma|=\delta\}$  via the Cauchy Theorem, we get
\begin{equation}\label{h6}
H_\beta(\x, \y; t)-H_{2\pi}(\x, \y; t)=\frac{1}{8\pi\beta
it}\exp\left\{\f{-(r^2+\rho^2)}{4t}\right\}
\int_{\Gamma^-_{\theta-\psi}\cup\Gamma^+_{\theta-\psi}}
\exp\left\{\f{r\rho\cos(\gamma+\psi-\theta)}{2t}\right\}\cot
\left(\frac{\pi\gamma}{\beta}\right)\,d\gamma
\end{equation}
and
\begin{equation}\label{h7}
\int_{C_\beta(R)}\left(H_\beta(\x, \x; t)-\frac{1}{4\pi
t}\right)d\x=\frac{1}{8\pi i t}\int_0^R\,dr\,
r\int_{\Gamma_0^-\cup\Gamma_0^+}\exp\left\{-\frac{r^2\sin^2(\gamma/2)}{t}\right\}
\cot\left(\frac{\pi\gamma}{\beta}\right)\,d\gamma\,.
\end{equation}
 The  integration over $r$ can be done
explicitly  and the right hand side of (\ref{h7}) reduces to
\begin{equation}\label{h8}\frac{1}{16\pi i}
\int_{\Gamma_0^-\cup\Gamma_0^+}\frac{\cot({\pi\gamma}/{\beta})}{\sin^2(\gamma/2)}\,
d\gamma+O(e^{-\epsilon/t})\;
\end{equation}
(one can assume that $\Re\sin^2(\gamma/2)$ is positive and separated from
zero when $\gamma\in\Gamma_0^-\cup\Gamma_0^+$).
 The contour of integration in (\ref{h8}) can be changed for a
negatively oriented circle centered at $\gamma=0$. Since ${\rm Res}\Big|_{\gamma=0}
\frac{\cot(\frac{\pi\gamma}{\beta})}{\sin^2(\gamma/2)}\,
=\frac{2}{3}(\frac{\beta}{2\pi}-\frac{2\pi}{\beta})$,
we arrive at (\ref{as1}).

\begin{remark}
{\rm The Laplacian $\Delta$ corresponding to the flat conical metric
  $(dr)^2+r^2(d\theta)^2,  0\leq\theta\leq
\beta$ on $C_\beta$ with domain $C^{\infty}_{0}(C_\beta\setminus {\cal O})$
has infinitely many self-adjoint extensions. Analyzing the asymptotics of
(\ref{h4}) near the vertex ${\cal O}$, one can show that for any $\y\in
C_\beta$ and $t>0$ the function $H_\beta(\cdot, \y; t)$ belongs to the domain of
the Friedrichs extension $\Delta_F$ of $\Delta$ and does not belong to the
domain of any other extension. Moreover, using Hankel transform, it is
possible to get an explicit spectral representation of $\Delta_F$ (this
operator has absolutely continuous spectrum of infinite multiplicity) and to
show that the Schwartz kernel of the operator $e^{t\Delta_F}$ coincides with
$H_{\beta}(\cdot, \cdot; t)$ (see, e. g., \cite{Taylor} formula (8.8.30)
together with \cite{Carslaw}, p. 370).
 }\end{remark}

\subsubsection{Heat asymptotics for compact polyhedral surfaces}

{\bf Self-adjoint extensions of conical Laplacian.} Let $\L$ be a compact
polyhedral surface with vertices (conical points) $P_1, \dots, P_N$. The
Laplacian $\Delta$ corresponding to the natural flat conical metric on $\L$
with domain $C^{\infty}_0(\L\setminus\{P_1, \dots, P_N\})$ (we remind the
reader that the Riemannian manifold $\L$ is smooth everywhere except the
vertices) is not essentially self-adjoint and  one has to fix one of its
self-adjoint extensions. We are to discuss now the choice of the self-adjoint
extension.

This choice is defined by the prescription of some particular asymptotical
behavior  near the conical points to functions from the domain of the
Laplacian; it is sufficient to consider a surface with  only one conical
point $P$ of the conical angle $\beta$. More precisely, assume that $\L$ is
smooth everywhere except the point $P$ and that some vicinity of $P$ is
isometric to a vicinity of the vertex ${\cal O}$ of the standard cone
$C_{\beta}$ (of course, now the metric on $\L$ can no more  be flat everywhere
in $\L\setminus P$ unless the genus $g$ of $\L$ is greater than one and
$\beta=2\pi(2g-1)$).

For $k\in {\mathbb N}_0$ introduce the functions  $V_\pm^k$ on $C_\beta$ by
$$V^k_\pm(r, \theta)=r^{{\pm}\frac{2\pi k}{\beta}}\exp\left\{i\frac{2\pi
k\theta}{\beta} \right\};\ \ k>0\;,$$
$$V_+^0=1,\ \  V_-^0=\log r\,.$$
Clearly, these functions are formal solutions to the homogeneous problem
$\Delta u=0$ on $C_\beta$. Notice that the functions $V_-^k$ grow near the
vertex but are still square integrable in its vicinity if
$k<\frac{\beta}{2\pi}$.

Let ${\cal D}_{\rm min}$ denote the graph closure of $C^\infty_0(\L\setminus
P)$, i. e.
$$U\in {\cal D}_{{\rm min}} \Leftrightarrow \exists u_m\in C^\infty_0(\L\setminus P),\; W\in L_2(\L):
u_m\rightarrow U \ {\rm and}\ \Delta u_m \rightarrow W\ \ {\rm in}\ \
L_2(\L).$$

Define the space $H_\delta^2(C_\beta)$ as the closure of
$C^\infty_0(C_\beta\setminus{\cal O})$ with respect to the norm
$$||u; H_\delta^2(C_\beta)||^2=\sum_{|\overrightarrow{\alpha} |\leq
2}\int_{C_\beta}r^{2(\delta-2+|
\overrightarrow{\alpha}|)}|D^{\overrightarrow{\alpha}}_xu({\bf
  x})|^2d\x\;.$$
(Here $\overrightarrow{\alpha}$ stands for the multi-index.)

Then for any $\delta\in {\mathbb R}$ such that $\delta-1\neq \frac{2\pi
k}{\beta}, k\in {\mathbb Z}$ one has the a priori estimate
\begin{equation}\label{S1}||u; H^2_\delta(C_\beta)||\leq c||\Delta u;
H^0_\delta(C_\beta)||\end{equation} for any $u\in
C^\infty_0(C_\beta\setminus{\cal O})$ and some constant $c$ being independent
of $u$ (see, e. g., \cite{NazPla}, Chapter 2, Proposition 2.5; here $||u;
H^0_\delta(C_\beta)||^2=\iint_{C_\beta}|u|^2|r|^{2\delta}dx$).

It follows from Sobolev's imbedding theorem (see, e. g.,  \cite{Maz1} or
\cite{Maz2}, eq. (2.30)) that for functions from $u\in H_\delta^2(C_\beta)$
one has the point-wise estimate
\begin{equation}
\label{S2}
r^{\delta-1}|u(r, \theta)|\leq c||v;
H_\delta^2(C_\beta)||\;.
\end{equation}

 Applying estimates (\ref{S1}) and (\ref{S2}),  we see that
functions $u$ from ${\cal D}_{\rm min}$ must obey the asymptotics $u(r,
\theta) =O(r^{1-\delta})$ as $r\to 0$ with any $\delta>0$.

Now the description of the set of all self-adjoint extensions of $\Delta$
looks as follows. Let $\chi$ be a smooth function on $\L$ which is equal to
$1$ near the vertex $P$ and such that in a vicinity of the support of
$\chi$ the Riemann surface
$\L$ is isometric to $C_\beta$. Denote by  ${\frak M}$ the linear subspace of
$L_2(\L)$ spanned by the functions $\chi V_\pm^k$ with $0\leq k<
{\beta}/{2\pi}$. The dimension, $2d$, of ${\frak M}$ is even. To get a
self-adjoint extension of $\Delta$ one chooses a subspace ${\frak N}$ of
${\frak M}$ of dimension $d$ such that
 $$(\Delta u, v)_{L_2(\L)}-(u, \Delta
v)_{L_2(\L)}=\lim_{\epsilon\to 0+}\oint_{r=\epsilon}\left(u\frac{\partial
v}{\partial r}-v\frac{\partial u}{\partial r}\right)=0$$
 for any $u, v\in {\frak N}$. To any such subspace ${\frak N}$
there corresponds a self-adjoint extension $\Delta_{\frak N}$ of $\Delta$
with domain ${\frak N}+{\cal D}_{{\rm min}}$.

The extension corresponding to the subspace ${\frak N}$ spanned by the
functions $\chi V_+^k$, $0\leq k<\frac{\beta}{2\pi}$ coincides with the
Friedrichs extension of $\Delta$. The functions from the domain of the
Friedrichs extension are bounded near the vertex.

>From now on we denote by $\Delta$ the Friedrichs extension of the Laplacian
on the polyhedral surface $\L$; other extensions will not be considered here.

{\bf Heat asymptotics.} The following theorem is the main result of this
section. Its first two statements open a way to define the determinant of the
Laplacian in an arbitrary polyhedral metric on a compact Riemann surface.
\begin{theorem}\label{Mainsp}
Let $\L$ be a compact polyhedral surface with vertices $P_1, \dots, \P_N$ of
conical angles $\beta_1, \dots, \beta_N$. Let $\Delta$ be the Friedrichs
extension of the Laplacian defined on functions from $C^\infty_0(\L\setminus
\{P_1, \dots, P_N\})$. Then

\begin{enumerate}\item The spectrum of the operator $\Delta$ is discrete,
all the eigenvalues of $\Delta$ have finite multiplicity.
\item Let ${\cal H}(\x, \y; t)$ be the heat kernel for $\Delta$. Then for some $\epsilon>0$
\begin{equation}\label{Mainasym}
{\rm Tr}\,e^{t\Delta}=\int_\L {\cal H}(\x, \x; t)\,d\x=\frac{{\rm
Area}(\L)}{4\pi
t}+\frac{1}{12}\sum_{k=1}^N\left\{\frac{2\pi}{\beta_k}-\frac{\beta_k}{2\pi}
\right\}+O(e^{-\epsilon/t}),
\end{equation}
as $t\to 0+$.
\item The counting function, $N(\lambda)$, of the spectrum of $\Delta$ obeys
the asymptotics $N(\lambda)=O(\lambda)$ as $\lambda\to +\infty$.
\end{enumerate}
\end{theorem}
{\it Proof.} 1) The proof of the first statement is a standard exercise (cf.
\cite{King}). We indicate only the main idea.
 Introduce the closure, ${\mathbb H}^1(\L)$,  of the
 $C^\infty_0(\L\setminus\{P_1, \dots, P_N\})$
with respect to the norm $|||u|||=||u; L_2||+||\nabla u; L_2||$. It is
sufficient to prove that any  bounded set $S$ in ${\mathbb H}^1(\L)$ is
precompact in $L_2$-topology (this will imply the compactness of the
self-adjoint operator $(I-\Delta)^{-1}$). Moreover, one can assume that the
supports of functions from $S$ belong to a small ball $B$ centered at a
conical point $P$. Now to prove the precompactness of $S$ it is sufficient to
make use of the expansion with respect to eigenfunctions of the Dirichlet
problem in $B$ and the diagonal process.

2) Let $\L=\cup_{j=0}^N K_j$, where $K_j$, $j=1, \dots, N$ is a neighborhood of
the conical point $P_j$ which is isometric to $C_{\beta_j}(R)$ with some
$R>0$, and $K_0=\L\setminus\cup_{j=1}^N K_j$.

Consider also extended neighborhoods $K^{\epsilon_1}_j\supset K_j$ such that
$K^{\epsilon_1}_j$ is isometric to $C_{\beta_j}(R+\epsilon_1)$ with some
$\epsilon_1>0$ and $j=1, \dots, N$.

Fixing $t>0$ and $\x, \y\in K_j$ with $j>0$, one has (cf. \cite{Cheeger}, p.
578-579)
\begin{equation}
\label{chasti}
\int_0^t\,ds\int_{K_j^{\epsilon_1}}\left(
\psi\{\Delta_{\z}-\partial_s\}\phi-\phi\{\Delta_{\z}+\partial_s\}\psi\right)\,d\z
\end{equation}
\begin{equation*}
=\int_0^t ds\int_{\partial K_j^{\epsilon_1}}\left(\phi\frac{\partial
\psi}{\partial n}-\psi\frac{\partial \phi}{\partial
n}\right)dl(\z)-\int_{K_j^{\epsilon_1}}\left(\phi(\z, t)\psi(\z,t)-\phi(\z,0)\psi(\z, 0)\right)\,d\z
\end{equation*}
 with $\phi(\z,
t)={\cal H}(\z, \y; t)-H_{\beta_j}(\z, \y; t)$ and $\psi(\z,
t)=H_{\beta_j}(\z, \x; t-s)$ (here it is important that we are working with the heat kernel of the
Friedrichs extension of the Laplacian, for other extensions the heat kernel
has growing terms in the asymptotics near the vertex and the right hand side
of (\ref{chasti}) gets extra terms). Therefore,
$$H_{\beta_j}(\x, \y; t)-{\cal H}(\x, \y; t)=\int_0^t ds\int_{\partial
  K_j^{\epsilon_1}}
\left({\cal H}(\y, \z; s)\frac{\partial H_{\beta_j}(\x, \z; t-s)}{\partial n(\z)}
-H_{\beta_j}(\z, \x; t-s)\frac{{\partial \cal H}(\z, \y; s)}{\partial
n(\z)}\right)\,dl(\z)
$$
$$
=O(e^{-\epsilon_2/t})$$ with some $\epsilon_2>0$ as
$t\to 0+$ uniformly with respect to $\x, \y\in K_j$. This
implies the asymptotics
\begin{equation}
\label{assy1}
\int_{K_j}{\cal H}(\x, \x; t)d\x=\int_{K_j} H_{\beta_j}(\x, \x;
t)d\x+O(e^{-\epsilon_2/t})\;,\hskip0.5cm {\rm as}\hskip0.5cm t\to 0^+
\end{equation}
Since the metric on $\L$ is flat in
a vicinity of  $K_0$, one has the asymptotics
$$\int_{K_0}{\cal H}(\x, \x; t)d\x=\frac{{\rm Area}(K_0)}{4\pi t}+O(e^{-\epsilon_3/t})$$
with some $\epsilon_3>0$ (cf. \cite{McKeanSinger}). Now (\ref{Mainasym})
follows from (\ref{as1}).

3) The third statement of the theorem follows from the second one due to the
standard Tauberian arguments.

\subsection{Determinant of Laplacian}

According to Theorem \ref{Mainsp} one can define the determinant, ${\rm
  det}\,\Delta$,
of the Laplacian on a compact polyhedral surface via the
standard Ray-Singer regularization. Namely, introduce the operator
$\zeta$-function
\begin{equation}
\label{oper}
\zeta_{\Delta}(s)=\sum_{\l_k>0}\frac{1}{\lambda_k^s}\;,
\end{equation}
where the summation goes over all strictly positive eigenvalues $\l_k$ of the
operator $-\Delta$ (counting  multiplicities). Due to the third statement of
Theorem \ref{Mainsp}, the function $\zeta_\Delta$ is holomorphic in the
half-plane $\{\Re s>1\}$. Moreover, due to the equality
\begin{equation}
\label{zeta}\zeta_\Delta(s)=\frac{1}{\Gamma(s)}\int_0^\infty\left\{{\rm
Tr}\,e^{t\Delta }-1\right\}t^{s-1}\,dt
\end{equation} and asymptotics
(\ref{Mainasym}), one has the equality
\begin{equation}
\label{predst}
\zeta_\Delta(s)=\frac{1}{\Gamma(s)}\left\{ \frac{{\rm Area}\,(\L)}{4\pi(s-1)}+\left[\frac{1}{12}\sum_{k=1}^N\left\{\frac{2\pi}{\beta_k}-\frac{\beta_k}{2\pi}
\right\}-1\right]\frac{1}{s}+e(s)
 \right\},
\end{equation}
where $e(s)$ is an entire function. Thus, $\zeta_\Delta$ is regular at $s=0$
and one can define the $\zeta$-regularized determinant of the Laplacian (cf.
\cite{Ray}) by
 \begin{equation}
\label{defdet}
{\rm det}\Delta:=\exp\{-\zeta'_\Delta(0)\}\,.
\end{equation}
 Moreover,
(\ref{predst}) and the relation $\sum_{k=1}^N b_k=2g-2$;
$b_k=\frac{\beta_k}{2\pi}-1$ yield
\begin{equation}
\label{euler}
\zeta_\Delta(0)=\frac{1}{12}\sum_{k=1}^N\left\{\frac{2\pi}{\beta_k}-\frac{\beta_k}{2\pi}
\right\}-1=\left(\frac{\chi(\L)}{6}-1\right)+\frac{1}{12}\sum_{k=1}^N\left\{\frac{2\pi}{\beta_k}+\frac{\beta_k}{2\pi}
-2\right\}\; ,
\end{equation}
where $\chi(\L)=2-2g$ is the Euler characteristics of $\L$.

It should be noted that the term $\frac{\chi(\L)}{6}-1$ at the right hand
side of (\ref{euler}) coincides with the value at zero of the operator
$\zeta$-function of the Laplacian corresponding to an arbitrary {\it smooth}
metric on $\L$ (see, e. g., \cite{Sarnak}, p. 155).

Let $\met$ and $\kappa\met$, $\kappa>0$ be two homothetic flat
metrics with the same conical points with conical angles $\beta_1, \dots, \beta_N$. Then  (\ref{oper}), (\ref{defdet}) and (\ref{euler}) imply the
following {\it rescaling property} of the conical Laplacian:
\begin{equation}
\label{rescaling} {\rm log}  \f{   {\rm det}\Delta^{{\kappa\met}}  }{  {\rm
det}\,\Delta^{{\met}}   } =\left\{ -\left(  \frac{\chi(\L)}{6}-1  \right)  -
\frac{1}{12}  \sum_{k=1}^N  \left(
\frac{2\pi}{\beta_k}+\frac{\beta_k}{2\pi}-2  \right) \right\} {\rm
log}\kappa\,.
\end{equation}

\subsection {Variation of the resolvent kernel}

For a pair $(\L, w)$ from
${\cal H}_{g}(k_1,\dots, k_M)$ introduce the Laplacian
$\Delta:=\Delta^{|w|^2}$ in flat conical metric $|w|^2$ on $\L$ (recall that
we always deal with the Friedrichs extensions). The corresponding resolvent
kernel $G(P, Q; \l)$, $\l\in {\mathbb C}\setminus {\rm sp}\,(\Delta)$
\begin{itemize}
\item satisfies $(\Delta_P-\l)G(P, Q; \l)=(\Delta_Q-\l)G(P, Q; \l)=0$ outside the diagonal
$\{P=Q\}$,
\item is bounded near the conical points i. e. for any $P\in \L\setminus
\{P_1, \dots, P_M\}$
$$G(P, Q; \l)=O(1)$$
as $Q\to P_k$, $k=1, \dots, M$,
\item obeys the asymptotics $$G(P,Q;
  \l)=\frac{1}{2\pi}\log|x(P)-x(Q)|+O(1)$$ as $P\to Q$,
where $x(\cdot)$ is
an arbitrary (holomorphic) local parameter near $P$.
\end{itemize}
The following proposition is an analog of the classical Hadamard formula for
the variation of the Green function of the Dirichlet problem in a plane
domain.

\begin{proposition}
\la{propvarres}
 There are the following variational formulas for the resolvent
kernel $G(P, Q; \l)$:
\begin{equation}\label{Had1}
\frac{\partial G(P, Q; \l)}{\partial A_\alpha}\Big|_{z(P),\,z(Q)}=2i\oint_{b_\alpha}\omega(P, Q;
\lambda)\,,\end{equation}

\begin{equation}\label{Had2}
\frac{\partial G(P, Q; \l)}{\partial
B_\alpha}\Big|_{z(P),\,z(Q)}=-2i\oint_{a_\alpha}\omega(P, Q;
\lambda)\,,\end{equation} where \be \omega(P, Q; \l)=G(P, z, \bar{z};
\l)G_{z\bar{z}}(Q, z, \bar{z}; \l) d\bar{z}+G_z(P, z, \bar{z}; \l)G_z(Q, z,
\bar{z}; \l)dz \la{omegaPQl} \ee is a closed $1$-form and $\alpha=1, \dots,
g$;
\begin{equation}\label{Had3}
\frac{\partial G(P, Q; \l)}{\partial
z_m}\Big|_{z(P),\,z(Q)}=-2i\lim_{\epsilon\to 0}\oint_{|z-z_m|=\epsilon}G_z(z,
\bar{z}, P; \l)G_z(z, \bar{z}, Q; \l)dz\,,\end{equation} where $m=2, \dots,
M$ and the circle of integration is positively oriented. It is assumed that
the coordinates $z(P)$ and $z(Q)$ are kept constant under variation of the
moduli $A_\alpha, B_\alpha, z_m$.
\end{proposition}
\begin{remark}{\rm One can unite the formulas (\ref{Had1}-\ref{Had3}) in a
single formula:
\begin{equation}\label{Had}
\frac{\partial G(P, Q; \l)}{\partial \ko_k}\Big|_{z(P),\,z(Q)}=-2i\left\{\oint_{\cy_k}\frac{G(R,
P; \l)\partial_R {\partial_{\bar R}} G(R,Q; \l)+\partial_RG(R, P; \l)\partial_R
G(R,Q; \l)}{w(R)}\right\}\,,
\end{equation}
where $k=1, \dots, 2g+M-1$. }
\end{remark}

{\it Proof of Proposition \ref{propvarres}.}
 We start with the following integral representation of a solution $u$
 to the homogeneous equation $\Delta u-\l u=0$ inside the fundamental polygon $\Lhat$:
\begin{equation}\label{green}
u(\xi, \bar \xi)=-2i\int_{\partial \Lhat}G(z, \bar z, \xi, \bar\xi;
\l)u_{\bar z}(z, \bar z)d\bar z+G_z(z, \bar z, \xi, \bar \xi; \l)u(z, \bar
z)dz\,.
\end{equation}
(We remind the reader that to get (\ref{green}) on has to rewrite the left
hand side of the equality
$$\iint_{\hat \L\setminus B_\epsilon(P)}(\Delta_Q-\lambda)G(P, Q; \l)u(Q)|dz(Q)|^2-\iint_{\hat L\setminus B_\epsilon(P)} G(P, Q;
\lambda)(\Delta_Q-\l)u(Q)|dz(Q)|^2=0$$ as an integral over the boundary
$\partial \hat\L\cup \partial (B_\epsilon(P))$ via the Stokes theorem (here
$B_\epsilon(P)$ is the disk of radius $\epsilon$ centered at $P$) and then
send $\epsilon$ to $0$.)

Let us first prove (\ref{Had2}). Cutting the surface $\L$ along the basic
cycles, we notice that the function $\partial_{B_\a} G(P,\  \cdot\ ; \l)$ is
a solution to the homogeneous equation $\Delta u-\l u=0$ inside the
fundamental polygon (the singularity of $G(P, Q; \l)$ at $Q=P$ disappears
after differentiation) and that the functions $\partial_{B_\a}{G}(P,\ \cdot\
; \l)$ and $\partial_{B_\a}{G}_{\bar{z}}(P,\  \cdot\ ; \l)$ have the jumps
$G_z(P,\ \cdot\ ; \l)$ and $G_{z\bar{z}}(P, \ \cdot\ ; \l)$ on the cycle
$a_\a$, respectively. (This follows from differentiation of the periodicity
relation $G(z+B_\alpha; \bar z+\bar B_\alpha; \lambda; \{A_\alpha, B_\alpha,
z_m\})=G(z, \bar z; \lambda; \{A_\alpha, B_\alpha, z_m\})$ with respect to
$B_\alpha$ and $\bar z$; cf.  the proof of Theorem \ref{varfow}, eq.
(\ref{fzBfz}).)

Applying  the formula (\ref{green}) with $u=\partial_{B_\a}{G}(P,\  \cdot\ ;
\l)$, we get the variational formula (\ref{Had2}). Formula (\ref{Had1}) can
be proved in the same manner.

The closedness of the form (\ref{omegaPQl}), $d\omega(P, Q;
\lambda)=0$, immediately follows from the equation for the resolvent kernel
$G_{z\bar z}(z, \bar z, P; \l)=\frac{\l}{4} G(z, \bar z, P; \l)$.

Let us prove (\ref{Had3}). From now on we assume for simplicity that $k_m=1$,
where $k_m$ is the multiplicity of the zero $P_m$ of the holomorphic
differential $w$ (the case $k_m>1$ differs only by a few details).

Applying  Green formula (\ref{green}) to the domain
$\Lhat\setminus\{|z-z_m|<\epsilon\}$ and $u={\partial
G}/{\partial z_m}$, one gets
\begin{equation}\label{tvetv}
\p_{z_m}{G}(P, Q; \l)=2i\lim_{\epsilon\to 0}\oint_{|z-z_m|=\epsilon}\p_{z_m}\{{G}_{\bar
z}(z, \bar z, Q; \l)\}G(z, \bar z, P; \l)d\bar{z}+\p_{z_m}\{{G}(z,
\bar z, Q;\l)\}G_z(z, \bar z, P; \l)dz\,.
\end{equation}
(Here the circle of integration is positively oriented.) Observe that the
function $x_m \mapsto G(x_m, \bar x_m, P; \l)$ (defined in a small
neighborhood of the point $x_m=0$) is a bounded solution to the elliptic
equation
$$\frac{\partial^2 G(x_m, \bar x_m, P; \l)}{\partial x_m\partial \bar{x}_m}
-\lambda|x_m|^2G(x_m, \bar x_m, P; \l)=0$$ with real analytic coefficients
and, therefore, is real analytic near $x_m=0$.

Recall that $x_m=\sqrt{z-z_m}$. Differentiating the
expansion
\begin{equation}\label{expan}
G(x_m, \bar{x}_m, P; \l)=a_0(P, \l)+a_1(P, \l)x_m +a_2(P,
\l)\bar{x}_m+a_3(P, \l)x_m\bar{x}_m
+\dots
\end{equation}
with respect to $z_m$, $z$ and $\bar z$, one gets the asymptotics
\begin{equation}
\label{assy11}
\p_{z_m}{G}(z, \bar z, Q; \l)=-\frac{a_1(Q,\l)}{2x_m}+O(1)\;,
\end{equation}
\begin{equation}
\label{assy2}
\p_{z_m}{G}_{\bar z}(z,\bar z, Q; \l)=\frac{\{\p_{z_m}{a}_2\}(Q,
  \l)}{2\bar{x}_m}-\frac{a_3(Q, \l)}{4x_m
\bar{x}_m}+O(1)\;,
\end{equation}
\begin{equation}
\label{assy3}
G_z(z, \bar z, P;\l)=\frac{a_1(P, \l)}{2x_m}+O(1)\;,
\end{equation}
Substituting (\ref{assy11}),
(\ref{assy2}) and (\ref{assy3}) into (\ref{tvetv}), we get the relation
$$\p_{z_m}{G}(P, Q, \l)=2\pi a_1(P, \l)a_1(Q, \l)\;.$$
On the other hand, calculation of the right hand side of formula (\ref{Had3})
via (\ref{assy3}) leads to the same result. $\square$

\subsection{Variation of the determinant of the Laplacian}

Introduce the notation
\be Q(\L, |w|^2):= \Big\{\frac{{\rm  {det}} \,
\Delta^{|w|^2}}{{\rm Vol}(\L, |w|^2)\,{\rm det}\Im \B} \Big\}\;,
\la{TLW0}
\ee
 where
${\rm Area}(\L,|w|^2)$ the area of the Riemann surface $\L$ in the metric
$|w|^2$ ($Q$ depends also on the choice of canonical basis of cycles on
$\L$).

The rest of this section is devoted to the proof of the following theorem.
\begin{theorem}\label{vartors}
The following variational formulas hold
\begin{equation}
\label{cor1}
\frac{\partial \log Q(\L, |w|^2)}{\partial \ko_k}=-\frac{1}{12\pi
i}\oint_{\cy_k}\frac{S_B-S_w}{w}\;,
\end{equation}
where $k=1,\dots,2g+M-1$;
 $S_B$ is the Bergman projective connection,
$S_w$ is the projective connection given by the Schwarzian derivative
$\Big\{\int^Pw, x(P) \Big\}$; $S_B-S_w$ is the meromorphic quadratic
differential with poles of the second order at the zeroes $P_m$ of $w$.
\end{theorem}

{\it Proof.} The following proof is based on the ideas of J. Fay  applied in
the context of flat metrics with conical singularities (cf. the proof of
Theorem 3.7 in \cite{Fay92}). In this case the calculations get shorter and
more elementary (in particular, the Ahlfors-Teichm\"uller theory is not used
here).

Due to Theorem \ref{Mainsp} one has
\begin{equation}\label{povtor}
{\rm Tr}\,e^{t\Delta}=\frac{c_0}{t}+c_1+O(t^{N})
\end{equation}
as $t\to 0+$, where $N$ is an arbitrary positive real number, $c_0=\frac{A}{4\pi}$,
and $$A:={\rm Area} (\L, |w|^2)=-\frac{1}{2i}\sum_{\a=1}^g(A_\a\bar
B_\a-\bar A_\a B_\a)$$ is the area of the surface $\L$. The coefficient
$c_1$ is independent of all moduli (we notice also that the coefficient $c_0$ is
independent of the moduli $z_2, \dots, z_M$).

Following \cite{Fay92}, consider the expression
$$J(\l, s)=\frac{1}{s\Gamma(s)}\int_0^{+\infty}e^{-\l t}t^{s-1}h(t)\,dt,$$
where
$$h(t)={\rm Tr}\,e^{t\Delta}-(1-e^{-t^2})-\frac{e^{-t}}{t}[(1+t)c_0+tc_1]\,.$$
Notice that $h(t)=O(t^{-N})$ as $t\to +\infty$ with any $N>0$ and
(\ref{povtor}) implies that $h(t)=O(t)$ as $t\to 0+$. Thus,
$$\frac{d}{d\l}J(\l, s)|_{s=0}=-\int_0^{+\infty}e^{-\l
t}h(t)\,dt=O(\frac{1}{\l^2})$$ as $\l\to +\infty$. From the calculations on
p. 42 of \cite{Fay92} it follows that
$$J(\l, s)=\frac{d}{ds}\zeta_\Delta(s; \l)|_{s=0}+\frac{\gamma}{2}-\int_0^\l\int_0^{+\infty}e^{-t^2-\l t}dt\,d\l+c_0(1+\l-\l\log(\l+1))+c_1\log(1+\l)+O(s),$$
as $s\to 0$, where $\gamma$ is the Euler constant and
$$\zeta_\Delta(s; \l)=\sum_{\l_n\in\, {\rm sp}\,
\Delta\setminus\{0\}}\frac{1}{(\l-\l_n)^s}\,.$$ This implies the relation
$$-\int_0^{+\infty}\frac{d}{d\l}J(\l, s)|_{s=0}\,d\l=J(0, 0)=\zeta_{\Delta}'(0)+\frac{\gamma}{2}+c_0\,$$
and, therefore, one has
\begin{equation}\label{Mainprep}
-\zeta_{\Delta}'(0)=\frac{\gamma}{2}+c_0-\int_0^{+\infty}d\l\int_0^{+\infty}e^{-\l
t}\left[{\rm Tr}\, e^{t\Delta}-(1-e^{-t^2})-\frac{e^{-t}}{t}((1+t)c_0+tc_1)
\right]\,dt\,.
\end{equation}
Consider the variation of (\ref{Mainprep}) with respect to $A_\alpha$.

We need the following Lemma.

\begin{lemma}\label{le123}

The following relation holds
\begin{equation} \p_{A_\a}\left[\iint_\L
F(P)dA(P)\right]=\iint_\L\p_{A_\a}\{F\}(P)dA(P)+\frac{i}{2}\oint_{b_\a}F(z,
\bar z){d\bar{z}},
\la{difintA}
\end{equation} where $dA(P)$ is the area element defined by
the metric $|w|^2$. The formula for differentiation with respect to
$B_\a$ looks similar; the only change is the sign in front of the
contour integral over $a_\a$ in the second term of the right-hand side.
\end{lemma}
{\it Proof.} The function $\L\in P\mapsto z= \int_{P_1}^P w$ is univalent in a
small vicinity $U(Q)$ of any point $Q$ of $\L$ except the zeroes, $P_1,
\dots, P_M$, of the differential $w$. Take a cover of $\L$ by small disks,
$B_m$, centered at the points $P_m$ and the vicinities $U(Q)$, $Q\in \L,
Q\neq P_m$. Let $\{U_j\}$ be a finite subcover and let $\{\chi_j\}$ be the
corresponding (smooth) partition of unity. Cutting $\L$ along the basic
cycles and giving to, say, $A_1$-coordinate a complex increment $\delta A$,
one gets

\begin{equation}\label{eee1}
\delta \iint\chi_j F d{\cal A}=\begin{cases} \iint\chi_j(z, \bar
z)\delta F(z, \bar z)|dz|^2,\ \ \  {\rm if}\ \ \  (w)\cap {\rm supp}\, \chi_j=\emptyset                     \\
4\iint\chi_j(x_m, \bar x_m)\delta F(x_m, \bar x_m)|x_m|^2|dx_m|^2, \ \ \
{\rm if} \ \ \ {\rm supp}\,\chi_j\ni P_m
\end{cases}
\end{equation}
for those $j$ for which the support of $\chi_j$ has no intersection with the
cycle $b_1$.

Let ${\rm supp\,}\chi_j\cap b_1\neq \emptyset$ and let $[0, 1]\ni t\mapsto
\gamma(t)$ be the parameterization of the part of contour $z(b_1)\subset
{\mathbb C}$ inside the support of the function $z\mapsto \chi_j(z, \bar z)$.
After variation of the coordinate $A_1$ this contour shifts to $t\mapsto
\gamma(t)+\delta A_1$. Setting
$$y=\Re z=\Re\gamma(t)+s\frac{\delta A_1+\overline{\delta A_1}}{2};\ \ \ x=\Im
z=\Im \gamma(t)+s\frac{\delta A_1-\overline{\delta A_1}}{2i},$$ with $0\leq s\leq
1$, for $z=x+iy$ in a vicinity of the contour $z(b_1)$ and  using the
relation
$$\frac{\partial(x, y)}{\partial (s, t)}=\Im \gamma'(t)\frac{\delta A_1+\overline{\delta A_1}}{2}-\Re\gamma'(t)\frac{\delta A_1-\overline{\delta A_1}}{2i}(>0!)\,,$$
one finds that
$$\delta \iint \chi_j Fd{A}=\iint\chi_j(z, \bar z)\delta F(z, \bar z)|dz|^2$$
$$+\int_0^1ds\int_0^1dt\, \chi_j(\gamma(t))F(\gamma(t))
\left(\frac{1}{2}\Im \gamma'(t)-\frac{1}{2i}\Re\gamma'(t)\right)
\delta A_1
$$
\begin{equation}
+\left(\frac{1}{2}\Im \gamma'(t)+\frac{1}{2i}\Re\gamma'(t)\right)\overline{\delta A_1}\,.\label{eee2}\end{equation}
where the second term coincides with
$$\left(\frac{i}{2}\int_{b_1}\chi_j F \,\overline{dz} \right)\delta A_1$$ and summing (\ref{eee1}),(\ref{eee2}) over all
$j$ one gets the lemma.

$\Box$

Using the formulas $\p_{A_\a}{c_1}=0$,
$\p_{A_\a}{c_0}=-{\overline{B}_\alpha}/{8\pi i}$  and Lemma \ref{le123}, we
get

\begin{equation}\label{eqn11}
\p_{A_\a}[-\zeta_\Delta '(0)]=-\frac{\overline{B}_\alpha}{8\pi
i}-\int_0^{+\infty}\,d\l\int_0^{+\infty}\,dt\,e^{-\l t}\left\{ \iint_\L(\p_{A_\a}
{\cal H}(P, P, t)+\frac{\p_{A_\a}{A}}{A^2}(1-e^{-t^2}))dA(P)+\right.
\end{equation}
\begin{equation*}
\left. \frac{i}{2}\oint_{b_\a}\left[{\cal H}(z,z,
t)-\frac{1}{A}(1-e^{-t^2})-\frac{e^{-t}}{4\pi t}(1+t) \right]{d\bar{z}}
\right\}\,.
\end{equation*}
(For brevity from now on we suppress  the antiholomorphic part $\bar{z}$ of
the argument $(z, \bar{z})$.)

Using the standard relation
$$G(x, y; \l)=-\int_0^{+\infty}e^{-\l t}{\cal H}(x, y, t)dt$$
between the resolvent and the heat kernels, we rewrite the right hand side of
(\ref{eqn11}) as
\begin{equation}\label{eqn12}
-\frac{\overline{B}_\alpha}{8\pi i}+\int_0^{+\infty}\,d\l\left\{
\iint_\L\{\p_{A_\alpha}{G}\}(P, P; \l)dA(P)-\frac{\p_{A_\alpha}{A}}{A}
I(\l)-\frac{i}{2}\oint_{b_\alpha}\widehat{G}(z, z; \l){d\bar{z}}\right\},
\end{equation}
where the derivative $\partial_{A_\alpha}G(P, Q; \l)$ is nonsingular at the
diagonal $P=Q$ due to (\ref{Had1});
$$I(\l)=\frac{1}{\l}-e^{\l^2/4}\int_{\l/2}^{+\infty}e^{-t^2}\,dt$$
as in (\cite{Fay92}, (2.34)) and $\widehat{G}(z, z; \l)$ is Fay's modified
resolvent
\begin{equation}\label{modres}
\widehat{G}(z, z; \l)=\int_{0}^{+\infty}e^{-\l
t}\left\{{\cal H}(z, z, t)-\frac{1}{A}(1-e^{-t^2})-\frac{e^{-t}}{4\pi
t}(1+t)\right\}\,dt\,
\end{equation}
(see \cite{Fay92}: the last formula on page 42, formulas (2.34),(2.35)
on page 38 and the first
two lines on page 39; to get (\ref{modres}) one has to make use of the fact that the
metric $|w|^2$ is Euclidean in a vicinity of the cycle $b_\a$ and, therefore,
the coefficients $H_0$ and $H_1$ in Fay's formulas are $1$ and $0$, respectively.) For
future reference notice that according to (\cite{Fay92}, p.38) one has the
relation
\begin{equation}\label{modres1}
\widehat{G}(z_1, z_2; \l)=G(z_1, z_2; \l)+\frac{1}{A}I(\l)-\frac{1}{2\pi}\left[
\log|z_1-z_2|+\gamma+\log\frac{\sqrt{\l+1}}{2}-\frac{1}{2(\l+1)} \right]\,,
\end{equation}
where the right hand side of (\ref{modres1}) is nonsingular at the diagonal
$z_1=z_2$. Now (\ref{Had1}) implies
$$\iint_\L\{\p_{A_\a}{G}\}(P, P; \l)dA(P)=\frac{i}{2}\oint_{b_\a}{d\bar{z}}\iint_{\L}\l G(z, P; \l)G(z, P;
\l)dA(P)$$
$$+2i\iint_\L dA(P)\oint_{b_\alpha}G_z(z, P; \l)G_z(z, P; \l)dz\;.$$
The interior contour integral in the last term  has $\delta$-type singularity
as $P$ approaches to the contour $b_\a$ and using Stokes formula and the
(logarithmic) asymptotics of the resolvent kernel at the diagonal, it is easy
to show that \begin{equation}\label{problema2}\iint_\L
dA(P)\oint_{b_\alpha}G_z(z, P; \l)G_z(z, P;
\l)dz=-\frac{1}{16\pi}\oint_{b_\a}{d\bar{z}}+\oint_{b_\a}dz\; {\rm
p.\,v.}\iint_{\L} G_z(z, P; \l)G_z(z, P; \l)\,dA(P)\,.\end{equation}

Indeed, choosing the same partition of unity as in Lemma \ref{le123}, one
rewrites the left hand side of (\ref{problema2}) as
\begin{equation}\label{hiphop}\frac{i}{2}\sum_k\sum_l\iint_\L \chi_k(z_1, \bar
z_1)\left(\oint_{b_\alpha}\chi_l(z, \bar z)\left(G_z(z, z_1;
\l)\right)^2dz\,\right)dz_1\wedge d\bar z_1\end{equation}

For a pair $(k, l)$ such that the ${\rm supp}\, \chi_k\cap
b_\alpha\neq\emptyset$ and ${\rm supp}\,\chi_k\cap {\rm
supp}\,\chi_l\neq\emptyset$  the corresponding term in (\ref{hiphop}) is
\begin{equation}\label{ch1}
\frac{i}{2}\iint_\L \chi_k(z_1, \bar z_1)\left(\oint_{b_k}\chi_l(z, \bar
z)\left(\frac{1}{16\pi^2}\frac{1}{(z-z_1)^2}+H(z, \bar z, z_1, \bar
z_1)\right)dz\right) dz_1\wedge d\bar z_1,\end{equation} where function $H$
has only the first order singularity at the diagonal. The iterated integral
with $\chi_k\chi_l H$ as integrand admits the change of order of integration,
whereas the remaining part of the right hand side of (\ref{ch1}) can be
rewritten as

$$\frac{i}{32\pi^2} \iint_\L  \chi_k(z_1, \bar
z_1)\partial_{z_1}\oint_{b_\alpha}\frac{\chi_l(z, \bar
z)\,dz}{z-z_1}dz_1\wedge d\bar z_1=\frac{i}{32\pi^2}\int_{\partial
\hat{\L}}\chi_k(z_1, \bar z_1)\oint_{b_\alpha}\frac{\chi_l(z, \bar
z)\,dz}{z-z_1}{d\bar z_1}-$$
\begin{equation}\label{ch2}
\frac{i}{32\pi^2}\iint_\L\left(\partial_{z_1}\chi_k(z_1, \bar
z_1)\right)\oint_{b_\alpha}\frac{\chi_l(z, \bar z)\,dz}{z-z_1}dz_1\wedge
{d\bar z_1}
\end{equation}
Due to Plemelj theorem on the jump of the Cauchy type integral the first
integral in (\ref{ch2}) is equal to
$$-\frac{1}{16\pi}\int_{b_\alpha}\chi_k\chi_ld\bar z\,.$$
Changing the order of integration in (\ref{hiphop}) for the remaining pairs
$(k, l)$ (since for these pairs the integrand in (\ref{hiphop})  is
nonsingular,  one can apply Fubini's theorem) and summing over all $k$ and
$l$ we arrive at (\ref{problema2}) (the second term in (\ref{ch2}) after
summation cancels out: $\sum_k\partial_{z_1}\chi_k=\partial_{z_1}1=0$).

Now from the resolvent identity
\begin{equation}\label{resid}\frac{G(Q, P;
\l)-G(Q, P; \mu)}{\l-\mu}=\iint_\L G(P, R; \l)\,G(Q, R; \mu)\,dA(R)
\end{equation}
it follows that the derivative $\partial_\l G(P, Q; \l)$ is nonsingular at
the diagonal $P=Q$ and
\begin{equation}\label{der1}\iint_\L G(z, P; \l)\,G(z, P;
\l)\,dA(P)=\{\partial_ \l G\}(z, z; \l)\;.
\end{equation}
Moreover, according to Lemma 3.3 from \cite{Fay92}  one has
$$\iint_\L
G_{z'}(z', P; \l)\,G_z(z, P; \l)\,dA(P)=-\frac{1}{16
\pi}\frac{{\bar{z}'-\bar{z}}}{z'-z} + {\rm p.\,v.\,}\iint_\L G_z(z, P; \l)G_z(z,
P; \l)dA(P)+O(z'-z)\,,$$ as $z\to z'$ and the resolvent identity (\ref{resid}) implies the
relation
\begin{equation}\label{der2}
{\rm p.v.}\iint G_z(z, P; \l)G_z(z, P;
\l)\,dA(P)=\frac{\partial}{\partial \l}\left\{G_{z'z}(z', z;
\l)-\frac{1}{4\pi}\frac{1}{(z'-z)^2}+\frac{\l}{16\pi}\frac{{\bar{z}'-\bar{z}}}{z'-z}
\right\}\Big|_{z'=z}\,.
\end{equation}
Thus, (\ref{eqn12}) can be rewritten as
\begin{equation}\label{promezh1}
-\frac{\overline{B}_\alpha}{8\pi
i}+\frac{i}{2}\int_{0}^{+\infty}\,d\l\oint_{b_\a}{d\bar{z}}\left[
\l\{\partial_\l G\}(z, z; \l)-\frac{1}{4\pi}+\widehat{G}(z, z;
\l)-\frac{1}{A}I(\l)\right]+
\end{equation}
\begin{equation*}
2i\int_0^{+\infty}\oint_{b_\a}\,dz\frac{\partial}{\partial
\l}\left\{G_{z'z}(z', z;
\l)-\frac{1}{4\pi}\frac{1}{(z'-z)^2}+\frac{\l}{16\pi}\frac{{\bar{z}'-\bar{z}}}{z'-z}
\right\}\Big|_{z'=z}\,.
\end{equation*}
Using (\ref{modres1}), rewrite the expression in the square brackets as
$$\frac{\partial}{\partial \l}\left(\l\widehat{G}-\frac{1}{4\pi}\frac{\l}{\l+1}-\frac{1}{A}\l
I(\l)\right)\,.$$
To finish our calculation we need several lemmas.

The first one is an analog of Corollary 2.8 from \cite{Fay92}.
\begin{lemma}\label{lem1}
In a vicinity of the cycle $b_\alpha$ the following relation holds
\begin{equation}\label{aslem11}
4\pi G_{z'z}(z', z;
\l)=\frac{1}{(z'-z)^2}-\frac{\l}{4}\frac{\bar{z}'-\bar{z}}{{z}'-{z}}+\alpha(z',z),
\end{equation}
where $\alpha(z, z')$ is $O(|z'-z|)$ as $z'\to z$ and $\l$ belongs to any
closed subinterval of  \ $(0, +\infty)$.
\end{lemma}
To prove the lemma we notice that the metric $|w|^2$ is flat in a vicinity of
a point $P\in b_\a$ and the geodesic local coordinates in this vicinity are
given by the local parameter $z$. Therefore, as it is explained on pp. 38-39
of \cite{Fay92} the asymptotical behavior of $4\pi G_{z'z}(z', z; \l)$
coincides with that of the second derivative with respect to $z'$ and $z$ of
the function
\begin{equation}\label{Fun1}F(z',\bar{z}', z,
  \bar{z})=\log|z'-z|^2+\frac{1}{4}\l|z-z'|^2\log|z'-z|^2
\end{equation}
(one has to put $H_0=1$ and $H_1=0$ in Fay's calculations on p.38 of
\cite{Fay92}). This immediately leads to (\ref{aslem11}).

The next two lemmas are classical (see \cite{Fay92}, p.25 and example 2.4
and the formula (2.18) on p.30).
\begin{lemma}\label{lem2}
There is the following Laurent expansion near the pole $\l=0$ of the
resolvent $G(P, Q; \l)$:
\begin{equation}\label{greenf}
G(P, Q; \l)=-\frac{1}{\l{\rm Area}\,(\L)}+G(P, Q)+O(\l)\,,
\end{equation}
as $\l\to 0$, where $G(z', z)$ is the Green function.
\end{lemma}

\begin{lemma}\label{lem3}
The following relation holds
\begin{equation}\label{svberg}
4\pi G_{\zeta'\zeta}(\zeta',
\zeta)=\frac{1}{(\zeta'-\zeta)^2}+\frac{1}{6}S_B(\zeta)-\pi\sum_{\a,
\b=1}^g(\Im \B)_{\a\b}^{-1}v_\a(\zeta)v_\b(\zeta)+O(\zeta'-\zeta),
\end{equation}
as $\zeta'\to \zeta$, where $G(\cdot, \cdot)$ is the Green function from
(\ref{greenf}), $S_B$ is the Bergman projective connection,
$\{v_\a\}_{\a=1}^g$ is the basis of normalized holomorphic differentials on
$\L$ and $\B$ is the matrix of $b$-periods of $\L$; $\zeta$ is an arbitrary
holomorphic local parameter and the functions $\zeta\mapsto v_\alpha(\zeta)$
are defined via $v_\alpha=v_\alpha(\zeta)d\zeta$.
\end{lemma}
It should be noted that the Green functions depends on the metric on $\L$
whereas its second derivative (\ref{svberg}) is independent of the
(conformal) metric.

The last lemma immediately follows from Rauch variational formula
(\ref{varB}) and the obvious relation $2i\p_{\zeta_k}[\log\,{\rm det}\,\Im\B]={\rm
Tr}\{(\Im \B)^{-1}\p_{\zeta_k}\B\}$\,.
\begin{lemma}\label{lem4}
The following relation holds
\begin{equation}\label{matrper}
\p_{A_\a}[\log{\rm det}\Im \B]=\frac{1}{2i}\sum_{\gamma, \b=1}^g(\Im
\B)_{\gamma\b}^{-1}\oint_{b_\a}\f{v_\b v_\gamma}{w}\,.
\end{equation}
\end{lemma}

Now using the asymptotics $I(\l)=O({\l^{-3}})$ as $\l\to+\infty$ and the
Lemmas (\ref{lem1})-(\ref{lem4}), one can perform the integration with
respect to $\l$ in (\ref{eqn12}). This leads to the relation
$$\p_{A_\a}[-\zeta_\Delta
'(0)]=\frac{1}{12\pi i}\oint_{b_\a}\f{S_B-S_w}{w}+\p_{A_\a}[\log \rm{det}\,\Im
\B]+\p_{A_\a}[log A]\;.$$

The latter relation is equivalent to (\ref{cor1}) for $k=1, ..., g$. The
proof of (\ref{cor1}) in the case $k=g+1, \dots, 2g$ is similar.

Consider now the variation of (\ref{Mainprep}) with respect to $z_m$.
Using the
equality $\p_{z_m}{c}_0=\p_{z_m}{c}_1=0$ and (\ref{Had3}),  we get
\begin{equation}
\label{eqn1}
\p_{z_m}[-\zeta_\Delta '(0)] =-2i\lim_{\epsilon\to
0}\int_{0}^{+\infty}d\l\iint_{\L}dA(P)\oint_{|z-z_m|=\epsilon}G_z(z, P;
\l)G_z(z, P; \l)\,dz\,.\end{equation}
After passing to local parameter
$x_m=\sqrt{z-z_m}$, the latter expression can be rewritten as
\begin{equation}\label{eqn2}-2i\lim_{\epsilon\to
0}\oint_{|x_m|=\sqrt{\epsilon}}\frac{dx_m}{2x_m}\int_0^{+\infty}d\l\iint_\L G_{x_m}(x_m,
P; \l)G_{x_m}(x_m, P; \l)dA(P)\,.
\end{equation}
Lemma 3.3 from \cite{Fay92}
implies the relation
$$
\iint_\L G_{x_m'}(x_m', P; \l)G_{x_m}(x_m, P;
\l)dA(P)
$$
\begin{equation}\label{eqn3}
=-\frac{1}{4\pi}|x_m|^2\frac{{\bar{x}_m'-\bar{x}_m}}{x_m'-x_m}+
\iint_\L G_{x_m}(x_m, P;\l)\,G_{x_m}(x_m, P; \l)\,dA(P)+O(|x_m'-x_m|)\,,
\end{equation}
as $x_m'\to x_m$. Using this relation rewrite the right hand side of (\ref{eqn2})
as
\begin{equation}
\label{eqn4} -2i \lim_{\epsilon \to
0}\oint_{|x_m|=\sqrt{\epsilon}}\frac{dx_m}{2x_m}\int_0^{+\infty}\,d\l\left\{
\iint_\L G_{x_m'}(x_m', P; \l)G_{x_m}(x_m, P; \l)dA(P)+\frac{1}{4\pi}|x_m|^2
\frac{\bar{x}_m'-\bar{x}_m}{x_m'-x_m}\right\}\Big|_{x_m=x_m'}\;.
\end{equation}
As before,
using the resolvent identity,  we rewrite the expression inside the braces as
a derivative with respect to $\l$ and see that the right hand side of
(\ref{eqn1}) equals
 \begin{equation}
\label{eqn5}
- 2i\lim_{\epsilon\to
0}\oint_{|x_m|=\sqrt{\epsilon}}\frac{dx_m}{2x_m}\int_0^{+\infty}d\l\frac{\partial}{\partial
\l}\left\{G_{x_m'\,x_m}(x_m', x_m;
\l)-\frac{1}{4\pi}\frac{1}{(x_m'-x_m)^2}+\frac{\l}{4\pi}|x_m|^2\frac{{\bar{x}_m'-\bar{x}_m}}{x_m'-x_m}
\right\}\Big|_{x_m'=x_m}\,.
\end{equation}
To further rewrite (\ref{eqn5}) we need the following two lemmas:
\begin{lemma}\label{lem111}
The following relation holds
\begin{equation}\label{aslem12}
4\pi G_{x_m'\,x_m}(x_m', x_m;
\l)=\frac{1}{(x_m'-x_m)^2}-\frac{1}{4x_m^2}-\l|x_m|^2\frac{{\bar{x}_m'-\bar{x}_m}}{x_m'-x_m}
+\alpha(x_m',x_m),
\end{equation}
where $\alpha(x_m, x_m')$ is $O(|x_m'-x_m|)$ as $x_m'\to x_m$ and $\l$ belongs to any
closed subinterval of  $(0, +\infty)$.
\end{lemma}
To prove the lemma we notice that  the geodesic local coordinates for the
flat metric $|w|^2$ in a vicinity of the point $P_m$ are given by the local
parameter $z=z_m+x_m^2$. Therefore, as it is explained on pp. 38-39 of
\cite{Fay92} the asymptotical behavior of $4\pi G_{x_m'x_m}(x_m', x_m; \l)$
coincides with that of the second derivative with respect to $x_m'$ and $x_m$
of the function
\begin{equation}
\label{Fun}F(x_m',\bar{x}_m', x_m,
\bar{x}_m)=\log|z'-z|^2+\frac{1}{4}\l|z-z'|^2\log|z'-z|^2\,,
\end{equation}
where $z'=z_m+(x'_m)^2$.

Using the Taylor expansion of $(x_m'-x_m)^2F_{x'_m x_m}(x_m',\bar{x}_m', x_m, \bar{x}_m)$ up to
the terms of the second order, we arrive at (\ref{aslem12}).

Further, one has the  following analog of Lemma \ref{lem4}, which is an
immediate consequence of variational formulas (\ref{varB}) for
$k=2g+1,\dots,2g+M-1$.
\begin{lemma}
The following relation holds
\begin{equation}\label{matrper1}
\f{\p}{\p z_m}[\log{\rm det}\Im \B]=\frac{1}{2i}\sum_{\a, \b=1}^g(\Im
\B)_{\a\b}^{-1}\oint_{s_{2g+m-1}}\frac{v_\a v_\b}{w}\,,\hskip0.7cm m=2,\dots,g
\end{equation}
\end{lemma}
These lemmas together with (\ref{eqn5}) and  formulas (\ref{greenf}) and
(\ref{svberg}) written in the local parameter $x_m$ imply the relation
$$\f{\p}{\p z_m}[-\zeta_\Delta'(0)]=-\frac{1}{12\pi
  i}\oint_{s_{2g+m-1}}\frac{S_B-S_w}{w}
+\f{\p}{\p z_m}[\log\rm{det}\,\Im \B]\,,
$$
where $s_{2g+m-1}$ is a small positively oriented circle around $P_m$. The latter
relation is equivalent to (\ref{cor1}) for $k=2g+m-1,\; m=2, \dots, M$.
$\square$

\subsubsection{Infinitesimal Polyakov type formula for the stratum ${\cal H}_g(1, \dots, 1)$.}

The following corollary of Theorem \ref{vartors} is an analog of classical
Polyakov formula for variation of the determinant of Laplacian under
infinitesimal variation of the smooth metric within a given conformal
class \footnote{This theorem gives an answer to the question posed by P. Zograf}.

\begin{theorem}
\la{infinitPolyakov}
Let $\omega$ be a holomorphic differential on $\L$ with
$M=2g-2$ simple zeros $P_1, \dots, P_M$, let $x_m$ be the corresponding
distinguished local parameter near $P_m$ and let $\phi$ be an arbitrary
holomorphic differential on $\L$. Define the function $x_m\mapsto\phi(x_m)$
via $\phi=\phi(x_m)dx_m$ and set $\phi'(P_m):=\phi'(x_m)|_{x_m=0}$. Then
\begin{equation}\label{poland}
\frac{d}{d\epsilon}\Big|_{\epsilon=0}\log\frac{{\rm
det\,}\Delta^{|\omega+\epsilon\phi|^2}}{{\rm Area}(\L,
|\omega+\epsilon\phi|^2)}=\frac{1}{16}\sum_{m=1}^M\phi'(P_m)
\end{equation}
\end{theorem}

{\it Proof.}  Consider the one-parametric family $\omega+\epsilon \phi$.
First, let us find the variational formulas for the coordinate $A_\alpha,
B_\alpha, z_m$ of the point $(\L, \omega+\epsilon\phi)\in H_{g}(1, \dots,
1)$. Obviously, one has
$$\dot{A_\alpha}:=\frac{d}{d\epsilon}A_\alpha\Big|_{\epsilon=0}=\frac{d}{d\epsilon}\oint_{a_\alpha}(\omega+\epsilon\phi)=\oint_{a_\alpha}\phi; \ \ \dot{B_\alpha}=\oint_{b_\alpha}\phi\,.$$
To find the variations of the coordinates $z_m$ one has to find the
derivative
$$\frac{d}{d\epsilon}\int_{\tilde{P}_1(\epsilon)}^{\tilde{P}_m(\epsilon)}(\omega+\epsilon\phi)\Big|_{\epsilon=0},$$
where  $\tilde{P}_m(\epsilon)$ are the zeroes of the differential
$\omega+\epsilon\phi$ (we have $2g-2$ one-parametric families of the zeroes
parameterized by $\epsilon\in[0,\delta]$ with sufficiently small $\delta>0$.)
One has
$$\int_{\tilde{P}_1(\epsilon)}^{\tilde{P}_m(\epsilon)}(\omega+\epsilon\phi)=\int_{P_1}^{P_m}(\omega+\epsilon\phi)+\int_{\tilde{P}_1}^{P_1}\omega+
\epsilon\int_{\tilde{P}_1}^{P_1}\phi+\int_{P_m}^{\tilde{P}_m}\omega+\epsilon\int_{P_m}^{\tilde{P}_m}\phi=$$
$$=\int_{P_1}^{P_m}(\omega+\epsilon\phi)+O(\epsilon^2)\,$$
since $\omega(P_1)=\omega(P_m)=0$. Therefore,
$$\dot{z}_m=\int_{P_1}^{P_m}\phi\,.$$

It is instructive to check the following property: the tangent vector
$${\bf V}=\sum_{\alpha=1}^g\left(\dot{A}_\alpha\frac{\partial}{\partial
A_\alpha}+\dot{B}_\alpha\frac{\partial}{\partial
B_\alpha}\right)+\sum_{m=2}^{M}\dot{z}_m\frac{\partial}{\partial z_m}$$ to
the space $H_g(1, \dots, 1)$ should annihilate any function on $H_{g}(1,
\dots, 1)$ which depends only on moduli of the underlying Riemann surface
$\L$. It is sufficient to show that ${\bf V}\{{\mathbb B}_{\gamma\delta}\}=0$,
for any entry of the matrix ${\mathbb B}$ of the $b$-periods of the surface
$\L$. Indeed,  formulas (\ref{varB}) imply
$${\bf V}\{{\mathbb B}_{\gamma \delta}\}=
-\sum_{\alpha=1}^g\oint_{a_\alpha}\phi\oint_{b_\alpha}\frac{v_\gamma
v_\delta}{\omega}+\sum_{\alpha=1}^g\oint_{b_\alpha}\phi\oint_{a_\alpha}\frac{v_\gamma
v_\delta}{\omega}+\sum_{m=1}^M\int_{P_1}^{P_m}\phi\oint_{P_m}\frac{v_\gamma
v_\delta}{\omega}=$$
$$=\int_{\partial\left[\widehat{\L}\setminus\cup_{m=1}^MB(P_m)\right]}\left(\int_{P_1}^P\phi\right) \frac{v_\gamma(P)
v_\delta(P)}{\omega(P)}=0 \,,
$$
where $B(P_m)$ are small disks centered at $P_m$. We have used Riemann's
bilinear relations and the equality
$$\oint_{P_m}\left(\int_{P_1}^{P}\phi\right)\frac{v_\gamma(P)v_\delta(P)}{\omega(P)}=\left(\int_{P_1}^{P_m}\phi\right)\oint_{P_m}\frac{v_\gamma
v_\delta}{\omega}\,;$$  the latter equality holds because the differential
$\frac{v_\gamma v_\delta}{\omega}$ has  a simple pole at $P_m$.

 Now Theorem \ref{vartors} implies
$$\frac{d}{d\epsilon}\Big|_{\epsilon=0}\log\frac{{\rm det\,}\Delta^{|\omega+\epsilon\phi|^2}}{{\rm Area}(\L, |\omega+\epsilon\phi|^2)}={\bf V}\left\{\log\frac{{\rm det\,}\Delta^{|\omega+\epsilon\phi|^2}}{{\rm
Area}(\L, |\omega+\epsilon\phi|^2){\rm det\,}\Im{\mathbb B}}\right\}=$$
\begin{equation}\label{E1}-\frac{1}{12\pi i}\left\{
-\sum_{\alpha=1}^g\oint_{a_\alpha}\phi\oint_{b_\alpha}\frac{S_B-S_\omega}{\omega}+
\sum_{\alpha=1}^g\oint_{b_\alpha}\phi\oint_{a_\alpha}\frac{S_B-S_\omega}{\omega}+\sum_{m=1}^M\int_{P_1}^{P_m}\phi\oint_{P_m}\frac{S_B-S_\omega}{\omega}
\right\}\,.\end{equation}

In the distinguished local parameter $x_m$ near $P_m$ one has
$\omega=2x_mdx_m$ and denoting by $F(x_m):=F(P):=\int_{P_1}^{P}\phi$, we obtain
$$\oint_{P_m}\left(\int_{P_1}^{P}\phi\right)\frac{S_B-S_\omega}{\omega}=\oint_{|x_m|=\delta}F(x_m)
\frac{(S_B(x_m)-\{x_m^2, x_m\})(dx_m)^2}{2x_mdx_m}$$
$$=\oint_{|x_m|=\delta}\left(\frac{F(x_m)S_B(x_m)}{2x_m}+\frac{3}{4}\frac{F(x_m)}{x_m^3}\right)dx_m$$
$$=F(0)\oint_{|x_m|=\delta}\frac{S_B(x_m)+\frac{3}{2}x_m^{-2}}{2x_m}dx_m+\frac{3}{4}\pi
iF''(0)$$
$$=\int_{P_1}^{P_m}\phi\oint_{P_m}\frac{S_B-S_\omega}{\omega}+\frac{3}{4}\pi i\phi'(P_m).$$
 Therefore, the right hand side of
(\ref{E1}) can be rewritten as

$$-\frac{1}{12\pi i}\left\{
-\sum_{\alpha=1}^g\oint_{a_\alpha}\phi\oint_{b_\alpha}\frac{S_B-S_\omega}{\omega}+
\sum_{\alpha=1}^g\oint_{b_\alpha}\phi\oint_{a_\alpha}\frac{S_B-S_\omega}{\omega}+
\sum_{m=1}^M\oint_{P_m}\left(\int_{P_1}^{P}\phi\right)\frac{S_B-S_\omega}{\omega}
-\frac{3}{4}\pi i\sum_{m=1}^M\phi'(P_m)\right\}=$$
$$=\frac{1}{16}\sum_{m=1}^M\phi'(P_m)-\frac{1}{12\pi
i}\int_{\partial\left[\widehat{\L}\setminus
\cup_{m=1}^MB(P_m)\right]}\left(\int_{P_1}^P\phi\right)\frac{S_B-S_\omega}{\omega}\;.$$
The contour integral in the last expression vanishes since the integrand is holomorphic
inside the union of integration contours, which implies (\ref{poland}).

$\Box$

It is instructive to check this result choosing $\phi=\omega$. Due to the
rescaling property (\ref{rescaling}) of the determinants in conical metrics
one has:
$${\rm det}\Delta^{|\omega+\epsilon\omega|^2}=
|1+\epsilon|^{2\left\{-(\frac{2-2g}{6}-1)-\frac{1}{12}\sum_{1}^{2g-2}(\frac{2\pi}{4\pi}+\frac{4\pi}{2\pi}-2)\right\}}{\rm
det}\Delta^{|\omega|^2}=|1+\epsilon|^{\frac{g-1}{2}+2}{\rm
det}\Delta^{|\omega|^2}$$ and
\be
\frac{d}{d\epsilon}\Big|_{\epsilon=0}\log\frac{{\rm
det\,}\Delta^{|\omega+\epsilon\omega|^2}}{{\rm Area}(\L,
|\omega+\epsilon\omega|^2)}=\frac{d}{d\epsilon}\Big|_{\epsilon=0}\log|1+\epsilon|^{\frac{g-1}{2}}=\frac{g-1}{4}\,
\la{resca}\ee
On the other hand for $\phi=\omega$ one has $\phi'(P_m)=2$ and
$$\frac{1}{16}\sum_{m=1}^M\phi'(P_m)=\frac{4g-4}{16}=\frac{g-1}{4}\,$$
in agreement with (\ref{resca}).

\subsection{Explicit formulas for $\det\Delta^{|w|^2}$}

The following theorem, which is the main result of the present paper, can be considered as a natural generalization of Ray-Singer formula (\ref{Ray}) to the higher genus case.

\begin{theorem}\label{mmaaiinn}
Let a pair $(\L, w)$ be a point of the space ${\cal H}_g(k_1, \dots, k_M)$.
Then the determinant of the Laplacian $\Delta^{|w|^2}$ acting in the trivial
line bundle over the Riemann surface $\L$ is given by the following
expression
\begin{equation}\label{M1111}
{{\rm  det}}\,\Delta^{|w|^2}=C\;{\rm Area}(\L,|w|^2)\;{{\rm  det}}\Im
\B\;|\tau(\L,w)|^2,
\end{equation}
where ${\rm Area}(\L, |w|^2):=\int_{\L}|w|^2$ is the area of $\L$; $\B$ is the
matrix of b-periods; constant $C$  is independent of a point of
connected
component of ${\cal H}_g(k_1, \dots, k_M)$. Here $\tau (\L, w)$ is the Bergman tau-function on the space  ${\cal H}_g(k_1, \dots, k_M)$ given by (\ref{*}).
\end{theorem}

{\it Proof.} The proof immediately follows from the definition of the Bergman
tau-function and Theorems \ref{tauquadr} and \ref{vartors}. $\Box$
\begin{remark} {\rm It can be easily verified that expression (\ref{M1111}) is
consistent with  rescaling property (\ref{rescaling})\footnote{We thank the
anonymous referee for this remark.}.}
\end{remark}

From (\ref{M1111}) we can deduce the ``integrated'' version of the
infinitesimal formula of Polyakov type given by Theorem
\ref{infinitPolyakov}. For simplicity we consider only the generic case of
differentials with simple zeros.

\begin{corollary}
Let $w$ and $\tilde{w}$ be two holomorphic differentials with simple zeros on
the same Riemann surface $\L$. Assume for convenience that none of zeros of
the differential $w$ coincides with a zero of the differential  $\tilde{w}$.
Then the following formula holds: \be \f{\det\Delta^{|w|^2}}{\det
\Delta^{|\tilde{w}|^2}}= \f{{\rm Area}(\L,|w|^2)}{{\rm
Area}(\L,|\tilde{w}|^2)}\; \prod_{k=1}^{2g-2}
\Big|  \f{\res|_{\Pt_m}\{w^2/\wt\}}
 {\res|_{\P_m}\{\wt^2/w\}}    \Big|^{1/12}\;,
\la{conAl} \ee where $\{P_k\}$ are zeros of
$w$; $\tilde{P}_k$ are zeros of $\tilde{w}$.
\end{corollary}

{\it Proof.}
The formula (\ref{conAl}) follows from the expression (\ref{M1111}) for
the determinant of
laplacian in the metric $|w|^2$ and the link
(\ref{zamena}) between Bergman tau-functions computed at the points
$(\L,w)$ and
$(\L,\tilde{w})$
of the space ${\cal H}_g(1,\dots,1)$.

$\Box$

\begin{remark}{\rm
For an arbitrary hermitian metric $\met$ on $\L$ the expression
\be
Q^{-1}:=\f{{\rm Area}(\L,\met)\;{\rm det}\Im\,\B}{{\rm det}\Delta^{\met}}:=
||1\otimes (\hd_1\wedge\dots\wedge \hd_g)||^2_{\met},
\la{Qm1}\ee
with
 $\{\hd_\a\}_{\a=1, \dots, g}$ being the basis of holomorphic 1-forms on $\L$ normalized
by $\oint_{a_\a} \hd_\b=\delta_{\a\b}$,
defines a  Quillen metric on the determinant line
$$
\lambda ({\Oc}_{\L})={\rm det}H^0(\L,\Oc_{\L})\otimes ({\rm det}
H^1(\L,\Oc_{\L}))^{-1}=
{\rm det}H^0(\L,\Oc_{\L})\otimes {\rm det}
H^0(\L,\O_{\L}^1)\,.$$

The formula  (\ref{mainres0}) shows that if $\met$ is chosen to be
flat singular metric with trivial holonomy given by $|w|^2$, then
corresponding function
$Q(\L,|w|^2)$ defined by  (\ref{TLW0}), (\ref{Qm1})  is {\it the modulus square of a holomorophic
function of moduli} (i.e. coordinates on the space of holomorphic
differentials). This property distinguishes singular flat metrics with
trivial holonomy among other metrics of a given conformal class,
in some sense their properties are nicer than those of the metric of constant negative curvature:
for the  Poincar\'e metric $\met$ the Belavin-Knizhnik theorem implies that the second holomorphic-antiholomorphic derivatives of
 $\log||1\otimes (\hd_1\wedge\dots\wedge\hd_g)||_{\met}$ with respect to (Teichm\"uller) moduli are nontrivial (see \cite{Fay92}).}
\end{remark}

It should be noted that some other generalizations of the Ray-Singer theorem
are already known. There exists an explicit formula for the determinant of
Laplacian in the Arakelov metric (see, e. g., \cite{Fay92}, formulas (1.29),
(4.58) and (5.23); see also  references in \cite{Fay92}). For Arakelov metric
$\me$ the property of holomorphic factorization also fails. Another higher
genus analog of the Ray-Singer formula was obtained by  Zograf,  Takhtajan
and  McIntyre (see \cite{McInt,Zograf} and references therein) for ${\rm det}\Delta$
in the Poincar\'e metric in the context of Schottky spaces; in the context of
Hurwitz spaces an analog of the Ray-Singer formula for the determinant of the
Laplacian in the Poincar\'e metric was found in \cite{IMRN1}.

It should be also noted that the results of the present work can be extended
to the case of arbitrary compact polyhedral surfaces (see \cite{K2007}).

{\bf Acknowledgments} We are grateful to R.Went\-worth, S.Zelditch, P. Zograf
and, especially, A.Zorich for important discussions. We thank anonymous
referee for numerous useful comments and suggestions. The work of DK was
partially supported by Concordia Research Chair grant, NSERC, NATEQ and
Humboldt foundation. AK thanks Max-Planck-Institut f\"ur Mathematik in den
Naturwissenschaften for the support.  We both thank Max-Planck-Institut f\"ur
Mathematik at Bonn where the main part of this work was completed for
hospitality and excellent working conditions.

\end{document}